\documentclass[titlepage,12pt,a4paper,twoside,dvips]{article}
\usepackage{latexsym, amssymb,german,diagrams}
\def \SLn{$\textrm{SL}_{n}(K)$ }
\def \GLn{$\textrm{GL}_{n}(K)$ }
\def \SL2{$\textrm{SL}_{2}(K)$ }
\def \SO{$\textrm{SO}_{2}(K)$ }
\def \GL2{$\textrm{GL}_{2}(K)$ }
\def\INVSL2{$K[V]^{SL_{2}(K)}$}
\def\INVSO2{$K[V]^{SO_{2}(K)}$}
\def\INVGL2{$K[V]^{GL_{2}(K)}$}
\def \CM{Cohen-Macaulay }
\def\Bew{\noindent \textit{Beweis. }}
\def\qed{\hfill $\Box$}
\def\Magma{{\sc Magma }} 
\def\caret{\symbol{94}}

\newtheorem{Lemma}{Lemma}[section]
\newtheorem{Satz}[Lemma]{Satz}
\newtheorem{Def}[Lemma]{Definition}
\newtheorem{Korollar}[Lemma]{Korollar}
\newtheorem{Bemerkung}[Lemma]{Bemerkung}
\newtheorem{SatzDef}[Lemma]{Satz und Definition}
\newtheorem{Prop}[Lemma]{Proposition}
\newtheorem{Hauptsatz}[Lemma]{Hauptsatz}
\newtheorem{Bsp}[Lemma]{Beispiel}
\newtheorem{Zus}[Lemma]{Zusammenfassung}

\begin{document}
\pagenumbering{roman}
\begin{titlepage}
{\noindent \Large \bf Technische Universit"at M"unchen\\[2ex]
Fakult"at f"ur Mathematik}\\[16ex]
{\LARGE \bf \mbox{Invarianten zusammenh"angender Gruppen}
und die Cohen-Macaulay Eigenschaft\\[8ex]}
{\Large Diplomarbeit von Martin Kohls\\[8ex]

\begin{center}
\begin{tabular}{ll}
Aufgabensteller: & Prof. Dr. Gregor Kemper\\
\\
Betreuer: & Prof. Dr. Gregor Kemper\\
\\
Abgabetermin: & 1. April 2005
\end{tabular}
\end{center}}
\newpage
\vspace*{2cm}
\noindent Ich erkl"are hiermit, dass ich die Diplomarbeit selbst"andig und nur mit den angegebenen Hilfsmitteln angefertigt habe.\\

\vspace*{2cm}
\noindent Garching, den 1. April 2005
\end{titlepage}
\tableofcontents

\newpage
\pagenumbering{arabic}
\section{Einleitung}
\begin{sloppypar}
Das wichtigste Ziel der Invariantentheorie ist es, den Invariantenring einer linearen algebraischen Gruppe, die linear auf einem Polynomring operiert, zu kennen. Dies ist nicht immer m"oglich, aber es gibt verschiedene Ans"atze, den Invariantenring zumindest bis zu einem bestimmten Grad zu bestimmen. Mit den Al\-go\-rithmen von  Bayer \cite{Bayer} oder  Kemper \cite{Kem1} kann man u.a. eine Vektorraumbasis des Invariantenrings zu gegebenem Grad berechnen. Die Berechnung des In\-var\-ian\-ten\-rings vereinfacht sich jedenfalls erheblich, wenn dieser Cohen-Macaulay ist. Dann ist er n"amlich frei "uber der von einem homogenen Pa\-rame\-ter\-system erzeugten Unteralgebra. W"ahlt man nun eine homogene Basis dieses freien Moduls, so erh"alt man eine Vektorraumbasis des Invarianten\-rings zu gegebenem Grad $d$ als die Menge aller Produkte aus einem homogenen Basiselement des Moduls mit Monomen in den Elementen des Para\-me\-ter\-systems, so dass der entstehende (Total-)grad gleich $d$ ist. In diesem Sinne ist dann also der Invariantenring vollst"andig bekannt (wobei die Bestimmung eines Parametersystems aber ebenso ein schwieriges Problem ist). Nach dem Satz von Hochster und Roberts \cite{HoRo} ist der Invariantenring einer linear reduktiven Gruppe stets Cohen-Macaulay. Dagegen bemerken Hochster und Eagon \cite{HoEa}, dass ihnen kein Beispiel eines nicht Cohen-Macaulay Invariantenrings einer zusammenh"angenden Gruppe bekannt ist:\\

\textit{We know of no example in which $G$ is reductive or connected and semireductive [...] and $R^{G}$ is not Cohen-Macaulay.}\\

Kemper \cite{Kem2} hat jedoch gezeigt, dass zu jeder reduktiven, aber nicht linear reduktiven Gruppe ein solcher \emph{existiert} - also eine gewisse Umkehrung des Satzes von Hochster und Roberts.\\

In dieser Arbeit werden nun zum ersten Mal "uberhaupt (konkret!) Invariantenringe zusammenh"angender algebraischer Gruppen angegeben, die nicht \CM sind.\\

Genau genommen habe ich in meinem Projekt \cite{Projekt} mit Hilfe des Beweises aus \cite{Kem2} f"ur die (zusammenh"angenden und reduktiven) Gruppen \SLn und \GLn (und ihre reduktiven Zwischengruppen) mit einem algebraisch abgeschlossenen K"orper $K$ in positiver Charakteristik $p$ die ersten konkreten Beispiele mit nicht \CM Invariantenring konstruiert. Der Vollst"andigkeit halber habe ich jedoch alle wesentlichen Resultate des Projekts (einschlie"slich der Beweise) in diese Arbeit mit aufgenommen.
Speziell f"ur $n=2$ ergab sich f"ur die Dimension des konstruierten Moduls $8p-5$, also $11$ f"ur $p=2$ und  $19$ f"ur $p=3$. Der Nachweis, dass die entstehenden Invariantenringe nicht \CM sind geht nach \cite{Kem2} "uber die Angabe eines partiellen homogenen Parametersystems aus drei Elementen, das keine regul"are Sequenz bildet, und daher der \CM  Eigenschaft widerspricht. Ein direktes Anwenden der Methode aus \cite{Kem2} (wie in meinem Projekt \cite{Projekt}) liefert ein partielles homogenenes Parametersystems im Grad 1. In dieser Arbeit dagegen wird versucht, das Parametersystem bzw. der von diesem annullierte Kozyklus in einen h"oheren Grad zu legen, wodurch sich eine Dimensionsersparnis ergibt. Gelungen ist dabei eine Dimensionsreduktion f"ur die Gruppen \SL2 und \GL2 von $11$ auf $10$ f"ur $p=2$, von $19$ auf $13$ f"ur $p=3$ und von $8p-5$ auf $7p-2$ f"ur $p>3$.\\
\end{sloppypar}

\noindent Kommen wir zur Gliederung. Die Arbeit beginnt mit den Grundlagen "uber Darstellungen und einer Einf"uhrung in die Invariantentheorie. Beweise werden nur angegeben, sofern sie elementar f"uhrbar sind. Danach wird das f"ur diese Arbeit geschriebene \Magma Programm \emph{IsNotCohenMacaulay} und der zu\-grun\-deliegende Algorithmus vorgestellt. \Magma wurde mir dabei freund\-lich\-er\-wei\-se vom \emph{Centre de Calcul MEDICIS} zur Verf"ugung gestellt. S"amtliche angegebenen Beispiele k"on\-nen mit diesem verifiziert werden. Falls man Com\-puter-Beweise zul"asst, ist ein \emph{true} dieses Programms gleichwertig mit einem Beweis auf Papier - den\-noch wird ein letzterer hier immer gegeben. Dann wird das Hauptverfahren an\-ge\-ge\-ben, welches aus der Arbeit \cite{Kem2} stammt. Der n"achste Abschnitt enth"alt die Zu\-sam\-men\-fassung aller Ergebnisse meines Projekts (mit Beweisen). Es folgt das Kern\-st"uck - die Konstruktion konkreter Beispiele in kleiner Dimension. Zu jedem Beispiel findet sich auch eine Datei auf Diskette, um es mit \Magma zu testen.

Der erste Anhang fasst die hier konstruierten Beispiele lediglich in Form kurzer Tabellen zusammen. Zwischentexte skizzieren nochmals, wie die ein\-zelnen Beispiele erhalten wurden. Der zweite Anhang beschreibt, wie der Bayer-Algorithmus, welcher zur Berechnung der Invarianten inner\-halb des Programms \emph{IsNotCohenMacaulay} verwendet wird, f"ur die uns interes\-sieren\-den Gruppen beschleunigt werden kann.\\

Damit bleibt nur noch, mich bei meinem Betreuer, Prof. Gregor Kemper,  herzlichst f"ur die hilfsbereite Beantwortung aller meiner Fragen und das interessante Thema zu bedanken.

\newpage
\section{Grundlagen}
Wir stellen im Folgenden die ben"otigten Grundlagen der Invariantentheorie zusammen. Auf die Beweise der zitierten, oft sehr tiefliegenden S"atze k"onnen wir hier - bis auf einige Ausnahmen - nicht eingehen. Der Einfachheit halber vereinbaren wir f"ur die {\bf gesamte Arbeit} folgende\\

{\noindent \bf Standardvoraussetzung. }Im Folgenden bezeichnet $K$ stets einen \emph{algebraisch abgeschlossenen} K"orper mit Charakteristik $p$.

\subsection{Variet"aten und Moduln}
\begin{Def}
Eine \emph{(affine) Variet"at} $V \subseteq K^{n}$ ist die Nullstellenmenge eines (endlichen) Systems von Polynomen $f_{1},...,f_{m} \in K[X_{1},...,X_{n}]$, d.h.
\[
V:=\left\{ x=(x_{1},...,x_{n}) \in K^{n}: \quad f_{i}(x)=0 \quad \forall i=1..m\right\}.
\]
\end{Def}

{\noindent }Da nach dem Hilbertschen Basissatz $K[X_{1},...,X_{n}]$ noethersch ist, sind auch Nullstellenmengen beliebig vieler Polynome oder von Idealen in $K[X_{1},...,X_{n}]$ affine Variet"aten.

\begin{Def}
Ein \emph{Morphismus von affinen Variet"aten} $V \subseteq K^{n}, W \subseteq K^{m}$ ist eine Abbildung $f: V \rightarrow W$, die durch Polynome $f_{1},....,f_{m} \in K[X_{1},...,X_{n}]$ gegeben ist, d.h. f"ur alle $x=(x_{1},...,x_{n}) \in V$ ist
\[
f(x)=(f_{1}(x),...,f_{m}(x)) \in W.
\]
\end{Def}

{\noindent }F"ur die folgende Definition bemerken wir, dass f"ur Variet"aten $V \subseteq K^{n}, W \subseteq K^{m}$, die Nullstellenmengen von  $f_{1},...,f_{r} \in K[X_{1},...,X_{n}]$ bzw. $g_{1},...,g_{s} \in K[Y_{1},...,Y_{m}]$  sind, auch $V \times W \subseteq K^{n+m}$ eine affine Variet"at ist, n"amlich Nullstellenmenge von $f_{1},...,f_{m},g_{1},...,g_{s} \in K[X_{1},...,X_{n},Y_{1},...,Y_{m}]$.

\begin{Def}
Eine \emph{lineare algebraische Gruppe} ist eine Variet"at $G \subseteq K^{r}$ zusammen mit  
Morphismen $\circ: G\times G \rightarrow G$ und $^{-1}: G \rightarrow G$, so dass $(G,\circ)$ zusammen mit der Inversenbildung $^{-1}$ eine Gruppe ist.
\end{Def}

\noindent Von nun an bezeichnen wir mit \emph{$G$ stets eine lineare algebraische Gruppe.}\\

\noindent Die wichtigsten Beispiele f"ur lineare algebraische Gruppen in dieser Arbeit sind
\[
\textrm{SL}_{2}(K)=\left\{ \left(
\begin{array}{cc}
a & b\\
c &  d
\end{array} \right) \in K^{4}:ad-bc=1 \right\}
\]
und
\[
\textrm{GL}_{2}(K)=\left\{ \left(
\begin{array}{cc}
a & b\\
c &  d
\end{array} \right)_{e}:=(a,b,c,d,e) \in K^{5}: (ad-bc)e=1 \right\}.
\]
In \GL2 muss man das Inverse $e$ der Determinante der Matrix hinzunehmen, um die Gruppenelemente algebraisch invertieren zu k"onnen; es ist
\[
\left(
\begin{array}{cc}
a & b\\
c &  d
\end{array} \right)_{e} ^{-1} = 
\left(
\begin{array}{cc}
ed & -eb\\
-ec &  ea
\end{array} \right)_{ad-bc}.
\]
\\
 
\begin{Def}
Ein \emph{$G$-Modul} ist ein endlichdimensionaler $K$-Vektorraum $V$ zusammen mit einer linearen, durch Polynome gegebene Operation $\cdot: G \times V \rightarrow V$; D.h. zu einer Basis $(e_{1},...,e_{n})$ von $V$ gibt es Polynome $a_{ij} \in K[X_{1},...,X_{r}], i,j=1..n$ (wobei $G \subseteq K^{r}$) so dass die Darstellungsmatrix der von $\sigma \in G$ durch Linksmultiplikation auf $V$  induzierten linearen Abbildung bzgl. der Basis $(e_{1},...,e_{n})$ durch die Matrix $A_{\sigma}:=\left(a_{ij}(\sigma)\right)_{i,j=1..n}$ gegeben ist:
\[
\sigma \cdot e_{j} = \sum_{i=1}^{n}a_{ij}(\sigma)e_{i} \quad \forall j=1..n, \sigma \in G.
\]
 Der Morphismus $G \rightarrow K^{n \times n}, \quad \sigma \mapsto A_{\sigma}$ hei"st dann eine \emph{Darstellung} von $G$.

Ein \emph{Untermodul} ist ein Untervektorraum $U$ von $V$, der zus"atzlich $G$-invariant ist, d.h. f"ur alle $u \in U$ und $\sigma \in G$ ist auch $\sigma \cdot u \in U$.
\end{Def}

\noindent Man beachte, dass die Forderung nach polynomialen Eintr"agen in $A_{\sigma}$ sehr stark ist. Beispielsweise ist zwar durch 
\[
{\mathbb C}^{\times} \rightarrow {\mathbb C}^{2 \times 2}, \quad a \mapsto
\left( \begin{array}{cc}
1 & \log |a| \\
0 & 1
\end{array} \right)
\] ein Gruppenhomomorphismus gegeben, aber da die Eintr"age der Matrix sicher keine Polynome in $a$ und $a^{-1}$ sind (als algebraische Gruppe ist 
$
{\mathbb C}^{\times}\cong \left\{ (a,b) \in {\mathbb C}^{2}: ab=1 \right\}
$),
handelt es sich hier nicht um die Darstellung eines ${\mathbb C}^{\times}$-Moduls.\\

\noindent In der Invariantentheorie interessiert man sich nun daf"ur, welche Elemente des Moduls unter der Operation von $G$ fix bleiben. Die Menge dieser Elemente,
\[
V^{G}:=\left\{ v \in V: \sigma \cdot v = v \quad \forall \sigma \in G\right\}
\]
bildet offenbar einen Untermodul von $V$, den \emph{Fixmodul}.\\

\subsection{Standardkonstruktionen mit Moduln}
\noindent Es gibt mehrere Verfahren, aus einem bzw. mehreren gegebenen Moduln einen neuen zu konstruieren. Am h"aufigsten verwendet wird nat"urlich die \emph{direkte Summe}, die genau wie f"ur Vektorr"aume mit zus"atzlich kom\-po\-nen\-ten\-wei\-ser $G$-Operation definiert ist und wohl keiner weiteren Erl"auterung bedarf. Wir besprechen hier noch Homo\-morphismen, den Dual, Faktor\-mo\-duln, Ten\-sor\-pro\-duk\-te und symme\-tri\-sche Potenzen.\\

\subsubsection{Homomorphismen von $G$- Moduln}
F"ur zwei $G$-Moduln $V$ und $W$ bildet die Menge aller $K$-Homomorphismen von $V$ nach $W$,
\[
\textrm{Hom}_{K}(V,W):=\left\{f: V \rightarrow W \textrm{ linear} \right\}
\]
einen $G$-Modul mit der Operation
\[
(\sigma \cdot f)(v):=\sigma \cdot f (\sigma^{-1} \cdot v) \quad \forall v \in V
\]
f"ur alle  $\sigma \in G$ und $f \in \textrm{Hom}_{K}(V,W)$. Bei Identifizierung von $\sigma \in G$ mit den induzierten linearen Abbildungen auf $V$ bzw. $W$ schreibt sich das schlampig als
\[
\sigma \cdot f = \sigma \circ f \circ \sigma^{-1}.
\]
Sind $A_{\sigma}$ bzw. $B_{\sigma}$ die Darstellungsmatrizen von $\sigma \in G$ auf $V$ bzw. $W$ bez"uglich der Basen $(v_{1},\ldots,v_{n})$ bzw. $(w_{1},\ldots,w_{m})$ und hat $f \in \textrm{Hom}_{K}(V,W)$ bez"uglich dieser Basen die Dar\-stellungs\-matrix $X \in K^{m \times n}$, so hat $\sigma \cdot f$ die Dar\-stellungs\-matrix
\begin{equation} \label{HomVW}
\sigma \cdot X:=B_{\sigma}XA_{\sigma^{-1}}.
\end{equation}
\\

\noindent Wir bestimmen den Fixmodul von $\textrm{Hom}_{K}(V,W)$: Es ist
\[
f \in \textrm{Hom}_{K}(V,W)^{G} \Leftrightarrow f \circ \sigma =\sigma \circ f  \forall \sigma \in G \Leftrightarrow f(\sigma v)=\sigma f (v) \quad \forall \sigma \in G, v\in V,
\]
d.h. es ist
\[
\textrm{Hom}_{K}(V,W)^{G}=:\textrm{Hom}_{G}(V,W)
\]
gleich der Menge der \emph{$G$-Homomorphismen} von $V$ nach $W$.\\

\begin{Def}
Zwei $G$-Moduln $V$ und $W$ hei\ss en \emph{isomorph}, in Zeichen $V \cong W$, wenn es einen bijektiven $G$-Homomorphismus $f: V \rightarrow W$ gibt, und $f$ hei\ss t dann ein \emph{Isomorphismus} von $G$-Moduln. Sind $\sigma \mapsto A_{\sigma}$ bzw. $\sigma \mapsto B_{\sigma}$ Darstellungen bzgl. Basen $\mathcal{B}_{V}$ bzw. $\mathcal{B}_{W}$ von $V$ bzw. $W$, so sind $V$ und $W$ offenbar genau dann isomorph, wenn es eine invertierbare Matrix $S$ mit
\[
S A_{\sigma} S^{-1} =B_{\sigma} \quad \forall \sigma \in G
\]
gibt, und $S$ ist dann Darstellungsmatrix eines Isomorphismus $f: V \rightarrow W$ bez"uglich der Basen $\mathcal{B}_{V}$ bzw. $\mathcal{B}_{W}$.
\end{Def}

\subsubsection{Der Dual eines $G$- Moduls}
Insbesondere f"ur $W=K$ mit trivialer (d.h. konstanter) Operation von $G$ hei\ss t
\[
V^{*}:=\textrm{Hom}_{K}(V,K)
\]
der \emph{Dual} von V. Mit der Basis $\{1\}$ von $W=K$ und obigen Bezeichnungen haben wir dann (wegen $B_{\sigma}=(1)$) nach Gleichung (\ref{HomVW})
\[
\sigma \cdot X:=XA_{\sigma^{-1}}
\]
mit $X \in K^{1 \times n}$ einem Zeilenvektor. Schreiben wir durch Transponieren diesen Koordinatenvektor eines Elements aus $V^{*}$ als Spalte, so erhalten wir aus $(XA_{\sigma^{-1}})^{T}=A_{\sigma^{-1}}^{T}X^{T}$           das wichtige

\begin{Lemma} \label{dual}
Hat $\sigma \in G$ die Darstellungsmatrix $A_{\sigma}$  bez"uglich einer gegebenen Basis von $V$, so ist $A_{\sigma^{-1}}^{T}$ die Darstellungsmatrix von $\sigma$ bez"uglich der zu\-ge\-h"origen Dualbasis von $V^{*}$. \hfill $\Box$
\end{Lemma}

\begin{Korollar}
Ein $G$-Modul $V$ ist zu seinem Bidual $V^{**}$ isomorph,
\[
V \cong V^{**}.
\]
\end{Korollar}

\noindent \textit{Beweis.} Sei $\sigma \mapsto A_{\sigma}$ eine Darstellung von $V$. Dann hat $V^{*}$ eine Darstellung $\sigma \mapsto B_{\sigma}:=A_{\sigma ^{-1}}^{T}$ und damit $V^{**}$ eine Darstellung
$\sigma \mapsto B_{\sigma^{-1}}^{T}=(A_{(\sigma^{-1})^{-1}}^{T})^{T}=A_{\sigma}$. Daher haben $V$ und $V^{**}$ zwei gleiche Darstellungen und sind daher isomorph. \hfill $\Box$

\subsubsection{Faktormoduln}
\noindent Neben der direkten Summe sind Faktormoduln die n"achst h"aufige Kon\-struk\-tion. Sei dazu $U$ ein Untermodul eines $G$-Moduls $V$. Der Faktorraum (Quo\-tien\-ten\-vek\-tor\-raum) $V/U$ wird zu einem $G$-Modul, dem \emph{Fak\-tor\-mod\-ul} von $V$ nach $U$ durch
\[
\sigma \cdot (v+U) := (\sigma v) + U \quad \forall \sigma \in G, v \in V.
\]
Da $\sigma (U) \subseteq U$, ist die Operation wohldefiniert. Bei geeigneter Wahl der Basis von $V$ kann man die Darstellung von $V/U$ leicht an der Darstellung von $V$ ablesen:

\begin{Lemma} \label{Faktormodul}
Sei $V$ ein $G$-Modul mit Basis $\mathcal{B}_{V}=(v_{1},...,v_{k},v_{k+1},...,v_{n})$ so, dass $\mathcal{B}_{U}=(v_{1},...,v_{k})$ eine Basis des Untermoduls $U$ von $V$ ist. Ist dann durch
\[
\sigma \mapsto 
\left( \begin{array}{cc}
A_{\sigma} & B_{\sigma} \\
0 & C_{\sigma}
\end{array} \right)  \in K^{n \times n}, \quad \forall \sigma \in G
\]
die Darstellung auf $V$ bzgl. $\mathcal{B}_{V}$ gegeben (wobei $A_{\sigma} \in K^{k \times k}, B_{\sigma} \in k^{k \times (n-k)}, C_{\sigma} \in K^{(n-k) \times (n-k)}$), so ist $\sigma \mapsto A_{\sigma}$ die Darstellung von $U$ bzgl. $\mathcal{B}_{U}$ und $\sigma \mapsto C_{\sigma}$ die Darstellung von $V/U$ bzgl. der Basis $(v_{k+1}+U,...,v_{n}+U)$.
\hfill $\Box$
\end{Lemma}

\subsubsection{Tensorprodukte}
\noindent In bekannter Weise wird f"ur zwei $G$-Moduln $V$ und $W$ das Tensorprodukt der Vektorr"aume $V \otimes_{K} W$ "uber $K$ konstruiert - dieser Vektorraum wird dann mittels
\[
\sigma \cdot (v \otimes w) := (\sigma v) \otimes (\sigma w) \qquad \forall \sigma \in G, v \in V, w \in W
\]
und linearer Fortsetzung zu einem $G$-Modul gemacht. F"ur die Details, ins\-besondere die Wohldefiniertheit, verweisen wir auf Lehrb"ucher der Dar\-stel\-lungs\-theorie. Informationen "uber die Darstellung auf $V \otimes W$ liefert das fol\-gende

\begin{Lemma} \label{Tensor}
Seien $V$ bzw. $W$ zwei $G$-Moduln mit Basen $(v_{1},\ldots,v_{m})$ bzw. $(w_{1},\ldots,w_{n})$ und Darstellungsmatrizen $A_{\sigma}$ bzw. $B_{\sigma}$ bez"uglich dieser Basen f"ur $\sigma \in G$. Schreibt man den Koordinatenvektor eines Elements
\[
u = x_{11} \cdot v_{1} \otimes w_{1} + \ldots + x_{mn} \cdot v_{m} \otimes w_{n} \in V \otimes W \]
mit  $x_{ij} \in K$ bez"uglich der Basis $(v_{1} \otimes w_{1}, \ldots, v_{m} \otimes w_{n})$ von $V \otimes W$ als Matrix $X=\left(x_{ij}\right)_{i=1..m, j=1..n} \in K^{m \times n}$, so hat $\sigma \cdot u$ bez"uglich dieser Basis die Koordinatenmatrix
\begin{equation} \label{TensMat}
\sigma \cdot X := A_{\sigma} X B_{\sigma}^{T} \quad \forall \sigma \in G.
\end{equation}
\end{Lemma}

\noindent \textit{Beweis. } Mit $A_{\sigma} =(a_{ij}) \in K^{m \times m}$  und $B_{\sigma}=(b_{kl}) \in K^{n \times n}$ ist
\begin{eqnarray*}
\sigma \cdot u &=& \sum_{j=1}^{m} \sum_{k=1}^{n} x_{jk} \cdot (\sigma v_{j}) \otimes (\sigma w_{k})\\
&=& \sum_{j=1}^{m} \sum_{k=1}^{n}x_{jk} \cdot \left(\sum_{i=1}^{m}a_{ij} v_{i} \right) \otimes \left( \sum_{l=1}^{n}b_{lk} w_{l} \right)\\
&=& \sum_{i,l} \left(\sum_{j=1}^{m} \sum_{k=1}^{n} a_{ij}x_{jk}b_{lk}\right) v_{i} \otimes w_{l},
\end{eqnarray*}
woran man die behauptete Koordinatenmatrix abliest. \hfill $\Box$\\

\begin{Korollar} \label{HomVW}
F"ur zwei $G$-Moduln $V$ und $W$ gilt die Isomorphie
\[
\textrm{Hom}_{K}(V,W) \cong W \otimes V^{*}.
\]
\end{Korollar}

\Bew Sind die Darstellungen von $V$ bzw. $W$ durch $A_{\sigma}$ bzw. $B_{\sigma}$ gegeben, so wird nach Gleichung (\ref{HomVW}) $\textrm{Hom}_{K}(V,W)$ beschrieben durch 
\[
\sigma \cdot X:=B_{\sigma}XA_{\sigma^{-1}}=B_{\sigma}X(A_{\sigma^{-1}}^{T})^{T}.
\]
Nach den Lemmata \ref{dual} und \ref{Tensor} wird so jedoch auch die Darstellung auf $W \otimes V^{*}$ beschrieben. \qed\\

\noindent Wir kommen nun zu der f"ur diese Arbeit wichtigsten Konstruktion. In Satz \ref{PKoz} wird sie zur Angabe einer Serie von nicht \CM Invarianten\-ring\-en f"uhren. Im folgenden Satz geben wir noch zus"atzlich eine Isomorphie an, mit deren Hilfe sich die Struktur der auf Basis dieses Satzes kon\-struier\-ten Beispiele f"ur $p>2$ wesentlich vereinfachen l"asst, indem der nichttriviale Kozyklus in Grad 2 geschoben wird, und die f"ur $p>3$ zu der Dimensions\-reduk\-tion von $8p-5$ auf $7p-2$ f"uhren wird.\\

\begin{Satz} \label{Hom0}
Sei $U$ Untermodul des $G$-Moduls $V$. Durch
\[
\textrm{Hom}_{K}(V,U)_{0}:=\left\{f \in \textrm{Hom}_{K}(V,U): f|_{U}=0 \right\}
\]
wird ein Untermodul von $\textrm{Hom}_{K}(V,U)$ definiert, und es gilt
\[
\textrm{Hom}_{K}(V,U)_{0} \cong U \otimes (V/U)^{*}
\]
\end{Satz}

\noindent {\textit Beweis.} Wir verwenden die Bezeichnungen von Lemma \ref{Faktormodul}. Eine lineare Abbildung $f \in \textrm{Hom}_{K}(V,U)_{0}$ hat dann bzgl. der Basis $\mathcal{B}_{V}$ eine Dar\-stell\-ungs\-ma\-trix der Form
\[
\left( \begin{array}{cc}
0_{k \times k} & X
\end{array} \right ) \quad \textrm{mit } X \in K^{k \times (n-k)}.
\]
Zu $\sigma \cdot f = \sigma \circ f \circ \sigma^{-1}$ geh"ort daher die Dar\-stell\-ungs\-matrix
\[
A_{\sigma} 
\left( \begin{array}{cc}
0_{k \times k} & X
\end{array} \right )
\left( \begin{array}{cc}
A_{\sigma^{-1}} & B_{\sigma^{-1}} \\
0 & C_{\sigma^{-1}}
\end{array} \right)
=
\left( \begin{array}{cc}
0_{k \times k} & A_{\sigma}XC_{\sigma^{-1}}
\end{array} \right),
\]
d.h. die Operation auf $\textrm{Hom}_{K}(V,U)_{0}$ l"asst sich mit Matrizen beschreiben durch
\[
\sigma \cdot X = A_{\sigma} X (C_{\sigma^{-1}}^{T})^{T}.
\]
Nach den Lemmata \ref{dual}, \ref{Faktormodul} und \ref{Tensor} ist dies eine Beschreibung der Operation auf $U \otimes (V/U)^{*}$. \hfill $\Box$\\

Wir wollen den hier konstruierten Isomorphismus noch basisfrei angeben: F"ur $v \in V$ sei $\rho: V \rightarrow V/U, \quad v \mapsto [v]$ der kanonische Epimorphismus. Durch 
\[
U \times (V/U)^{*} \rightarrow \textrm{Hom}_{K}(V,U)_{0}, \quad (u,\varphi) \mapsto u \cdot \varphi \circ \rho: v \mapsto \varphi([v])u 
\]
ist dann eine wohldefinierte bilineare Abbildung gegeben. Nach der uni\-ver\-sel\-len Eigenschaft des Tensorproduktes gibt es dann genau eine lineare Abbildung $U \otimes (V/U)^{*} \rightarrow \textrm{Hom}_{K}(V,U)_{0}$ mit $u\otimes\varphi \mapsto u \cdot \varphi \circ \rho$, und diese ist der gesuchte Isomorphismus, wie man sich leicht "uberzeugen kann.

Alternativ kann man sich auch "uberlegen, dass durch $\textrm{Hom}_{K}(V/U,U)\rightarrow \textrm{Hom}_{K}(V,U)_{0}, \quad f \mapsto f \circ \rho$ ein Isomorphismus gegeben ist. Mit Korollar \ref{HomVW} angewendet auf $\textrm{Hom}_{K}(V/U,U)$ hat man dann ebenfalls
\[
\textrm{Hom}_{K}(V,U)_{0} \cong \textrm{Hom}_{K}(V/U,U) \cong U \otimes (V/U)^{*}.
\]

\subsubsection*{Das Tensorprodukt von Matrizen}
Wir k"onnen zwar nach Gleichung (\ref{TensMat}) bereits mit Tensorprodukten rechnen, aber wir haben noch keine Darstellungsmatrix des Tensorproduktes. Dazu definieren wir f"ur Matrizen $A=(a_{ij}) \in K^{m \times m}$ und $B=(b_{ij}) \in K^{n \times n}$ ihr \emph{Tensorprodukt} oder \emph{Kronecker-Produkt}
\[
A \otimes B:= \left(
\begin{array}{ccc}
a_{11}B & \dots & a_{1m}B \\
\vdots &&\vdots\\
a_{m1}B & \dots & a_{mm}B
\end{array} \right) \in K^{mn \times mn}.
\]
F"ur $A,C \in K^{m \times m}$ sowie $B,D \in K^{n \times n}$ und $E \in K^{k \times k}$ gelten dann die folgenden Regeln:
\begin{eqnarray} \label{MaTens}
(A \otimes B) (C \otimes D) &=& (AC) \otimes (BD) \nonumber \\
(A \otimes B)^{T} & = & A^{T}\otimes B^{T}\\
(A \otimes B)^{-1} & = & A^{-1} \otimes B^{-1} \quad \quad \textrm{f"ur }  A,B \textrm{ invertierbar} \nonumber\\
(A \otimes B)\otimes E&=& A \otimes (B \otimes E). \nonumber
\end{eqnarray}
Der Beweis des n"achsten Lemmas ergibt sich durch einfaches Nachrechnen.

\begin{Lemma} \label{TensDarst}
Seien $V$ bzw. $W$ zwei $G$-Moduln mit Basen $(v_{1},\ldots,v_{m})$ bzw. $(w_{1},\ldots,w_{n})$ und Darstellungsmatrizen $A_{\sigma}$ bzw. $B_{\sigma}$ bez"uglich dieser Basen f"ur $\sigma \in G$. Dann ist
\[
A \otimes B: \sigma \mapsto A_{\sigma} \otimes B_{\sigma}
\]
die Darstellung von $V \otimes W$ bez"uglich der Basis $(v_{1} \otimes w_{1},\ldots,v_{1} \otimes w_{n},\ldots,v_{m} \otimes w_{1},\ldots,v_{m} \otimes w_{n})$. \hfill $\Box$
\end{Lemma}

\noindent Damit und mit der Gleichungsgruppe (\ref{MaTens}) kann man auch sehen, wie zueinander "ahnliche Darstellungen von $V$ und $W$ eine "ahnliche Darstellung von $V \otimes W$ liefern: F"ur invertierbare Matrizen $S$ und $T$ gilt
\[
(S^{-1}A_{\sigma}S) \otimes (T^{-1}B_{\sigma}T)=(S \otimes T)^{-1}(A_{\sigma} \otimes B_{\sigma})(S \otimes T).
\]
 
\begin{Korollar} \label{UVW}
F"ur $G$-Moduln $U,V,W$ gilt die Isomorphie
\[
\left(U\otimes V\right) \otimes W \cong U\otimes \left(V \otimes W\right)
\]
und ein Isomorphismus ist gegeben durch lineare Fortsetzung von
\[
\left(u\otimes v\right) \otimes w \mapsto u\otimes \left(v \otimes w\right) \quad \forall u\in U,v \in V, w\in W.
\]
\end{Korollar}

\Bew Dies folgt aus der Regel $(A_{\sigma} \otimes B_{\sigma})\otimes C_{\sigma}= A_{\sigma} \otimes (B_{\sigma} \otimes C_{\sigma})$ f"ur Darstellungen $A,B,C$ von $U,V,W$. \qed

\begin{Korollar} \label{TensIsom}
F"ur zwei $G$-Moduln $V$ und $W$ gilt
\[
\left(V \otimes W \right)^{*} \cong V^{*} \otimes W^{*}.
\]
\end{Korollar}

\noindent \textit{Beweis.} Seien durch $A_{\sigma}$ bzw. $B_{\sigma}$ Darstellungen auf $V$ bzw. $W$ gegeben. Dann hat $\left(V \otimes W \right)^{*}$ eine Darstellung, die durch $(A \otimes B)_{\sigma^{-1}}^{T} \stackrel{(\ref{MaTens})}{=}A_{\sigma^{-1}}^{T} \otimes B_{\sigma^{-1}}^{T}$ gegeben ist, und eine solche Darstellung hat nach den Lemmata \ref{dual} und \ref{TensDarst} auch der Modul $V^{*} \otimes W^{*}$.
\hfill $\Box$

\subsubsection{Verheftung an Untermoduln} \label{Verheft}
Manchmal ben"otigt man verschiedene Moduln als Untermoduln in einem gr"o\ss eren Modul. In der Regel wird man den gr"o\ss eren dann einfach als direkte Summe der ben"otigten Untermoduln definieren. Es kann jedoch sein, dass man von den Untermoduln nur einen bestimmten Teil explizit braucht, und den anderen Teil, der vielleicht auch bei den anderen Untermoduln vorkommt, irgendwie gemeinsam nutzen will, um Dimension einzusparen. Daf"ur ist der folgende Satz hilfreich.

\begin{Satz}
Seien $V$ und $W$ $G$-Moduln, die beide einen zu $U$ isomorphen Untermodul besitzen. Dann gibt es einen $G$-Modul $X$, der (isomorphe Bilder von) $V$ und $W$ als Untermoduln enth"alt, von $V$ und $W$ erzeugt wird, und es gilt $V \cap W\cong U$. Sind durch
\[
\sigma \mapsto
\left( \begin{array}{cc}
A_{\sigma} & B_{\sigma}\\
0 & C_{\sigma}
\end{array} \right)
\quad \textrm{ bzw. } \quad
\sigma \mapsto
\left( \begin{array}{cc}
A_{\sigma} & D_{\sigma}\\
0 & E_{\sigma}
\end{array} \right)
\]
Darstellungen auf $V$ bzw. $W$ gegeben, wobei die Basen so gew"ahlt sind, dass sie eine Basis des jeweils zu $U$ isomorphen Untermodul enthalten und so, dass die Darstellungen auf diesen isomorphen Untermoduln gleich sind, so ist
$X$ gegeben durch die Darstellung
\[
\sigma \mapsto
\left( \begin{array}{ccc}
A_{\sigma} & B_{\sigma}&D_{\sigma}\\
0 & C_{\sigma}&0\\
0&0&E_{\sigma}\\
\end{array} \right).
\]
\end{Satz}

\Bew Es gen"ugt zu zeigen, dass die angegebene Abbildung ein Gruppen-Homomorphismus ist, da der zu dieser Darstellung geh"orige Modul dann offenbar alle behaupteten Eigenschaften besitzt. Wir verwenden dazu die Homomorphie-Eigenschaft der Darstellungen von $V$ und $W$, also
\[
\left( \begin{array}{cc}
A_{\sigma} & B_{\sigma}\\
0 & C_{\sigma}
\end{array} \right)
\left( \begin{array}{cc}
A_{\tau} & B_{\tau}\\
0 & C_{\tau}
\end{array} \right)
=
\left( \begin{array}{cc}
A_{\sigma}A_{\tau} & A_{\sigma}B_{\tau}+B_{\sigma}C_{\tau}\\
0 & C_{\sigma} C_{\tau}
\end{array} \right)
=
\left( \begin{array}{cc}
A_{\sigma\tau} & B_{\sigma\tau}\\
0 & C_{\sigma\tau}
\end{array} \right)
\]
und
\[
\left( \begin{array}{cc}
A_{\sigma} & D_{\sigma}\\
0 & E_{\sigma}
\end{array} \right)
\left( \begin{array}{cc}
A_{\tau} & D_{\tau}\\
0 & E_{\tau}
\end{array} \right)
=
\left( \begin{array}{cc}
A_{\sigma}A_{\tau} & A_{\sigma}D_{\tau}+D_{\sigma}E_{\tau}\\
0 & E_{\sigma} E_{\tau}
\end{array} \right)
=\left( \begin{array}{cc}
A_{\sigma\tau} & D_{\sigma\tau}\\
0 & E_{\sigma\tau}
\end{array} \right) 
\]
f"ur alle $\sigma,\tau \in G$. Damit berechnen wir
\[
\begin{array}{cl}
&
\left( \begin{array}{ccc}
A_{\sigma} & B_{\sigma}&D_{\sigma}\\
0 & C_{\sigma}&0\\
0&0&E_{\sigma}\\
\end{array} \right)
\left( \begin{array}{ccc}
A_{\tau} & B_{\tau}&D_{\tau}\\
0 & C_{\tau}&0\\
0&0&E_{\tau}\\
\end{array} \right)\\
=&
\left( \begin{array}{ccc}
A_{\sigma}A_{\tau} & A_{\sigma}B_{\tau}+B_{\sigma}C_{\tau}&A_{\sigma}D_{\tau}+D_{\sigma}E_{\tau}\\
0 & C_{\sigma}C_{\tau}&0\\
0&0&E_{\sigma}E_{\tau}\\
\end{array} \right)\\
\stackrel{\textrm{s.o.}}{=}&
\left( \begin{array}{ccc}
A_{\sigma\tau} & B_{\sigma\tau}&D_{\sigma\tau}\\
0 & C_{\sigma\tau}&0\\
0&0&E_{\sigma\tau}\\
\end{array} \right).
\end{array}
\]
Also liefert die Konstruktion tats"achlich eine Darstellung.\qed\\

Wir wollen noch zwei weitere Konstruktionsm"oglichkeiten f"ur die Ver\-heftung skizzieren, die wir jedoch im Folgenden nicht ben"otigen (da wir mit Ver\-heftungen rechnen wollen, ist f"ur uns vor allem die im Satz angegebene Darstellung von $X$ wichtig, weniger eine abstrakte Konstruktion):

Zum einen kann man $X$ mittels $\tilde{U}:=\left\{(u,0)-(0,u): u\in U\right\} \le V \oplus W$ als Faktorraum $X:= (V \oplus W)/\tilde{U}$ definieren.

Alternativ kann man $X$ als \emph{pushout} oder \emph{Fasersumme} (siehe etwa \cite{Kowalsky}, S. 258) des Diagramms
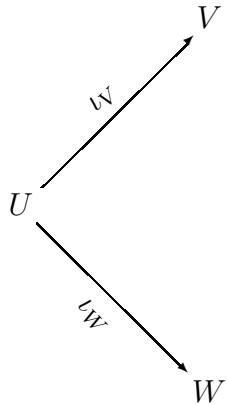
\begin{diagram}
        &               & V\\
        &\ruTo^{\iota_{V}}&\\
U       &               &\\
        &\rdTo_{\iota_{W}}      &\\
        &               &W\\
\end{diagram}
mit den Injektionen $\iota_{V}$ bzw. $\iota_{W}$ von $U$ in $V$ bzw. $W$ definieren. Dann ist der pushout dieses Diagramms ein Modul $X$ zusammen mit Homomorphismen von $V$ bzw. $W$ in $X$ derart, dass das folgende Diagramm kommutativ wird,
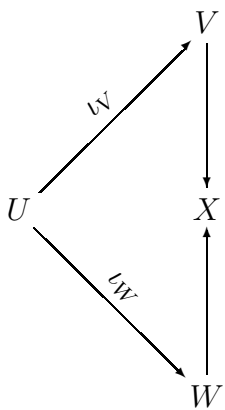
\begin{diagram}
        &               & V\\
        &\ruTo^{\iota_{V}}      & \dTo\\
U       &               &X\\
        &\rdTo^{\iota_{W}}      &\uTo\\
        &               &W\\
\end{diagram}
und es f"ur jeden weiteren Modul $X'$ mit einem solchen kommutativen Dia\-gramm genau einen Homomorphismus von $X$ in $X'$ gibt, der das folgende Diagramm kommutativ macht:
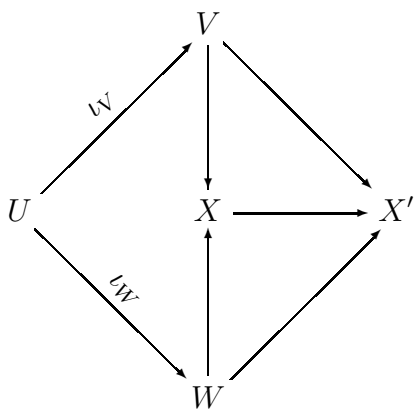
\begin{diagram}
        &                                       & V             &       &\\
        &\ruTo^{\iota_{V}}      & \dTo  &\rdTo&\\
U       &                                       &X              &\rTo&X'\\
        &\rdTo^{\iota_{W}}      &\uTo   &\ruTo  &\\
        &                                       &W              &               &       \\
\end{diagram}
Geht man den Beweis des Satzes "uber die Existenz der Fasersumme durch, sieht man, dass $X$ dann genau als $(V \oplus W)/\tilde{U}$ wie oben definiert wird.

\subsubsection{Symmetrische Potenzen}
Gegeben sei ein $G$-Modul $V$ mit Basis $X=(X_{1},...,X_{n})$. Wir wollen f"ur $k \in \mathbb N_{0}$ einen $G$-Modul $S^{k}(V)$ mit einer Basis $\left(X_{1}^{i_{1}}X_{2}^{i_{2}} \dots X_{n}^{i_{n}}:  i_{1} + \ldots +i_{n}=k\right)$ konstruieren, so dass die Operation auf $S^{k}(V)$ durch
\[
\sigma \cdot X_{1}^{i_{1}}X_{2}^{i_{2}} \dots X_{n}^{i_{n}}:= (\sigma X_{1}^{i_{1}}) (\sigma X_{2}^{i_{2}}) \dots (\sigma X_{n}^{i_{n}})
\]
und formales ausmultiplizieren wie bei Polynomen der rechts auftretenden Linearkombinationen gegeben ist. Wie bei Polynomen "ublich soll dabei die Reihenfolge der Variablen $X_{i}$ keine Rolle spielen. Man fasst also die Elemente aus $V$ als Polynome ersten Grades in $K[X]$ auf und die Elemente aus $S^{k}(V)$ als homogene Polynome $k$-ten Grades in $K[X]$. (Insbesondere soll $S^{0}(V):=K$ mit der trivialen Operation sein). Dabei werden die Basis\-ele\-men\-te $X_{i} \in V$  mit den Variablen eines Polynomringes gleichgesetzt.\\
Die obigen Ziele lassen sich durch folgende Definition realisieren:

\begin{Def}
Die \emph{k-te symmetrische Potenz} $S^{k}(V)$ mit $k \ge 1$ ist der Faktormodul
\[
S^{k}(V):=\left(\bigotimes_{i=1}^{k}V\right)/N
\]
mit dem von allen symmetrischen Differenzen
\[
v_{1}\otimes \ldots \otimes v_{k} -  v_{\pi(1)}\otimes \ldots \otimes v_{\pi(k)}
\quad \textrm{\mit} v_{1},\ldots,v_{n} \in V, \pi \in S_{k}
\]
erzeugten Untermodul $N$.
\end{Def}

\noindent Damit ist die Existenz symmetrischer Potenzen gesichert; im
wesentlichen werden wir mit ihnen aber einfach formal so rechnen, wie
oben beschrieben, wobei wir auf den Beweis, dass obige Definition dies
tats"achlich erlaubt, verzichten wollen.\\

\noindent Wir geben noch ein h"aufig verwendete Isomorphie an: F"ur $G$-Moduln
$U$ und $V$ gilt
\[
S^{2}(U \oplus V) \cong S^{2}(U) \oplus S^{2}(V) \oplus U \otimes V.
\]

\newpage
\subsection{Beispiele f"ur Moduln} \label{BspMod}
Wir geben hier die im Folgenden h"aufig verwendeten Grundmoduln an. Dabei  beschr"anken wir uns auf den Fall, dass $G$ eine (algebraische) Untergruppe von \GL2 ist. Wir geben einen Modul meist in Form $V=\langle  X_{1},\ldots,X_{n}\rangle $ mit einer Darstellungsmatrix $A_{\sigma}$ an, wobei die Notation andeuten soll, dass $A_{\sigma}$ bez"uglich der geordneten Basis $(X_{1},\ldots,X_{n})$ gegeben ist. Ist $W=\langle Y_{1},\ldots,Y_{m}\rangle $ ein weiterer Modul mit Darstellung $B_{\sigma}$, so geben wir die Konstruktionen des letzten Abschnitts oft durch Kombination der Basiselemente an, also bezeichnen wir beispielsweise mit $\langle X_{1},\ldots,X_{n},Y_{1},\ldots,Y_{m}\rangle$ die direkte Summe $V \oplus W$ mit Darstellungsmatrix $\left( \begin{array}{cc} A_{\sigma} & 0\\ 0 & B_{\sigma} \end{array} \right)$.\\

\noindent Analog bezeichnen wir z.B. f"ur $m=n=2$, also $V=\langle X_{1},X_{2} \rangle$ und $W=\langle Y_{1},Y_{2} \rangle$, mit $\langle X_{1} \otimes Y_{1},X_{1} \otimes Y_{2},X_{2} \otimes Y_{1},X_{2} \otimes Y_{2} \rangle$ das Tensorprodukt $V \otimes W$ mit Darstellungsmatrix $A_{\sigma} \otimes B_{\sigma}$. Manchmal kann es aber von Vorteil sein, die Darstellungsmatrix des Tensorprodukts bez"uglich einer anderen Basis zu berechnen. Dann schreiben wir z.B. $V \otimes W = \langle X_{1} \otimes Y_{1}, X_{2} \otimes Y_{2}, X_{1}\otimes Y_{2}-X_{2}\otimes Y_{1},X_{2}\otimes Y_{1} \rangle$, wenn wir andeuten wollen, dass die Darstellungsmatrix des Tensorprodukts bez"uglich der angedeuteten Basis berechnet wird. In diesem Fall w"are also die Darstellungsmatrix, mit der im Folgenden gerechnet werden soll, gegeben durch
\begin{equation} \label{XXYY}
\left(
\begin{array}{cccc}
1 & 0 & 0 & 0\\
0 & 0 & 1 & 0 \\
0 & 0 & -1 & 1 \\
0 & 1 & 0 & 0 \\
\end{array} \right)^{-1}
(A_{\sigma} \otimes B_{\sigma})
\left(
\begin{array}{cccc}
1 & 0 & 0 & 0\\
0 & 0 & 1 & 0 \\
0 & 0 & -1 & 1 \\
0 & 1 & 0 & 0 \\
\end{array} \right).
\end{equation}\\

\noindent Genauso gehen wir bei symmetrischen Potenzen vor. So bezeichnen wir z.B. die zweite symmetrische Potenz von $V=\langle X,Y \rangle$ mit $S^{2}(V)=\langle X^{2},Y^{2},XY \rangle$. Falls es Verwechslungen bzgl. der entstehenden Notation f"ur die Basiselemente geben kann, muss geklammert werden, also z.B.
\[
S^{2}\left(\langle X^{2},Y^2,XY \rangle \right)=\langle (X^{2})^{2},(Y^{2})^{2},(XY)^{2},(X^{2})(Y^{2}),(X^{2})(XY),(Y^{2})(XY) \rangle.
\]

\noindent Wir reservieren nun einige Bezeichnungen von Basis-Elementen, die immer zu den gleichen Darstellungen geh"oren sollen (bis auf seltene Ausnahmen, wo wir das vermerken). Es sei also $G$ eine Untergruppe von \GL2. Mit $\langle X,Y \rangle$ bezeichnen wir ab jetzt stets den $G$-Modul mit der nat"urlichen Darstellung 
$
\sigma \mapsto
A_{\sigma}
=
\left( \begin{array}{cc} a & b\\ c & d \end{array} \right)$ 
f"ur $\sigma=\left( \begin{array}{cc} a & b\\ c & d \end{array} \right)_{e} \in G$, d.h.
\[
\begin{array}{l}
\sigma \cdot X = aX+cY\\
\sigma \cdot Y = bX+dY.
\end{array}
\] 
(Wir geben diese und einige weitere im Folgenden verwendeten Be\-zeichnung\-en f"ur Moduln in Form einer Tabelle am Ende des Abschnitts an).\\

\noindent F"ur den Dual $\langle X,Y \rangle^{*}=:\langle X^{*},Y^{*} \rangle$ erhalten wir aus
\[
\left(
\begin{array}{cc}
a & b\\
c &  d
\end{array} \right)_{e} ^{-1} = 
\left(
\begin{array}{cc}
ed & -eb\\
-ec &  ea
\end{array} \right)_{ad-bc}
\]
und transponieren die Darstellung
\[
\sigma \mapsto
\left(
\begin{array}{cc}
ed & -ec\\
-eb &  ea
\end{array} \right). 
\]
Speziell f"ur \SL2, also $e=1$ erhalten wir so bzgl. der Basis $(-Y^{*},X^{*})$ wieder die Darstellungsmatrix $\left(
\begin{array}{cc}
a & b\\
c &  d
\end{array} \right)$ und k"onnen daher notieren

\begin{Bemerkung}
F"ur die Gruppe \SL2 ist der nat"urliche Modul $\langle X,Y \rangle$ selbstdual,
\[
\langle X,Y \rangle \cong \langle -Y^{*},X^{*} \rangle
\]
\end{Bemerkung}
\qed\\

\noindent Im Fall einer Charakteristik $p>0$ enth"alt die $p$-te symmetrische Potenz $S^{p}\left( \langle X,Y \rangle \right)$ den Untermodul $\langle X^{p},Y^{p} \rangle$ mit der Darstellung $\sigma \mapsto \left(
\begin{array}{cc}
a^{p} & b^{p}\\
c^{p} &  d^{p}
\end{array} \right)$. In der Tat gilt mit dem Frobenius-Homomorphismus
\[
\begin{array}{l}
\sigma \cdot X^{p} = (aX+cY)^{p}=a^{p}X^{p}+c^{p}Y^{p}\\
\sigma \cdot Y^{p} = (bX+dY)^{p}=b^{p}X^{p}+d^{p}Y^{p}.
\end{array}
\] 
Dieser Modul ist ebenfalls selbstdual f"ur \SL2.

\noindent Das Tensorprodukt $\langle X,Y \rangle \otimes \langle X,Y \rangle =: \langle X \otimes X,X \otimes Y,Y \otimes X,Y \otimes Y \rangle$ hat eine Darstellung gegeben durch
\begin{equation} \label{tensXY}
A_{\sigma}=\left( \begin{array}{cc} a & b\\ c & d \end{array} \right) \otimes \left( \begin{array}{cc} a & b\\ c & d \end{array} \right) =
\left( \begin{array}{cccc}
a^{2} & ab & ab & b^{2}\\
ac & ad & bc & bd\\
ac & bc & ad & bd\\
c^{2} & cd & cd & d^{2}
\end{array} \right).
\end{equation}
Nach der letzten Bemerkung und Korollar \ref{TensIsom} ist dieses Tensorprodukt f"ur \SL2 ebenfalls selbstdual.

F"ur sp"ater ben"otigen wir die Darstellung bzgl. der Basis
\[
\mathcal{B}=\left( X \otimes X, Y \otimes Y, X\otimes Y-Y\otimes X,Y\otimes X \right)
\]
gem"a\ss{} Gleichung (\ref{XXYY}) explizit:
\begin{equation} \label{XXYYbasis}
\begin{array}{cl}
&
\left( \begin{array}{cccc}
1 & 0 & 0 & 0\\
0 & 0 & 1 & 0 \\
0 & 0 & -1 & 1 \\
0 & 1 & 0 & 0 
\end{array} \right)^{-1}
\left( \begin{array}{cccc}
a^{2} & ab & ab & b^{2}\\
ac & ad & bc & bd\\
ac & bc & ad & bd\\
c^{2} & cd & cd & d^{2}
\end{array} \right)
\left(
\begin{array}{cccc}
1 & 0 & 0 & 0\\
0 & 0 & 1 & 0 \\
0 & 0 & -1 & 1 \\
0 & 1 & 0 & 0 
\end{array} \right)\\
=&
\left( \begin{array}{cccc}
1 & 0 & 0 & 0\\
0 & 0 & 0 & 1 \\
0 & 1 & 0 & 0 \\
0 & 1 & 1 & 0 
\end{array} \right)
\left( \begin{array}{cccc}
a^{2} & b^{2} & 0 & ab\\
ac & bd & ad-bc & bc\\
ac & bd & -(ad-bc) & ad\\
c^{2} & d^{2} & 0 & cd
\end{array} \right)\\
=&
\left( \begin{array}{cccc}
a^{2} & b^{2} & 0 & ab\\
c^{2} & d^{2} & 0 & cd\\
ac & bd & ad-bc & bc\\
2ac & 2bd & 0 & ad+bc\\
\end{array} \right).
\end{array}
\end{equation}\\
 
\noindent F"ur die zweite symmetrische Potenz $S^{2}\left( \langle X,Y \rangle \right)=\langle X^{2},Y^2,XY \rangle$ erhalten wir aus
\[
\begin{array}{l}
\sigma \cdot X^{2}=(aX+cY)^{2}=a^{2}X^{2}+c^{2}Y^{2}+2acXY\\
\sigma \cdot Y^{2}=(bX+dY)^{2}=b^{2}X^{2}+d^{2}Y^{2}+2bdXY\\
\sigma \cdot XY=(aX+cY)(bX+dY)=abX^{2}+cdY^{2}+(ad+bc)XY\\
\end{array}
\]
die Darstellungsmatrix 
\begin{equation} \label{SXY}
A_{\sigma}=
\left( \begin{array}{ccc}
a^{2} & b^{2} & ab\\
c^{2} & d^{2} & cd\\
2ac & 2bd & ad+bc
\end{array} \right) \quad
\textrm{mit } \sigma = \left( \begin{array}{cc} a & b\\ c & d \end{array} \right)_{e}.
\end{equation}
Man kann die Darstellungsmatrix einer zweiten symmetrischen Potenz auch aus der Ihres Tensorprodukts mit sich selbst erhalten, was an diesem Beispiel demonstriert sei: Wir erhalten die zweite Potenz ja aus dem Tensorprodukt durch Identifizierung von $X^{2}\leftrightarrow X\otimes X$ sowie $Y^{2}\leftrightarrow Y\otimes Y$ und  $XY \leftrightarrow X\otimes Y, \,  Y\otimes X$. Daher erhalten wir durch zusammenfassen (Addieren) der zweiten und dritten Spalte und streichen der zweiten oder dritten Zeile aus (\ref{tensXY}) die Matrix (\ref{SXY}) (allerdings bzgl. der Basis $(X^2,XY,Y^2)$). Genau diese "Uberlegung hat "ubrigens den Basiswechsel motiviert, der zur Darstellung (\ref{XXYYbasis}) gef"uhrt hat. Die Untermatrix von (\ref{XXYYbasis}) die zu den Zeilen/Spalten 1,2,4 geh"ort, liefert daher ebenfalls eine Darstellung der zweiten Potenz als Faktormodul $\langle X,Y \rangle \otimes \langle X,Y \rangle/\langle X\otimes Y- Y \otimes X \rangle$, vgl. auch Lemma \ref{Faktormodul}.\\

\noindent Speziell im Fall der Charakteristik $p=2$ und der Gruppe \SL2 wird (\ref{SXY}) zu
\[
A_{\sigma}=
\left( \begin{array}{ccc}
a^{2} & b^{2} & ab\\
c^{2} & d^{2} & cd\\
0 & 0 & 1
\end{array} \right)
\quad \textrm{mit } \sigma = \left( \begin{array}{cc} a & b\\ c & d \end{array} \right).
\]
Die Darstellung des h"aufig verwendeten Duals, der die extra Bezeichnung
$\langle \mu,\nu,\pi \rangle := \langle X^{2},Y^2,XY \rangle^{*}$ mit
$\mu:=(X^{2})^{*}, \nu=(Y^{2})^{*}, \pi:=(XY)^{*}$ bekommt, erhalten wir dann durch transponieren obiger Matrix und auswerten an der Stelle $\sigma^{-1} = \left( \begin{array}{cc} d & b\\ c & a \end{array} \right)$ zu
\[
A_{\sigma^{-1}}^{T}=
\left( \begin{array}{ccc}
d^{2} & c^{2} & 0\\
b^{2} & a^{2} & 0\\
bd & ac & 1
\end{array} \right).
\]
Damit ist $\pi$ eine Invariante von $\langle \mu,\nu,\pi \rangle$. Da $\langle X^{2},Y^2,XY \rangle$ keine Invariante hat, sieht man, dass die symmetrische Potenz eines selbstdualen Moduls nicht mehr selbstdual sein muss.\\
\newpage
\subsubsection{Tabelle mit Darstellungen verschiedener Moduln}
\begin{table}[here] 
\begin{tabular}{ccccc} 
Bezeichnung & Darstellungsmatrix & $G$ & p &\\
\hline \tabularnewline[0.5ex]
$\langle X,Y \rangle$ & $\left( \begin{array}{cc} a & b\\ c & d \end{array} \right)$ &\GL2& &* \tabularnewline[3ex]

$\langle X^{p},Y^{p} \rangle$ & $\left(
\begin{array}{cc}
a^{p} & b^{p}\\
c^{p} &  d^{p}
\end{array} \right)$ & \GL2  & $p>0$&* \tabularnewline[3ex]

$\langle X \otimes X,X \otimes Y,Y \otimes X,Y \otimes Y \rangle$ &
$\left( \begin{array}{cccc}
a^{2} & ab & ab & b^{2}\\
ac & ad & bc & bd\\
ac & bc & ad & bd\\
c^{2} & cd & cd & d^{2}
\end{array} \right)$ &\GL2 & &*\tabularnewline[6ex]

$\langle X^{2},Y^2,XY \rangle$
&$\left( \begin{array}{ccc}
a^{2} & b^{2} & ab\\
c^{2} & d^{2} & cd\\
2ac & 2bd & ad+bc
\end{array} \right)$
& \GL2 &  \tabularnewline[5ex]

$\langle X^{2},Y^2,XY \rangle$
&$\left( \begin{array}{ccc}
a^{2} & b^{2} & ab\\
c^{2} & d^{2} & cd\\
0 & 0 & 1
\end{array} \right)$
& \SL2 & 2 \tabularnewline[5ex]

$\langle \mu,\nu,\pi \rangle=\langle X^{2},Y^2,XY \rangle^{*}$&
$\left( \begin{array}{ccc}
d^{2} & c^{2} & 0\\
b^{2} & a^{2} & 0\\
bd & ac & 1
\end{array} \right)$&
\SL2 & 2 \tabularnewline[4ex]
\hline
\end{tabular}
\caption{Einige Moduln mit ihren Darstellungen. Ein * in der letzten Spalte bedeutet, dass der Modul f"ur die Gruppe \SL2 selbstdual ist.} \label{Darstellungen}
\end{table} 

\newpage
\subsection{Invariantentheorie}
\noindent In der Invariantentheorie interessiert man sich weniger f"ur die Invarianten eines einzelnen Moduls, sondern m"oglichst gleich f"ur die Invarianten aller symmetrischen Potenzen dieses Moduls (bzw. genauer seines Duals). Diese bilden einen Ring, den Invariantenring. Jedoch kann man dies in ein allge\-mein\-er\-es Konzept einf"ugen.

\subsubsection{Algebren, Polynomringe, Invariantenringe}
\begin{Def}
\begin{enumerate}
\renewcommand{\labelenumi}{(\alph{enumi})}
\item Eine \emph{$K$-Algebra} ist ein $K$-Vektorraum $R$ mit einer zus"atzlichen Multiplikation $\cdot$, so dass $R(+,\cdot)$ zu einem \emph{kommutativen} Ring mit Eins $1$ wird, so dass gilt
\[
\lambda (v \cdot w) = (\lambda v) \cdot w = v \cdot (\lambda w) \quad \forall \lambda \in K, v,w \in R.
\]
\item $R$ hei\ss t \emph{graduiert}, wenn es eine direkte Summenzerlegung in $K$-Vek\-tor\-r"aume
\[
R=R_{0} \oplus R_{1} \oplus R_{2} \oplus \ldots
\]
mit $R_{i} \cdot R_{j} \subseteq R_{i+j} \quad \forall i,j \ge 0$ gibt.\\
Ein Element $f \in R$ hei\ss t \emph{homogen}, wenn es ein $i$ mit $f \in R_{i}$ gibt. Jedes $f \in R$ besitzt eine eindeutige Zerlegung in seine homogenen Komponenten, $f=f_{0}+f_{1}+f_{2}+\ldots$ mit $f_{i} \in R_{i} \; \forall i$, wobei nur endlich viele $f_{i} \ne 0$ sind. Ist $f \ne 0$ und $d$ die gr"o\ss te Zahl mit $f_{d} \ne 0$, so schreibt man auch $\deg f=d$ und nennt $d$ den \emph{Grad} von $f$.
\item $R$ hei\ss t \emph{zusammenh"angend} wenn zus"atzlich $R_{0}=K \cdot 1$ gilt.
\item Ist $R$ eine graduierte Algebra mit $\dim_{K} (R_{d})<\infty \quad \forall d$, so hei\ss t die formale Potenzreihe
\[
H(t):=\sum_{d=0}^{\infty} \dim_{K} (R_{d})t^{d} \in {\mathbb Z}[[t]]
\]
die \emph{Hilbertreihe} von $R$.
\end{enumerate}
\end{Def}

\noindent Das Standardbeispiel ist der Polynomring
$K[X]=K[X_{1},\ldots,X_{n}]$ in $n$ Va\-ri\-ab\-len. Die Graduierung ist dabei durch Zerlegung in homogene Polynome bez"uglich des (total-)Grads gegeben. 

Sei nun $V$ ein $G$-Modul und $X=(X_{1},\ldots,X_{n})$  eine Basis des \emph{Duals} $V^{*}$. Wie im Abschnitt "uber symmetrische Potenzen beschrieben, fasst man nun die Basisvektoren $X_{i}$ als Variablen eines Polynomrings $K[V]:=K[X]=K[X_{1},\ldots,X_{n}]$  auf, und l"asst hierauf $G$ wie induziert linear operieren (da $\dim K[V]=\infty$  handelt es sich aber um keinen $G$-Modul). Die Gradu\-ier\-ung ist dabei durch die symmetrischen Potenzen des Duals  gegeben, also
\[
K[V]=\bigoplus_{d=0}^{\infty} S^{d}(V^{*})=\bigoplus_{d=0}^{\infty} K[V]_{d}
\]
wobei $S^{0}(V^{*}):=K$ die triviale Operation tr"agt. Man nennt $K[V]$ den \emph{Polynomring von $V$}. Entsprechend definieren wir den \emph{Invariantenring von $V$} als
\begin{eqnarray*}
K[V]^{G}&:=&\left\{ f \in K[V]: \sigma \cdot f =f \quad \forall \sigma \in G \right\}\\
&=&\bigoplus_{d=0}^{\infty} S^{d}(V^{*})^{G} = \bigoplus_{d=0}^{\infty} K[V]_{d}^{G}.
\end{eqnarray*}
Ein wichtiges Ziel der Invariantentheorie ist es, den Invariantenring bzw. seine Struktur m"oglichst gut zu bestimmen.

\begin{Bemerkung}
Eine zusammenh"angende graduierte Algebra $R$ ist genau dann noethersch, wenn sie endlich erzeugt ist, d.h. wenn es $f_{1},\ldots,f_{n} \in R$ gibt mit $R=K[f_{1},\ldots,f_{n}]$.
\end{Bemerkung}

\Bew Siehe Lang \cite{Lang}, Proposition X.5.2, S. 427. \qed

\begin{SatzDef}[Noether-Normalisierung] \label{Noether}
Sei $R$ eine noe\-ther\-sche, zusammenh"angende graduierte (und kommutative!) $K$-Algebra. Dann gibt es homogene Elemente $f_{1},\ldots,f_{n} \in R$, die algebraisch unabh"angig "uber $K$ sind, so dass  $R$ endlich erzeugt ist als Modul "uber $A:=K[f_{1},\ldots,f_{n}]$, d.h. es gibt $g_{1},...,g_{m} \in R$ (O.E.) homogen mit
\[
R=\sum_{i=1}^{m}Ag_{i}.
\]
Dann hei\ss t $\left\{f_{1},\ldots,f_{n}\right\}$ ein \emph{homogenes Parametersystem} (hsop), wobei die Zahl $n=:\dim(R)$ eindeutig bestimmt ist und die \emph{Krull-Dimension} von $R$ hei\ss t. 

Eine Menge von Elementen aus $R$ hei\ss t ein \emph{partielles homogenes Para\-meter\-syst\-em} (phsop) wenn sie sich zu einem hsop erg"anzen l"asst.

Falls $R$ ein Integrit"atsring ist, gilt
$n=trdeg\left(Q\left(R\right)/K\right)$, wobei $Q(R)$ der Quotien\-tenk"or\-per zu $R$ ist. 
\end{SatzDef}

\Bew Siehe Derksen/Kemper \cite{Kem0}, Definition 2.4.6, Corollary 2.4.8 auf S. 61 und Proposition 3.3.1 auf S. 80, sowie Eisenbud \cite{Eisenbud}, Chapter 13, Theorem A, S. 290. Siehe auch Kemper \cite{Kem4} S. 4-5. \qed\\

\noindent Damit wir diesen Satz stets voll verwenden k"onnen, vereinbaren wir ab jetzt\\

\noindent {\bf Standardvoraussetzung.} Im Folgenden bezeichnen wir mit $R$ stets eine noethersche, zusammenh"angende graduierte und nullteilerfreie (kommutative) $K$-Algebra.\\

\noindent Die Voraussetzungen dieses Satzes werden z.B. offenbar von $K[V]$ erf"ullt. Eine hinreichende Bedingung, dass  $K[V]^{G}$ als Algebra endlich erzeugt ist (Hilberts 14. Problem, vgl. Kemper \cite{Kem4}, S. 2), ist dass $G$ reduktiv (siehe sp"ater) ist (und in dieser Arbeit interessieren wir uns nur f"ur redukive Grup\-pen). Nach obiger Be\-mer\-kung ist dann jedenfalls $K[V]^{G}$ noethersch, so dass der letzte Satz anwendbar ist. Man nennt $f_{1},\ldots,f_{n}$ dann auch \emph{Prim"ar\-in\-va\-ri\-an\-ten} und $g_{1},\ldots,g_{m}$ \emph{Se\-kun\-d"ar\-in\-va\-rian\-ten}. Kennt man beide, so kann man f"ur jeden gegebenen Grad $d$ ein Erzeugen\-den\-syst\-em von $K[V]_{d}^{G}$ als $K$-Vektorraum angeben: Offenbar ist
\[
K[V]_{d}^{G}=\langle f_{1}^{\alpha_{1}}\cdot\ldots\cdot f_{n}^{\alpha_{n}}g_{i}: \sum_{j=1}^{n}\alpha_{j}\deg(f_{j}) + \deg(g_{i})=d, \quad \alpha_{j} \ge 0,i=1..m \rangle.
\]
Dann hat man nat"urlich erst recht ein Erzeugendensystem von $K[V]^{G}$ als $K$-Algebra, es folgt dann n"amlich
\[
K[V]^{G}=K[f_{1},\ldots,f_{n},g_{1},\ldots,g_{m}].
\]
\noindent (Dasselbe gilt nat"urlich entsprechend auch f"ur eine
beliebige Algebra $R$. Wir formulieren die meisten Aussagen ab jetzt
meist f"ur $K[V]^{G}$, da wir nur an diesem Fall interessiert sind,
vermerken aber, wenn sie auch f"ur $R$ beliebig im Sinne unserer
Standardvoraussetzung gelten. Sprechweisen wie \emph{Polynome} und
\emph{Prim"ar-} bzw. \emph{Sekund"arinvarianten} sind dann nat"urlich
entsprechend zu modifizieren). Damit kennt man den Invariantenring
also schon recht gut. Obiges Er\-zeu\-gen\-densyst\-em ist jedoch im
Allgemeinen keine Basis. Damit dies der Fall ist, ben"otigt man eine
zus"atzliche Eigenschaft:

\newpage
\subsubsection{Die Cohen-Macaulay Eigenschaft}
\begin{SatzDef} \label{CMDef}
In der Situation von Satz \ref{Noether} hei\ss t $R$ (z.B. $K[V]^{G}$) \emph{Cohen-Macaulay} (CM) wenn $R$ sogar frei ist als Modul "uber $A$, d.h. man kann die  (Sekun\-d"ar\-in\-varian\-ten) $g_{i}$ so w"ahlen, dass
\[
R=\bigoplus_{i=1}^{m}Ag_{i}
\]
gilt. Gilt diese Eigenschaft f"ur \emph{ein} $A$ definierendes homogenes Parameter\-system, so auch f"ur \emph{jedes} andere (also ist der Begriff \CM wohl\-definiert).
\end{SatzDef}

\Bew Siehe Derksen/Kemper \cite{Kem0}, Proposition 2.5.3 (a),(c) und (d), S. 63. \qed\\

\begin{Satz} \label{CMBasis}
Ist $K[V]^{G}$ (bzw. $R$) Cohen-Macaulay, so ist mit obigen Bezeichnungen
\[
\left\{f_{1}^{\alpha_{1}}\cdot\ldots \cdot f_{n}^{\alpha_{n}}g_{i}: \sum_{j=1}^{n}\alpha_{j}\deg(f_{j}) + \deg(g_{i})=d, \quad \alpha_{j} \ge 0,i=1..m \right\}
\]
eine $K$-Basis von $K[V]_{d}^{G}$ (bzw. $R_{d}$).
\end{Satz}

\Bew  W"are das nicht der Fall, so g"abe es Polynome $F_{1},\ldots,F_{m} \in K[X_{1},\ldots,X_{n}]$ nicht alle gleich $0$, mit
\[
\sum_{i=1}^{m}F_{i}(f_{1},\ldots,f_{n})g_{i}=0.
\]
Da die Summe $\bigoplus_{i=1}^{m}Ag_{i}$ direkt ist folgt daraus aber
\[
F_{i}(f_{1},\ldots,f_{n})=0 \quad \forall i=1..m.
\]
Da ein $F_{i} \ne 0$ w"are das ein Widerspruch zur algebraischen Unabh"angigkeit von $f_{1},\ldots,f_{n}$. \qed\\

\begin{Korollar} \label{HReihe}
Ist der Invariantenring $K[V]^{G}$ (bzw. $R$) Cohen-Macaulay, so berechnet sich seine Hilbertreihe mit obigen Bezeichnungen zu
\[
H(t)=\frac{t^{\deg (g_{1})}+\ldots+t^{\deg (g_{m})}}{\left(1-t^{\deg(f_{1})} \right) \cdots \left(1-t^{\deg(f_{n})} \right)}.
\]
\end{Korollar}

\Bew Durch entwickeln in die geometrische Reihe erhalten wir f"ur die rechte Seite
\[
\left( t^{\deg (g_{1})}+\ldots+t^{\deg (g_{m})}\right)\cdot  \left( \sum_{\alpha_{1}=0}^{\infty} t^{\alpha_{1} \deg(f_{1})} \right) \cdots  
\left( \sum_{\alpha_{n}=0}^{\infty} t^{\alpha_{n} \deg(f_{n})} \right).
\]
Der Koeffizient von $t^{d}$ dieser ausmultiplizierten Reihe ist gleich der M"achtig\-keit der in Satz \ref{CMBasis} angegebenen Basis von $K[V]_{d}^{G}$, also gleich der Koeffizient der Hilbertreihe.  \qed\\

\noindent Ein vollst"andiges hsop l"asst sich f"ur einen
Invariantenring nur in den wenigsten F"allen bestimmen, meist kennt
man nur ein phsop. Daher ist folgende Aussage von Bedeutung, um
gegebenenfalls (bei Nichterf"ullung) die \CM Eigenschaft zu widerlegen.

\begin{Bemerkung}
Ist $R$ Cohen-Macaulay und ist mit den Bezeichnungen von Definition \ref{CMDef} und Satz \ref{Noether} $f_{1},\ldots,f_{k}$ ein phsop, so ist $R$ auch frei (aber f"ur $k<n$ nicht mehr endlich erzeugt) "uber $\tilde{A}:=K[f_{1},\ldots,f_{k}]$, und zwar mit Basis
\[
\left\{f_{k+1}^{\alpha_{k+1}}\cdot\ldots \cdot f_{n}^{\alpha_{n}}g_{i}:  \quad \alpha_{j} \ge 0,\,  j=(k+1),\ldots,n, \, i=1,\ldots,m \right\}
\]
\end{Bemerkung}
\vspace{0.2cm}
\noindent Die Umkehrung hiervon gilt nat"urlich nicht - z.B. ist f"ur
$k=0$ der Vektorraum $R$ stets frei "uber $\tilde{A}=K$, ohne dass $R$
dazu Cohen-Macaulay sein m"usste.\\

\Bew Dass es sich um ein Erzeugendensystem handelt, ist klar. Die Annahme, dass das angegebene System nicht linear unabh"angig "uber $\tilde{A}$ ist, f"uhrt man  genau wie im Beweis zum Satz \ref{CMBasis} zum Widerspruch. \qed\\

\noindent Aufgrund dieser und weiterer Eigenschaften ist man nat"urlich an Kriterien f"ur die \CM Eigenschaft interessiert. Da wir ja in dieser Arbeit nicht \CM Invariantenringe konstruieren wollen, sind wir vor allem auch an notwendigen Bedingungen interessiert, um dann Beispiele anzugeben, die diese notwendigen Bedingungen verletzen.

\newpage
\subsubsection{Regul"are Sequenzen}

\begin{Def}
Eine endliche Folge von Elementen $f_{1},\ldots,f_{r} \in R$ der L"ange $r$ hei\ss t \emph{regul"are Sequenz}, falls f"ur das von ihnen erzeugte Ideal  $(f_{1},\ldots,f_{r}) \ne R$ gilt und f"ur alle $i=1..r$ die Multiplikation mit $f_{i}$ injektiv auf
$R/(f_{1},\ldots,f_{i-1})$ ist, d.h. f"ur $g\in R$ gilt
\[
f_{i}\cdot g \in (f_{1},\ldots,f_{i-1}) \Rightarrow  g \in (f_{1},\ldots,f_{i-1}).
\]
\end{Def}

\noindent Um zu entscheiden, ob eine gegebene homogene Sequenz regul"ar ist, ist folgende Bemerkung hilfreich.

\begin{Bemerkung} \label{homreg}
Bildet die Folge aus \emph{homogenen} Elementen $f_{1},\ldots,f_{r} \in R$ 
mit $(f_{1},\ldots,f_{r}) \ne R$
keine regul"are Sequenz, dann gibt es \emph{homogene} Elemente $a_{1},\ldots,a_{i}\in R, \, i \le r$ mit
\[
a_{1}f_{1}+\ldots+a_{i}f_{i}=0 \quad \textrm{und } a_{i} \not\in (f_{1},\ldots,f_{i-1}).
\]
\end{Bemerkung}

\Bew Nach Voraussetzung gibt es $i \le r, b_{i} \in R$ mit $b_{i} \not\in (f_{1},\ldots,f_{i-1})$ und $f_{i}\cdot b_{i} \in (f_{1},\ldots,f_{i-1})$, d.h. es gibt $b_{1},\ldots,b_{i-1}$ mit $b_{1}f_{1}+\ldots+b_{i}f_{i}=0$. Ist $b_{i}=h_{1}+\ldots+h_{d}$ die Zerlegung von $b_{i}$ in homogene Elemente $h_{j} \in R_{j}$, so gibt es wegen $b_{i}=h_{1}+\ldots+h_{d} \not\in(f_{1},\ldots,f_{i-1})$ einen Index $j$ mir $h_{j} \not\in (f_{1},\ldots,f_{i-1})$. Setze dann $a_{i}:=h_{j}$ und w"ahle nun f"ur $k<i$ jeweils $a_{k}$ als diejenige homogene Komponente von $b_{k}$ mit $\deg a_{k} + \deg f_{k} = \deg a_{i} + \deg f_{i}$ bzw. $a_{k}=0$, falls dies nicht m"oglich ist. Da $R=R_{0} \oplus R_{1} \oplus \ldots$ (direkte Summe) gilt dann wegen $b_{1}f_{1}+\ldots+b_{i}f_{i}=0$ auch f"ur die homogenen Komponente
$a_{1}f_{1}+\ldots+a_{i}f_{i}=0$.

\qed

\begin{Satz} \label{regPerm}
Sei $n=\dim(R)$ die Krull-Dimension von $R$.
\begin{enumerate}
\renewcommand{\labelenumi}{(\alph{enumi})}
\item  F"ur jede regul"are Sequenz aus \emph{homogenen} Elementen $f_{1},\ldots,f_{r}$ und jedes $\pi \in S_{r}$ ist auch die Permutation
$f_{\pi(1)},\ldots,f_{\pi(r)}$ eine regul"are Sequenz
\item Jede \emph{homogene} regul"are Sequenz  ist ein phsop. Insbesondere gibt es keine homogene regul"are Sequenz einer L"ange $>n$.
\item $R$ ist genau dann Cohen-Macaulay, wenn jedes phsop eine regul"are Sequenz ist.
\item $R$ ist Cohen-Macaulay genau dann, wenn eine homogene regul"are Sequenz  der L"ange $n$ existiert.
\end{enumerate}
\end{Satz}

\Bew (a) Kunz \cite{Kunz}, Korollar E.16 und Satz E.17, S. 257.

 (b) Kemper \cite{Kem5}, Lemma 1.5b, S. 337.

 (c) Derksen/Kemper \cite{Kem0}, Proposition 2.5.3 (a) und (b), S. 63.

 (d) ist eine h"aufig verwendete alternative Definition der \CM Eigenschaft, siehe etwa Benson \cite{Benson}, Definition 4.3.1 und Theorem 4.3.5, S. 50f. \qed\\

\noindent Interessant f"ur uns ist Teil (c) dieses Satzes. Wir werden einen Invariantenring mit einem phsop konstruieren, das keine regul"are Sequenz ist - nach (c) kann dieser dann nicht \CM sein. Dabei muss man jedoch vorsichtig sein: Wir ben"otigen ein phsop in $K[V]^{G}$ (und nicht in $K[V]$), das dort keine regul"are Sequenz ist. F"ur sogenannte reduktive Gruppen ist dieses Problem aber handhabbar - die Eigenschaft \emph{reduktiv} wollen wir jedoch hier nicht genauer definieren (siehe hierf"ur etwa Derksen/Kemper \cite{Kem0}, Abschnitt 2.2, S. 44-50), sondern wir begn"ugen uns damit, dass die uns interessierenden Gruppen SL$_{n}(K)$ und GL$_{n}(K)$ diese Eigenschaft haben (ebenda, S. 44 unten), so dass f"ur sie folgendes Resultat gilt :

\begin{Lemma} \label{redphsop}
Sei $G$ eine reduktive Gruppe und  $V$ ein $G$-Modul.\\
Wenn $a_{1},...,a_{k} \in K[V]^{G}$ ein phsop in $K[V]$ bilden, dann auch in $K[V]^{G}$. 
\end{Lemma}

\Bew Siehe Kemper \cite{Kem2}, Lemma 4, S. 89. \qed\\

\noindent Ob nun aber homogene Elemente $a_{1},\ldots,a_{k}$ ein phsop in $K[V]$ bilden, l"asst sich relativ leicht bestimmen:\\

\begin{Satz} \label{Krull}
Seien $a_{1},\ldots,a_{k} \in R$ homogen vom Grad $\ge 1$. Dann sind "aquivalent:
\begin{enumerate}
\renewcommand{\labelenumi}{(\alph{enumi})}
\item $\left\{a_{1},\ldots,a_{k} \right\}$ ist ein phsop in $R$.
\item F"ur die Krulldimensionen gilt
\[
\dim(a_{1},\ldots,a_{k}):=\dim(R/(a_{1},\ldots,a_{k}))=\dim(R)-k.
\]
\end{enumerate}
Diese Dimension ist genau dann gleich $0$, wenn es sich sogar um ein hsop handelt. 
\end{Satz}
  
\Bew Siehe Kemper \cite{Kem4}, Proposition 2.3, S. 4.\qed\\

\noindent F"ur die Variablen $X_{1},\ldots,X_{k}$ eines Polynomrings $K[X_{1},\ldots,X_{n}]$ gilt bei\-spiels\-wei\-se
\[
K[X_{1},\ldots,X_{n}]/(X_{1},\ldots,X_{k}) \cong K[X_{k+1},\ldots,X_{n}]
\]
und daher
\[
\dim(X_{1},\ldots,X_{k})=\dim K[X_{k+1},\ldots,X_{n}] = \textrm{trdeg}K(X_{k+1},\ldots,X_{n})/K=n-k
\]
d.h. $X_{1},\ldots,X_{k}$ ist ein phsop. F"ur $k=n$ haben wir ein hsop und wir sehen, dass ein Polynomring \CM ist (frei "uber sich selbst mit Basis $g_{1}=1$).\\

Wie an diesem Beispiel angedeutet, liegt der Nutzen der letzten beiden Resultate darin, dass man damit im Polynomring relativ leicht entscheiden kann, ob ein phsop vorliegt. (Und da wir die Gruppe $G$ stets als reduktiv voraussetzen wollen, haben wir dann aufgrund des Lemmas  auch ein phsop im Invariantenring. In dieser Arbeit ist dies die einzig verwendete Methode, in diesem ein solches zu konstruieren.)

Es l"asst sich n"amlich die Dimension in (b) f"ur einen Polynomring $R=K[V]$ algorithmisch bestimmen, siehe etwa Derksen/Kem\-per \cite{Kem0}, Algorithm 1.2.4, S. 17. 

Ein solcher Algorithmus ist auch in vielen Com\-pu\-ter\-algebra\-syst\-em\-en im\-plementiert, etwa kann man mit dem Befehl {\tt Dimension} des Computer\-algebra\-syst\-ems \Magma die Dimension eines Ideals direkt be\-rech\-nen, so dass man die Suche nach phsops mit relativ wenig Aufwand automatisieren kann. (Wir geben daf"ur im n"achsten Abschnitt ein sehr einfaches  "`Greedy"'-Verfahren an).

Auch die "Uberpr"ufung auf Regularit"at des gefundenen phsops werden wir automatisieren. Dabei stellt sich sp"ater eine Frage, deren Antwort wir schon jetzt in Form eines Lemmas geben:

\begin{Lemma} \label{regPot}
Seien $f_{1},\ldots,f_{k} \in R$ \emph{homogen} und $p_{1},\ldots,p_{k} \ge 1$. 
\begin{enumerate}
\item Dann sind "aquivalent:
\begin{enumerate}
\item $\{f_{1},\ldots,f_{k}\}$ ist ein phsop.
\item $\{f_{1}^{p_{1}},\ldots,f_{k}^{p_{k}}\}$ ist ein phsop.
\end{enumerate}
\item Ebenfalls "aquivalent sind:
\begin{enumerate}
\item $f_{1},\ldots,f_{k}$ ist eine regul"are Sequenz.
\item $f_{1}^{p_{1}},\ldots,f_{k}^{p_{k}}$ ist eine regul"are Sequenz.
\end{enumerate}
oder negativ ausgedr"uckt:
\[
f_{1},\ldots,f_{k} \textrm{ nicht regul"ar} \Leftrightarrow f_{1}^{p_{1}},\ldots,f_{k}^{p_{k}} \textrm{ nicht regul"ar.}
\]
\end{enumerate}
\end{Lemma}

\Bew 
{\it Zu 1.} Es gelte (a). Wir erg"anzen $\{f_{1},\ldots,f_{k}\}$ zu einem hsop\\
$\{f_{1},\ldots,f_{k},f_{k+1},\ldots,f_{n}\}$, d.h. diese Elemente sind algebraisch unabh"angig "uber $K$ und $R$ ist endlich erzeugt als Modul "uber $A:=K[f_{1},\ldots,f_{k},f_{k+1},\ldots,f_{n}]$, wobei $\mathcal{B}$ das Erzeu\-genden\-syst\-em sei. Dann sind aber auch $f_{1}^{p_{1}},\ldots,f_{k}^{p_{k}},f_{k+1},\ldots,f_{n}$ algebraisch unabh"angig "uber $K$, und $R$ ist endlich erzeugt als Modul "uber 
\[
\tilde{A}:=K[f_{1}^{p_{1}},\ldots,f_{k}^{p_{k}},f_{k+1},\ldots,f_{n}]
\]
 mit dem Erzeu\-gen\-den\-syst\-em 
\[
\left\{f_{1}^{e_{1}}\cdot \ldots\cdot f_{k}^{e_{k}}\cdot b: \quad 0\le e_{i} < p_{i},\, i=1..k, \, b\in \mathcal{B}\right\}
\]
Also ist $f_{1}^{p_{1}},\ldots,f_{k}^{p_{k}},f_{k+1},\ldots,f_{n}$ ein hsop und daher $f_{1}^{p_{1}},\ldots,f_{k}^{p_{k}}$ ein phsop.

Nun gelte (b). Wir erg"anzen das phsop $f_{1}^{p_{1}},\ldots,f_{k}^{p_{k}}$ zu einem hsop 
\[
f_{1}^{p_{1}},\ldots,f_{k}^{p_{k}},f_{k+1},\ldots,f_{n},
\]
d.h. diese Elemente sind algebraisch unabh"angig "uber $K$ und $R$ ist endlich erzeugt als Modul "uber
\[
A:=K[f_{1}^{p_{1}},\ldots,f_{k}^{p_{k}},f_{k+1},\ldots,f_{n}].
\]  
Dann ist $R$ erst recht endlich erzeugt als Modul "uber
\[
\tilde{A}:=K[f_{1},\ldots,f_{k},f_{k+1},\ldots,f_{n}],
\]  
und es sind (nichttrivialerweise) auch 
\[
f_{1},\ldots,f_{k},f_{k+1},\ldots,f_{n}
\]
algebraisch unabh"angig "uber $K$. Da n"amlich die K"orpererweiterung
\[
K\left(f_{1},\ldots,f_{n}\right)/K\left(f_{1}^{p_{1}},\ldots,f_{k}^{p_{k}},f_{k+1},\ldots,f_{n}\right)
\]
algebraisch ist, sind die Transzendenzgrade beider K"orper "uber $K$ gleich, und der kleinere K"orper wird von einer Transzendenzbasis erzeugt, also
\[
n=\textrm{trdeg}_{K}K\left(f_{1}^{p_{1}},\ldots,f_{k}^{p_{k}},f_{k+1},\ldots,f_{n}\right) =\textrm{trdeg}_{K} K\left(f_{1},\ldots,f_{n}\right).
\]
Da sich aber aus dem Erzeugendensystem $f_{1},\ldots,f_{n}$ einer K"orpererweiterung eine (algebraisch unabh"angige) Transzendenzbasis (der L"ange des Transzen\-denz\-grades $n$) ausw"ahlen l"asst, sind $f_{1},\ldots,f_{n}$ ebenfalls algebraisch unab\-h"ang\-ig. (Vgl. \cite{Bosch}, Abschnitt 7.1, insb. Lemma 4). Also ist $\{f_{1},\ldots,f_{n}\}$ ein hsop und damit $\{f_{1},\ldots,f_{k}\}$ ein phsop.

{\it Zu 2. }Gelte (a), sei also $f_{1},\ldots,f_{k}$ eine regul"are Sequenz. Dann ist auch $f_{1},\ldots,f_{k-1},f_{k}^{p_{k}}$ eine regul"are Sequenz: Sei n"amlich $g\in R$ mit
\[
g f_{k}^{p_{k}} = (g f_{k}^{p_{k}-1}) f_{k} \in (f_{1},\ldots,f_{k-1}).
\]
Da $f_{1},\ldots,f_{k}$ regul"ar, folgt also auch $(g f_{k}^{p_{k}-1}) \in (f_{1},\ldots,f_{k-1})$ und damit sukzessive $\left(g f_{k}^{p_{k}-2}\right),\ldots,g \in (f_{1},\ldots,f_{k-1})$. Da $f_{1},\ldots,f_{k-1}$ regul"ar bleibt, ist also auch $f_{1},\ldots,f_{k-1},f_{k}^{p_{k}}$ eine regul"are Sequenz, und nach Satz \ref{regPerm} (a) dann auch $f_{k}^{p_{k}},f_{1},\ldots,f_{k-1}$. Dieselbe Argumen\-ta\-tion hierauf wiederholt angewandt liefert (b).

Gelte nun also (b). Dann ist auch $f_{1}^{p_{1}},\ldots,f_{k-1}^{p_{k-1}},f_{k}$ regul"ar: Denn f"ur $g\in R$ mit
\[
g f_{k} \in \left( f_{1}^{p_{1}},\ldots,f_{k-1}^{p_{k-1}} \right)
\]
gilt nat"urlich erst recht
\[
g f_{k}^{p_{k}} \in \left( f_{1}^{p_{1}},\ldots,f_{k-1}^{p_{k-1}} \right)
\]
und wegen der Regularit"at dann
\[
g \in \left( f_{1}^{p_{1}},\ldots,f_{k-1}^{p_{k-1}} \right).
\]
Damit ist also $f_{1}^{p_{1}},\ldots,f_{k-1}^{p_{k-1}},f_{k}$ regul"ar und wegen der Unabh"angigkeit dieser Eigenschaft von der Reihenfolge auch $f_{k},f_{1}^{p_{1}},\ldots,f_{k-1}^{p_{k-1}}$. Die obige Argumen\-ta\-tion hierauf wiederholt angewandt liefert schlie\ss lich (a). \qed
\newpage

\subsubsection{Linear reduktive Gruppen}

\noindent Das st"arkste Resultat f"ur die \CM Eigenschaft stammt von Hochster und Roberts f"ur eine besondere Klasse von Gruppen.

\begin{Def}
Eine lineare algebraische Gruppe $G$ hei\ss t \emph{linear reduktiv}, wenn f"ur jeden $G$-Modul $V$ zu jedem Untermodul $U$ von $V$ ein $G$-invariantes Komplement existiert, d.h. ein Untermodul $W$ von $V$ mit $V=U \oplus W.$
\end{Def}

\noindent Es gilt nun der ber"uhmte
\begin{Satz}[Hochster und Roberts \cite{HoRo}]
Ist $G$ linear reduktiv, so ist $K[V]^{G}$ Cohen-Macaulay f"ur jeden $G$-Modul $V$. \qed
\end{Satz}

\noindent Damit ist klar, dass wir nicht \CM Invariantenringe nur f"ur nicht linear reduktive Gruppen konstruieren k"onnen. Da wir aber bei der Konstruktion in jedem Fall Lemma \ref{redphsop} verwenden wollen, muss die unter\-such\-te Gruppe zumindest reduktiv sein. Daher ist folgendes Resultat von Bedeutung.

\begin{Bemerkung}
In Charakteristik 0 ist eine lineare algebraische Gruppe genau dann linear reduktiv, wenn sie reduktiv ist.
\end{Bemerkung}

\Bew Siehe Derksen/Kemper \cite{Kem0}, Theorem 2.2.13, S. 50. \qed\\

\noindent Aufgrund dieser Bemerkung und des Satzes von Hochster und Roberts k"onnen wir also mit unseren Methoden in Charakteristik 0 keinen nicht \CM Invariantenring konstruieren.

\newpage
\section{Ein  Test auf die Cohen-Macaulay Eigenschaft}
\noindent Wir wollen in diesem Abschnitt einen Algorithmus angeben, der $K[V]^{G}$ auf die Cohen-Macaulay Eigenschaft testet. Der Algorithmus kann dabei auf zwei Weisen enden:
\begin{itemize}
\item Es kann keine Aussage gemacht werden, oder
\item $K[V]^{G}$ ist (definitiv) nicht Cohen-Macaulay.
\end{itemize}
Der Algorithmus basiert auf folgendem Prinzip: Die Hilbertreihe wird bis zu einem gegebenen Grad $d_{\max}$ sowohl berechnet als auch gesch"atzt. F"ur \CM Ringe stimmen die gesch"atzte und tats"achliche Hilbertreihe "uberein - bei nicht "ubereinstimmen ist der Ring also nicht Cohen-Macaulay. 

Die Sch"atzung geschieht dabei durch Berechnung eines phsops und danach durch Berechnung eines Erzeugendensystems bis zum Grad $d_{max}$ von $K[V]^{G}$ als $A=K[phsop]$-Modul. Im Falle eines \CM Ringes wird das bestimmte Erzeugendensystem sogar linear unabh"angig "uber $A$ sein, so dass sich die ersten $d_{max}+1$ Koeffizienten der Hilbertreihe wie folgt berechnen lassen:
\begin{Lemma} \label{HilbSchaetz}
Sei $K[V]^{G}$ Cohen-Macaulay mit einem phsop $f_{1},\ldots,f_{r}$. Ist $\left\{g_{i}: i=1,2,3,\ldots \right\}$ eine homogene Basis von $K[V]^{G}$ als $A=K[f_{1},\ldots,f_{r}]$-Modul mit  $\deg g_{1} \le \deg  g_{2} \le \deg g_{3} \le \ldots$
 und ist $d_{max} < \deg g_{m+1}$, so stimmen die Koeffizienten von $1,t,t^{2},\ldots,t^{d_{max}}$ der Reihe
\begin{equation} \label{pHilb}
\frac{t^{\deg (g_{1})}+\ldots+t^{\deg (g_{m})}}{\left(1-t^{\deg(f_{1})} \right) \cdots \left(1-t^{\deg(f_{r})} \right)}
\end{equation}
(Hilbertreihe des $A$-Moduls $A{g_{1}} \oplus \ldots \oplus A{g_{m}}$)
mit denen der Hilbertreihe von $K[V]^{G}$ "uberein.
\end{Lemma}

\Bew Wie im Beweis von Korollar \ref{HReihe} sieht man, dass der Koeffizient von $t^{d}, d=0,\ldots,d_{max}$ der angegebenen Reihe gleich der M"achtigkeit der Menge 
\[
\left\{f_{1}^{\alpha_{1}}\cdot\ldots \cdot f_{r}^{\alpha_{n}}g_{i}: \sum_{j=1}^{r}\alpha_{j}\deg(f_{j}) + \deg(g_{i})=d, \quad \alpha_{j} \ge 0,i=1..m \right\}
\]
ist, und diese ist aufgrund der Voraussetzungen eine $K$-Basis von $K[V]_{d}^{G}$. \qed\\

%\newpage
\subsection{Das Ger"ust des Algorithmus \emph{IsNotCohenMacaulay}} \label{AlgINC}
\noindent Damit k"onnen wir nun das Ger"ust des angesprochenen Algorithmus formulieren. Die einzelnen Schritte werden sp"ater noch genauer erl"autert.\\
{\noindent \bf Algorithmus} \emph{IsNotCohenMacaulay}\\
{\bf Eingabe:} 
\begin{itemize}
\item Eine reduktive Gruppe $G$ (in Form des definierenden Gleichungssystems von $G$ als affine Variet"at bzw. das von diesen Gleichungen erzeugte Ideal)
\item Ein $G$-Modul $V$ (in Form einer Darstellungsmatrix des \emph{Duals} $V^{*}$).
\item Der Grad $d_{max}$ bis zu welchem die Hilbertreihe gesch"atzt werden soll.
\end{itemize}
{\bf Ausgabe:} 
\begin{itemize}
\item \emph{true}, falls $K[V]^{G}$ als (definitiv!) nicht Cohen-Macaulay erkannt wurde.
\item \emph{false}, falls keine Aussage m"oglich.
\end{itemize}
{\bf BEGIN}
\begin{enumerate}
\item Berechne f"ur $d=1,\ldots,d_{max}$ ein $K$-Basis $\textrm{Inv}_{d}$ von $K[V]^{G}_{d}$ mit Hilfe des Bayer-Algorithmus (siehe Anhang).
\item W"ahle ein (m"oglichst gro\ss es) phsop $\left\{f_{1},\ldots,f_{r}\right\}$  von $K[V]$ aus $\textrm{Inv}_{1} \cup \ldots \cup \textrm{Inv}_{d_{max}}$ aus.
\item Berechne Elemente $g_{1},\ldots,g_{m}$, die den Ring $K[V]^{G}$ bis zum Grad $d_{max}$ als $A:=K[f_{1},\ldots,f_{r}]$ Modul erzeugen, auf folgende Weise:
\begin{enumerate}
\item Setze $g_{1}:=1 \in K[V]_{0}^{G}=K, m:=1$.
\item F"ur $d=1,\ldots,d_{max}$ durchlaufe $f$ die Basis $\textrm{Inv}_{d}$:
\begin{itemize}
\item Falls $f \not\in Ag_{1}+\ldots+Ag_{m}$, setze $m:=m+1, g_{m}:=f$.
\end{itemize} 
\end{enumerate}
\item Vergleiche die Koeffizienten von $1,t,t^{2},\ldots,t^{d_{max}}$ der durch Gleichung (\ref{pHilb}) gegebenen gesch"atzten Hilbertreihe mit denen der tats"achlichen; diese sind gegeben durch $\dim K[V]_{d}^{G} = |\textrm{Inv}_{d}|$.
\begin{itemize}
\item Sind alle Koeffizienten gleich, so ist keine Aussage m"oglich. Gib \emph{false} aus.
\item Stimmt ein Koeffizient der gesch"atzten nicht mit dem der echten Hilbertreihe "uberein, so ist $K[V]^{G}$ nicht Cohen-Macaulay. Gib \emph{true} aus. 
\end{itemize}
\end{enumerate}
{\bf END}\\

\noindent Die Korrektheit dieses Algorithmus folgt nun aus der

\begin{Bemerkung}
Ist $K[V]^{G}$ Cohen-Macaulay, so ist nach dem homogenen Nakayama Lemma (siehe Derksen/Kemper \cite{Kem0}, Lemma 3.5.1, S. 89) jedes minimale homogene Erzeugendensystem von $K[V]^{G}$ als $A$-Modul eine Basis - insbesondere ist das in Schritt 3. des Algorithmus Grad f"ur Grad konstruierte (partielle) Erzeugendensystem $A$-linear unab\-h"ang\-ig. Daher sind die Voraus\-setzungen von Lemma \ref{HilbSchaetz} erf"ullt, und die ge\-sch"atz\-te und tats"achliche Hilbert\-reihe stim\-men bis zum Grad $d_{max}$ "uberein.

 \qed
\end{Bemerkung}

\subsection{Bemerkungen zur Implementierung in {\sc Magma}} \label{ImplMagma}
\noindent Der angegebene Algorithmus l"asst sich sehr einfach in \Magma imple\-men\-tier\-en. Allerdings kann \Magma nur mit endlichen (also nicht algebraisch abge\-schlos\-sen\-en) K"orpern rechnen. Aufgrund der Funktionsweise des Algorithmus macht dies jedoch jedoch letztlich f"ur das Endergebnis keinen Unterschied, solange die Darstellungsmatrix nur Koeffizienten aus einem endlichen K"orper - in der Regel dem Primk"orper - enth"alt. Den K"orper mit $2$ Elementen definiert man z.B. in \Magma durch\\
{\tt K:=GF(2);}

Als n"achstes wird die algebraische Gruppe $G$ als mittels des von den die Variet"at definierenden Polynomen erzeugten Ideals angegeben. F"ur die Gruppe \SL2 definiert man also zun"achst einen Polynomring in vier Variablen, z.B.\\
{\tt Pa<a,b,c,d>:=PolynomialRing(K,4);}\\
und in diesem dass die Gruppe definierende Ideal\\
{\tt IG:=ideal<Pa|[a*d-b*c-1]>;}

Der Modul $V^{*}$ wird nun durch Angabe seiner Darstellungsmatrix in Form einer Liste ihrer Spalten definiert. Die Elemente dieser Liste sind dann Poly\-no\-me aus {\tt Pa}. So definiert man also den Modul $V^{*}=\langle X^{2},Y^2,XY \rangle$ durch (siehe Ta\-bel\-le \ref{Darstellungen})\\
{\tt A:=[[a\caret 2,c\caret 2,0],[b\caret 2,d\caret 2,0],[a*b,c*d,1]];}

Da \emph{IsNotCohenMacaulay} noch einige weitere Ausgaben macht, m"ochte man vielleicht die Variablen des Polynomrings $K[V]$ benennen. Dazu definiert man einen weiteren Polynomring, der mit $K[V]$ identifiziert werden soll. In diesem Beispiel empfiehlt sich z.B.\\
{\tt P<x2,y2,xy>:=PolynomialRing(K,3);}.

Der Aufruf von von \emph{IsNotCohenMacaulay} geschieht dann durch\\
{\tt IsNotCohenMacaulay(IG,A,dmax: P:=P);},\\
wobei {\tt dmax} angibt, bis zu welchem Grad Invarianten berechnet werden.
({\tt IsNotCohenMacaulay} befindet sich in der Datei {\tt CMTest.txt}, die vorher in \Magma geladen werden muss).

Es sei noch auf eine Besonderheit f"ur Untergruppen von \GL2 hin\-ge\-wie\-sen, die den Torus
\[
T:=\left\{ 
\left( \begin{array}{cc}
a & 0\\
0 & a^{-1}
\end{array} \right): a\in K\setminus\{0\} \right\}
\]
umfassen. In den meisten F"allen reduziert sich f"ur $\sigma= \left( \begin{array}{cccc}
a &0\\
0 & a^{-1}
\end{array} \right) \in T$ die Darstellungsmatrix $A_{\sigma}$ auf die Form
\[
A_{\sigma}=\left( \begin{array}{cccc}
a^{w_{1}} & 0 & \cdots &0 \\
0 & a^{w_{2}} & &\vdots\\
\vdots & &\ddots &0\\
0 & \cdots &0 & a^{w_{n}}
\end{array} \right) \quad \textrm{mit } w_{1},\ldots,w_{n} \in {\mathbb Z} \textrm{ (fest)}.
\]
Man kann sogar genauer zeigen, dass sich die Darstellung auf dem Torus durch eine geeignete "Ahnlich\-keits\-trans\-formation \emph{immer} auf diese Form bring\-en l"asst, doch in keinem der von mir untersuchten Beispiele wurde dies n"otig.
Bei zus"atzlicher "Ubergabe des \emph{Gewichtsvektors} $w:=(w_{1},w_{2},\ldots,w_{n})$ an \emph{IsNotCohenMacaulay} wird ein modifizierter Algorithmus zur Berechnung der Invarianten aufgerufen (wir beschreiben sein Prinzip im Anhang), der erheblich (bis zu 500 mal) schneller ist. F"ur $V^{*}=\langle X^{2},Y^2,XY \rangle$ mit $w=(2,-2,0)$ erfolgt der Aufruf dann durch\\
{\tt IsNotCohenMacaulay(IG,A,dmax: P:=P, w:=[2,-2,0]);}.

Der vollst"andige Kopf von \emph{IsNotCohenMacaulay} lautet\\
{\tt IsNotCohenMacaulay:=function(IG,A,dmax: P:=1, mdp:=1,\\ maxdp:=dmax, w:=[ ]);}

Der optionale Parameter {\tt maxdp} gibt dabei an, bis zu welchem Grad "uberhaupt phsop-Elemente bestimmt werden sollen. Normalerweise m"ochte man hier nat"urlich in allen berechneten Invarianten suchen (also Standardwert {\tt dmax}), da diese Suche aber bei gro\ss en Beispielen recht lange dauern kann und weil die interessierenden Grade f"ur das phsop meist $\le 3$ sind, kann man hier den Grad runtersetzen.

Der optionale Parameter {\tt mdp} kann in der Regel auf {\tt 1} gesetzt werden. Das Programm macht bei der phsop-Bestimmung {\tt mdp} Durchl"aufe $i=1,\ldots,mdp$, wobei das phsop dann erst aus Invarianten vom Grad $i$ aufw"arts ausgew"ahlt wird. Damit sollte gegebenfalls eine h"ohere Vielfalt an phsops erreicht werden. Meine Experimente haben jedoch zu keinem Beispiel gef"uhrt, bei dem diese Variation notwendig wurde. Ein Grund hierf"ur liegt auch in Lemma \ref{regPot}, worauf wir noch zur"uckkommen werden.\\

\noindent Wir gehen nun kurz auf die Implementierung der einzelnen Schritte des Algorithmus ein.

\subsubsection{Berechnung homogener Invarianten}
\noindent Die genauere Beschreibung dieses ersten Schritts verschieben wir auf den Anhang, wo wir das Prinzip des Algorithmus von Bayer \cite{Bayer} mit einer Modi\-fika\-tion beschrei\-ben. Vorerst sehen wir diesen als "`Black-Box"' an und demon\-strieren die Funktions\-weise nur an einem Beispiel. Wir wollen eine Basis der In\-varianten zweiten Grades von $K[V]$ mit $V^{*}=\langle X_{1},Y_{1} \rangle \oplus \langle X_{2},Y_{2} \rangle$ f"ur die Gruppe \SL2 in Charakteristik $2$ berechnen:\\
{\tt
K:=GF(2);\\
Pa<a,b,c,d>:=PolynomialRing(K,4);\\
IG:=ideal<Pa|[a*d-b*c-1]>;\\

\noindent // <X1,Y1> + <X2,Y2>\\
A:=[[a,c,0,0],[b,d,0,0],[0,0,a,c],[0,0,b,d]];\\
P<X1,Y1,X2,Y2>:=PolynomialRing(K,4);\\
HomogeneousInvariantsBayerTORUS(IG,A,2: P:=P, w:=[1,-1,1,-1]);\\
}\\
liefert als Ausgabe\\
{\tt
[\\
    X1*Y2 + Y1*X2\\
].\\}
Der vollst"andige Kopf der Funktion lautet\\
{\tt HomogeneousInvariantsBayerTORUS:=function(IG,A,d: P:=1, w:=[])}\\
wobei die Parameter wie bei {\tt IsNotCohenMacaulay} sind: {\tt IG} gibt die Gruppe an, {\tt A} den Modul $V^{*}$, {\tt d} den Grad, in welchem eine $K$-Basis von Invarianten berechnet werden soll und mit $P$ (optional) benennt man die Variablen von $K[V]$. Die Angabe des Gewichtsvektors {\tt w} ist ebenfalls optional, f"uhrt aber zum Aufruf des schnelleren, modifizierten Algorithmus (was bei diesem kleinen Beispiel nat"urlich nicht sp"urbar ist).

Die gefundene Invariante $X_{1}Y_{2}+X_{2}Y_{1}$ wird in dieser Arbeit noch ein gro\ss es Abenteuer erleben.

\subsubsection{Bestimmung eines phsop}
\noindent Wir bestimmen ein phsop in $K[V]$ nach folgendem "`Greedy"'- Vefahren:\\

\noindent {\bf Algorithmus} \emph{ChoosePhsop}\\
{\bf Eingabe:} 
\begin{itemize}
\item Eine endliche Folge $M=(f_{1},\ldots,f_{m})$ von (paarweise verschiedenen) Elementen aus einem Polynomring K[V] in $n$ Variablen, aus der ein phsop in $K[V]$ ausgew"ahlt werden soll. Elemente mit kleinerem Index sollen dabei bevorzugt werden.
\item In unserem Fall besteht $M$ aus einer nach Grad geordneten Auflistung der im Schritt 1. des Algorithmus berechneten Invarianten $\textrm{Inv}_{1},\ldots,\textrm{Inv}_{d}$.\\
\end{itemize}
%\newpage
{\bf Ausgabe}
\begin{itemize}
\item Ein phsop in K[V]
\end{itemize}
{\bf BEGIN}
\begin{enumerate}
\item Setzte $phsop:=\{\}$
\item F"ur $i=1,\ldots,m$
\begin{itemize}
\item Falls $phsop \cup \{f_{i}\}$ noch ein phsop ist, also falls
\[
\dim\left( \langle phsop \cup \{f_{i}\} \rangle \right) = n-|phsop|-1
\]
(vgl. Lemma \ref{Krull}), dann setze
\[
phsop:=phsop \cup \{f_{i}\}.
\]
\end{itemize}
\end{enumerate}
{\bf END}\\

\noindent Es ist klar, dass dieses Verfahren ein phsop liefert. Die Dimensionen in Schritt zwei des Algorithmus berechnet man in \Magma f"ur eine Sequenz {\tt phsop:=[p1,p2,...,pr]} von Elementen in einem Polynomring {\tt P} ($=K[V]$) (Bezeichnungen wie bisher) mit dem Befehl {\tt Dimension(ideal<P|phsop>)}. Da wir die Gruppe $G$ als reduktiv vorausgesetzt haben und in unserem Fall die $f_{i}$ Invarianten sind, erhalten wir nach Lemma \ref{redphsop} mit diesem Algorithmus sogar ein phsop in $K[V]^{G}$.

\subsubsection{Berechnung des Erzeugendensystems}
Dieser Schritt wurde schon im Algorithmus \emph{IsNotCohenMacaulay}
genau be\-schrie\-ben. Allerdings bleibt noch die Frage, wie man f"ur gegebene homogene Elemente (\Magma Notation)  {\tt g:=[g1,...gm]},  ein gegebenes phsop {\tt phsop:=[f1,...,fr]} und ein homogenes Element {\tt f} aus {\tt P}, pr"ufen kann, ob
\[  
f \in Ag_{1}+\ldots+Ag_{m} \quad \textrm{mit } A:=K[f_{1},\ldots,f_{r}]
\]
gilt. Dies geht mit dem Befehl {\tt HomogeneousModuleTest(phsop,g,[f])}, der \emph{true} liefert, falls ja, sonst \emph{false}. Intern wird dabei $f$ als $K$-Linear\-kombi\-nation von Elementen der Form $f_{1}^{p_{1}}\cdot\ldots\cdot f_{r}^{p_{r}}\cdot g_{i}$ angesetzt, die denselben Grad wie $f$ haben. Dies sind nur endlich viele, und durch Koeffizienten\-vergleich erh"alt man ein inhomogenes lineares Gleichungs\-system, welches genau dann eine L"osung hat, wenn $f$ im von den $g_{i}$ erzeugten Untermodul liegt.

\subsubsection{Berechnung der Hilbertreihe}
Es bleibt noch die Frage, wie man die Koeffizienten der gesch"atzten Hilbert\-reihe (\ref{pHilb}) berechnen kann. Dies geht deshalb, weil man formale Potenzreihen nach einem gegebenen Grad $d$ abschneiden kann. Mit den verbleibenden Ko\-effi\-zien\-ten kann man dann wie fr"uher weiter rechnen. F"ur $f=\sum_{i=0}^{\infty}a_{i}X^{i} \in K[[X]]$ schreibt man dann
\[
f=a_{0}+a_{1}X+a_{2}X^{2}+\ldots+a_{d}X^{d}+O(X^{d+1}).
\]
Kennt man nun $g \in K[[X]]$ ebenfalls nur bis zum Grad $d$, so kann man
trotzdem $f+g$ und $f \cdot g$ berechnen, aber auch nur bis zum Grad
$d$. Abstrakt gesprochen rechnet man also in dem Restklassenring
$K[[X]]/(X^{d+1})$, und obige Schreibweise wird korrekt, wenn man
"`$=$"' durch "`$\in$"' ersetzt und
$O(X^{d+1}):=(X^{d+1})$ liest, also das von $X^{d+1}$ in $K[[X]]$
erzeugte Ideal.

  In \Magma kann man dies f"ur $d=d_{max}$ mittels\\
{\tt Qt<t>:=PowerSeriesAlgebra(Rationals(),dmax+1); }\\
realisieren. Man berechnet nun die gesch"atzte Reihe (\ref{pHilb}) als abgeschnittene Potenzreihe {\tt H} in {\tt Qt}. Den $i$-ten Koeffizient von {\tt H} erh"alt man dann mit dem Befehl {\tt Coefficient(H,i)}. 

\newpage
\subsection{Erkennung nicht-regul"arer Sequenzen}
Mit \emph{IsNotCohenMacaulay} konnten alle von mir konstruierten nicht \CM In\-va\-rian\-ten\-ringe als solche erkannt werden. Denn wenn das ge\-fun\-de\-ne phsop keine regul"are Sequenz ist, so wird dies erkannt falls $d_{max}$ ge\-n"u\-gend gro\ss{} ist (Bezeichnungen wie in Abschnitt \ref{AlgINC}).

\begin{Satz} \label{Erkenn}
Falls das in Schritt 2 von \emph{IsNotCohenMacaulay} bestimmte phsop keine regul"are Sequenz in $K[V]^{G}$ ist, d.h. wenn es (O.E. homogene, vgl. Bemerkung \ref{homreg}) Elemente $a_{1},\ldots,a_{i} \in K[V]^{G}, i\le r$ gibt mit
\[
a_{1}f_{1}+\ldots+a_{i}f_{i}=0,\quad a_{i} \not\in (f_{1},\ldots,f_{i-1})_{K[V]^{G}}
\]
wobei O.E. $\deg a_{1}f_{1} = \ldots = \deg a_{i}f_{i}$, dann stimmen die gesch"atzte Hilbertreihe (\ref{pHilb}) mit der tats"achlichen nicht "uberein, und zwar sp"atestens nicht mehr ab dem Grad $\deg a_{i}f_{i}$. Setzt man also $d_{max}$ gr"o\ss er gleich diesem Grad, so erkennt \emph{IsNotCohenMacaulay} $K[V]^{G}$ als nicht Cohen-Macaulay.
\end{Satz}

\Bew Nach der Arbeitsweise von \emph{IsNotCohenMacaulay} ist das System
\begin{equation} \label{ES}
f_{1}^{j_{1}}\cdots f_{r}^{j_{r}}g_{j}: \quad 1\le j \le m,\quad \deg (f_{1}^{j_{1}}\cdots f_{r}^{j_{r}}g_{j}) \le d_{max}
\end{equation}
jedenfalls ein $K$-Erzeugendensystem von $K[V]_{0}^{G} \oplus \ldots \oplus K[V]_{d_{max}}^{G}$, und die gesch"atzte Hilbertreihe stimmt mit der echten bis zum Grad $d_{max}$ genau dann "uberein, wenn dieses System sogar eine $K$-Basis ist. Wir zeigen, dass es unter den gemachten Voraussetzungen aber $K$-linear abh"angig ist: Da $a_{j} \in K[V]^{G}$ ($j=1..i$) und $\deg a_{j} \le \deg a_{j}f_{j} \le d_{max}$, liegt $a_{j}$ in dem von dem System (\ref{ES}) erzeugten $K$-Vektorraum, d.h. es gibt $h_{jk} \in K[f_{1},\ldots,f_{r}],\quad k=1..m$ mit
\begin{equation} \label{rel}
a_{j}=h_{j1}g_{1}+\ldots+h_{jm}g_{m}.
\end{equation}
Da nach Voraussetzung $a_{i} \not\in (f_{1},\ldots,f_{i-1})_{K[V]^{G}}$, gibt es wenigstens einen Index $k$, so dass das Polynom $h_{ik}$ (in den algebraisch unabh"angigen Elementen $f_{1},\ldots,f_{r}$) wenigstens einen Term der Form $f_{1}^{0}\cdots f_{i-1}^{0}f_{i}^{s_{i}}\cdots f_{r}^{s_{r}}$ enth"alt. Damit enth"alt auch $h_{ik}f_{i}$ einen solchen Term, jedoch keines der Polynome $h_{jk}f_{j}$ f"ur $j \le i-1$ (denn $f_{j}$ kommt mit mindestens Grad $1$ vor), und damit auch $a_{j}f_{j}$ nicht f"ur $j\le i-1$. In der Linearkombination $a_{1}f_{1}+\ldots+a_{i}f_{i}=0$ kommt also (nach Zerlegung der $a_{j}$ gem"a"s (\ref{rel})) der Vektor $f_{1}^{0}\cdots f_{i-1}^{0}f_{i}^{s_{i}+1}\cdots f_{r}^{s_{r}}g_{k}$ des Erzeugendensystems (\ref{ES}) nur im Term zu $a_{i}f_{i}$ mit Koeffizient $\ne 0$ vor. Damit ist also das System (\ref{ES}) linear abh"angig. Die gesch"atzte Hilbertreihe stimmt dann mit der tats"achlichen nicht "uberein und \emph{IsNotCohenMacaulay} gibt \emph{true} aus.\qed

\newpage
\section{Das Hauptkonstruktionsverfahren}
In diesem Abschnitt beschreiben wir das Verfahren aus \cite{Kem2} zur Konstruktion eines pshop (in einem Invariantenring), das keine regul"are Sequenz bildet, und damit das Verfahren zur Konstruktion eines nicht Cohen-Macaulay In\-va\-ri\-an\-ten\-rings. Die ben"otigten Resultate aus \cite{Kem2} habe ich bereits in meinem Projekt \cite{Projekt} zusammengefasst, f"ur eine selbsttragende Darstellung soll dies hier aber nochmal wiederholt werden, teilweise mit etwas anderen, unserer einfacheren Situation besser angepassten Beweisen.

\subsection{Exakte Sequenzen}
In dieser Arbeit betrachten wir nur \emph{exakte Sequenzen} von $G$-Moduln; Dies ist eine Folge von $G$-Homomorphismen von $G$-Moduln
\[
\ldots \longrightarrow U \stackrel{\varepsilon}{\longrightarrow} V \stackrel{\pi}{\longrightarrow} W \longrightarrow \ldots
\]
so, dass das Bild des Vorg"angers stets der Kern des Nachfolgers ist, also
\[
\textrm{Bild } \varepsilon = \textrm{Kern } \pi.
\]
Wir interessieren uns hier vor allem f"ur \emph{kurze exakte Sequenzen}; das sind exakte Sequenzen der Form
\[
0 \longrightarrow U \stackrel{\varepsilon}{\longrightarrow} V \stackrel{\pi}{\longrightarrow} W \longrightarrow 0.
\]
Der triviale Homomorphismus $0 \longrightarrow U$ hat dabei Bild $\{0\}$, und damit hat $\varepsilon$ trivialen Kern und muss daher injektiv sein. Man kann daher $\varepsilon$ auch als Injektion $U\hookrightarrow V$ ansehen. Der triviale Homomorphismus $W \longrightarrow 0$ hat ganz $W$ als Kern, und dies ist das Bild des damit surjektiven Homomorphismus $\pi$. Aufgrund des Homomorphiesatzes gilt dann $W\cong V/U$.

\begin{SatzDef} \label{zerfall}
Man sagt, dass eine kurze exakte Sequenz von $G$-Moduln
\[
0 \longrightarrow U \stackrel{\varepsilon}{\longrightarrow} V \stackrel{\pi}{\longrightarrow} W \longrightarrow 0
\]
\emph{zerf"allt}, wenn eine der folgenden drei "aquivalenten Eigenschaften gilt:
\begin{enumerate}
\renewcommand{\labelenumi}{(\alph{enumi})}
\item Es gibt einen $G$-Homomorphismus $\varphi: W \rightarrow V$ mit $\pi \circ \varphi=\textrm{id}_{W}$.
\item Es gibt ein Komplement zu $\varepsilon(U)$ in $V$.
\item Es gibt einen $G$-Homomorphismus $\psi: V \rightarrow U$ mit $\psi \circ \varepsilon=\textrm{id}_{U}$.
\end{enumerate}
\end{SatzDef}

\Bew Wir zeigen die "Aquivalenz $(a) \Rightarrow (b) \Rightarrow (c) \Rightarrow (a)$.

$(a) \Rightarrow (b)$. Wir zeigen $V=\varepsilon(U) \oplus \textrm{Bild } \varphi$. F"ur $v \in V$ ist $v=\left( v- \varphi \circ \pi (v) \right) + \varphi \circ \pi (v)$ mit 
$\pi \left( v- \varphi \circ \pi (v) \right)=\pi(v)-\pi(v)=0$, also liegt der erste Summand in $\textrm{Kern }\pi = \textrm{Bild } \varepsilon.$ Der zweite Summand liegt offenbar in $\varphi(W)$, also $v \in \varepsilon(U) + \textrm{Bild } \varphi$.\\
Ist $v \in \varepsilon(U) \cap \textrm{Bild } \varphi$, so ist $v=\varphi(w)$ mit $w \in W$, und wegen $v \in \varepsilon(U)=\textrm{Kern } \pi$ ist $0=\pi(v)=\pi \circ \varphi(w) =w$, also $v=\varphi(0)=0$.

$(b) \Rightarrow (c)$. Sei $V=\varepsilon(U) \oplus W'$. F"ur $v=\varepsilon(u)+w'$ mit $u \in U, w' \in W'$ definiere $\psi(v):=u$. Da die Summe direkt, $\varepsilon$ injektiv und ein $G$-Homo\-morphismus ist, ist $\psi$ wohldefiniert, erf"ullt $\psi \circ \varepsilon = \textrm{Id}_{U}$ und ist ebenfalls $G$-Homomorphismus.

$(c) \Rightarrow (a)$. Da $\pi$ surjektiv, gibt es zu $w \in W$ ein Urbild $v \in V$ mit $\pi(v)=w$. Wir definieren $\varphi(w):=v-\varepsilon \circ \psi (v)$ und zeigen zun"achst die Wohldefiniertheit: Ist auch $\pi(v')=w$ mit $v' \in V$, so ist $v-v' \in \textrm{Ker } \pi = \varepsilon(U)$, also $v-v'=\varepsilon(u)$ mit $u \in U$. Es folgt $(v-v')-\varepsilon \circ \psi (v-v')=(v-v')-\varepsilon \circ \psi (\varepsilon(u))=(v-v')-\varepsilon(u)=0$, also die Wohldefiniertheit. Ferner gilt dann wegen $\pi \circ \varepsilon=0$ auch $\pi \circ \varphi(w)=\pi(v- \varepsilon \circ \psi (v))=\pi(v)=w$, also $\pi \circ \varphi =\textrm{Id}_{W}$. Man pr"uft auch leicht nach, dass $\varphi$ ein $G$-Homomorphismus ist.

\qed

\subsection{Erste Kohomologie algebraischer Gruppen}
Im Folgenden wollen wir besonders einfache kurze exakte Sequenzen unter\-suchen. Als Hilfsmittel verwenden wir daf"ur sogenannte Kozyklen. 
Sei dazu $V$ ein $G$-Modul. 

Ein Morphismus von affinen Variet"aten $g: G \rightarrow V, \quad \sigma \mapsto g_{\sigma}$ hei\ss t \textit{1-Kozyklus}, falls 
\[
g_{\sigma\tau}=\sigma(g_{\tau})+g_{\sigma}\quad \textrm{f"ur alle } \sigma,\tau \in G.
\]
Die additive Gruppe aller 1-Kozyklen (die zugleich ein $K$-Vektorraum ist) wird mit $Z^{1}(G,V)$ bezeichnet. 

F"ur ein $v \in V$ ist durch $\sigma \mapsto (\sigma - 1)v:=\sigma(v)-v$ ein spezieller 1-Kozyklus gegeben. Die Untergruppe dieser so gebildeten 1-Kozyklen, die Menge der \emph{1-Kor"ander}, wird mit $B^{1}(G,V)$ bezeichnet, und die zugeh"orige Faktorgruppe $Z^{1}(G,V)/B^{1}(G,V)$ mit $H^{1}(G,V)$.

In dieser Arbeit bezeichnen wir einen Korand manchmal als \emph{trivialen Kozyklus}, wobei wir einen Kozyklus, der kein Korand ist dann auch als \emph{nichttrivialen Kozyklus} bezeichnen.\\

\newpage
\begin{Bemerkung} \label{Kozykl}
Sei
\[
0 \rightarrow V  \hookrightarrow  \tilde{V}  \stackrel{\pi}{\rightarrow} K \rightarrow 0
\]
eine kurze exakte Sequenz von $G$-Moduln ($K$ soll dabei stets triviale $G$-Operation haben). Diese zerf"allt genau dann, wenn f"ur ein (und dann alle) $v_{0} \in \pi^{-1}(1)$ und den durch $g_{\sigma}:=(\sigma-1)v_{0}$ definierten Kozyklus $g \in Z^{1}(G,V)$ gilt, dass $g$ sogar in $B^{1}(G,V)$ liegt.
\end{Bemerkung}

\Bew Wenn die Sequenz zerf"allt, so hat $V$ ein Komplement $W$ in $\tilde{V}$, das wegen $\textrm{dim }\tilde{V} = \textrm{dim Bild }\pi + \textrm{dim Kern }\pi = 1 + \textrm{dim  }V$ eindimensional ist. Ist $v_{0}=v+w$ mit $v \in V, w \in W$, so gibt es also zu $\sigma \in G$ ein $\lambda \in K$ mit $\sigma w=\lambda w$. Aus $1=\pi(v_{0})=\pi(w)=\sigma \pi(w)=\pi(\sigma w)=\pi(\lambda w)=\lambda \pi(w)$ folgt $\lambda=1$, und $w$ ist $G$-invariant. Es folgt $g_{\sigma}=(\sigma -1)v$, also $g \in B^{1}(G,V)$.

Ist umgekehrt $g \in B^{1}(G,V)$, also $g_{\sigma}=(\sigma-1)v$ mit $v \in V$, so ist $K(v_{0}-v)$ ein $G$-invariantes Komplement zu $V$. \qed
\\

\noindent Es ist n"utzlich, sich mit der Darstellung von $G$ auf $\tilde{V}$ vertraut zu machen: Wir erg"anzen $v_{0}$ mit Hilfe einer Basis von $V$ zu einer Basis von $\tilde{V}$. Ist dann $A_{\sigma}$ die Darstellung von $\sigma$ auf $V$ und identifizieren wir $g_{\sigma}$ mit seinem Koordinatenvektor bzgl. der Basis von $V$, so hat $\sigma$ bzgl. der Basis von $\tilde{V}$ die Darstellung  
$\left( \begin{array}{cc}
A_{\sigma} & g_{\sigma} \\
0 & 1
\end{array} \right)$. Insbesondere finden wir in der Identit"at
\[
\left( \begin{array}{cc}
A_{\sigma\tau} & g_{\sigma\tau} \\
0 & 1
\end{array} \right)
=
\left( \begin{array}{cc}
A_{\sigma} & g_{\sigma} \\
0 & 1
\end{array} \right)
\cdot
\left( \begin{array}{cc}
A_{\tau} & g_{\tau} \\
0 & 1
\end{array} \right)
=
\left( \begin{array}{cc}
A_{\sigma}A_{\tau} & A_{\sigma}g_{\tau}+g_{\sigma} \\
0 & 1
\end{array} \right)
\]
die Kozyklus-Eigenschaft von $g$ wieder. Daher kann man aus gegebenem $g \in Z^{1}(G,V)$ durch $\tilde{V}:=V \oplus K$ (hier direkte Summe von $K$-Vektorr"aumen) und  $\sigma (v,\lambda):=(\sigma v+\lambda g_{\sigma},\lambda) \cong 
\left( \begin{array}{cc}
A_{\sigma} & g_{\sigma} \\
0 & 1
\end{array} \right)
\left( \begin{array}{c}
v\\ \lambda
\end{array} \right)
$
sowie $\pi(v,\lambda):=\lambda$ eine kurze exakte Sequenz definieren. Wir sagen, $\tilde{V}$ bzw. die kurze exakte Sequenz wird von $g$ \emph{induziert}. Ist $(v,1) \in \pi^{-1}(1)$, so ist $(\sigma-1)(v,1)=(g_{\sigma}+\sigma (v) -v,0)$, so dass man aus dieser Sequenz die Restklasse von $g$ in $H^{1}(G,V)$ zur"uckgewinnen kann.\\

\noindent Was passiert wenn man $V$ mittels eines trivialen Kozyklus $g_{\sigma}=(\sigma-1)v_{0}$ mit $v_{0} \in V$ zu $\tilde{V}$ erweitert? Um dieser Frage nachzugehen, erg"anzen wir $v_{0}$ zu einer Basis von $V$. Ist dann $\sigma \mapsto A_{\sigma}$ die Darstellung bez"uglich dieser Basis, so ist die Koordinatendarstellung des Kozyklus bez"uglich dieser Basis durch $(A_{\sigma}-I)e_{1}$ gegeben, wobei $e_{1}=(1,0,\ldots,0)^{T}$ der erste Einheitsvektor der L"ange $n=\dim V$ ist. Ist $a_{\sigma}:=A_{\sigma}e_{1}$ die erste Spalte von $A_{\sigma}$, so hat also $\tilde{V}$ die Darstellung
\[
\left( \begin{array}{cc}
A_{\sigma} & g_{\sigma} \\
0 & 1
\end{array} \right)
=
\left( \begin{array}{cc}
A_{\sigma} & a_{\sigma}-e_{1} \\
0 & 1
\end{array} \right).
\]
Mit Hilfe der "Ahnlichkeitstransformation (Basiswechsel)
\[
\left(
\begin{array}{cc}
I_{n} &e_{1}\\
0&1 
\end{array} \right)
\left( \begin{array}{cc}
A_{\sigma} & a_{\sigma}-e_{1} \\
0 & 1
\end{array} \right)
\left(
\begin{array}{cc}
I_{n}&-e_{1}\\
0&1
\end{array} \right)
=
\left( \begin{array}{cc}
A_{\sigma} & 0\\
0 & 1
\end{array} \right),
\]
wobei $I_{n}$ die $n \times n$ Einheitsmatrix, erkennt man die Isomorphie $\tilde{V}\cong V\oplus K$ (direkte Summe (von Moduln) mit dem trivialen $G$-Modul) - die Erweiterung mit einem Korand bringt also nichts neues.

Man beachte auch, dass man einen Modul nur einmal nichttrivial mit dem gleichen Kozyklus erweitern kann, denn
\[
\left(
\begin{array}{c}
g_{\sigma}\\
0
\end{array} \right)
=
\left(
\left(
\begin{array}{cc}
A_{\sigma} & g_{\sigma}\\
0&1
\end{array} \right)
-I_{n+1} \right) e_{n+1},
\]
d.h. der Kozyklus $\sigma \mapsto (g_{\sigma},0)$ wird in $\tilde{V}$ trivial.\\

\noindent Man kann Kozyklen bzw. kurze exakte Sequenzen benutzen, um zu entscheiden ob ein Untermodul ein Komplement besitzt:

\begin{Prop} \label{kompl}
Sei $W$ Untermodul eines $G$-Moduls $V$, sowie $\iota \in \textrm{Hom}_{K}(V,W)$ mit $\iota|_{W}=\textrm{id}_{W}$. Dann ist durch $\sigma \mapsto g_{\sigma}:=(\sigma-1)\iota$ ein Kozyklus in $\textrm{Hom}_{K}(V,W)_{0}$ (siehe Satz \ref{Hom0}) gegeben, welcher genau dann ein Korand ist, wenn $W$ ein ($G$-invariantes) Komplement hat.

"Aquivalent dazu ist dann nat"urlich auch die Bedingung, dass die von $g$ induzierte kurze exakte Sequenz
\[
0 \rightarrow \textrm{Hom}_{K}(V,W)_{0} \hookrightarrow \widetilde{\textrm{Hom}_{K}}(V,W)_{0}=\textrm{Hom}_{K}(V,W)_{0} \oplus K\iota \stackrel{\pi}{\rightarrow} 0
\]
mit $\pi(f+\lambda \cdot \iota):=\lambda$ f"ur $f \in \textrm{Hom}_{K}(V,W)_{0}, \quad \lambda \in K$ zerf"allt. 
\end{Prop}

\Bew Zun"achst ist $g_{\sigma} \in \textrm{Hom}_{K}(V,W)_{0}$ f"ur $\sigma \in G$ zu zeigen: F"ur $w \in W$ ist
\begin{equation} \label{g0}
g_{\sigma}(w)=((\sigma-1)\iota)(w)=\sigma\left (\iota(\sigma^{-1}w) \right) -w \stackrel{\iota|_{W}=\textrm{id}_{W}}{=}\sigma\sigma^{-1}w-w=0,
\end{equation}  
also $g_{\sigma}|_{W}=0$ und damit $g_{\sigma} \in \textrm{Hom}_{K}(V,W)_{0}$.

Sei nun $g$ ein Korand, d.h. es gibt $f \in \textrm{Hom}_{K}(V,W)_{0}$ mit
$g_{\sigma}=(\sigma-1)\iota=(\sigma-1)f$. F"ur $h:=\iota-f$ folgt dann $\sigma \cdot h=h$, also $h \in \textrm{Hom}_{G}(V,W)$, und $\ker h$ ist damit ein Untermodul ($G$-invarianter Untervektorraum) von $V$. F"ur $w \in W$ gilt ferner $h(w)=\iota(w)-f(w)=w-0=w$, also $h|_{W}=id_{W}$. Da dann $h(V) = W$ folgt also $h(h(v))=h(v) \forall v \in V$. Also ist $h$ Projektion auf $W$, und in "ublicher Weise folgt nun $V=W \oplus \ker h$:
\[
\forall v \in V: \quad v=\underbrace{(v-h(v))}_{\in \ker h}+\underbrace{h(v)}_{\in W},
\] 
denn h(v-h(v))=h(v)-h(h(v))=h(v)-h(v)=0. Ferner gilt f"ur $v \in W \cap \ker h$, dass $v\stackrel{h|_{W}=\textrm{id}_{W}}{=}h(v)=0$, also ist die Summe direkt.

Gilt umgekehrt $V=W \oplus U$ mit einem $G$-invarianten Teilraum $U$, so w"ahle  
\[
f \in \textrm{Hom}_{K}(V,W)_{0} \textrm{ mit } f|_{W}=0, \; f|_{U}=\iota|_{U}.
\]
Wir zeigen $g_{\sigma}=(\sigma-1)f$: F"ur $w\in W, u \in U$ ist
\[
\begin{array}{l}
g_{\sigma}(w+u)\stackrel{(\ref{g0})}{=}g_{\sigma}(u)=((\sigma-1)\iota)(u)=\sigma(\iota(\sigma^{-1}u))-\iota(u)\stackrel{f|_{U}=\iota|_{U}}{=}\\
\sigma(f(\sigma^{-1}u))-f(u)=((\sigma-1)f)(u)\stackrel{f|_{W}=0}{=}((\sigma-1)f)(w+u),
\end{array}
\]
also Gleichheit auf $V=W \oplus U$, und $g$ ist ein Korand. \qed\\

\noindent Im Zusammenhang mit diesem Satz ist folgendes Ergebnis n"utzlich:

\begin{Bemerkung} \label{kleinKompl}
Sei $ U \le V \le W$ eine Kette von Moduln. Wenn $U$ kein Komplement in $V$ hat, so auch nicht in $W$.
\end{Bemerkung}

\Bew Angenommen, $W=U \oplus U'$. Dann gilt auch $V=U \oplus (U' \cap V)$, im Widerspruch zur Voraussetzung: F"ur $v \in V \le W$ gibt es n"amlich $u \in U, u' \in U'$ mit $v=u+u'$, also $u'=v-u \in U' \cap V$. Ausserdem ist $U \cap (U' \cap V) \subseteq U \cap U' =\{0\}$. \qed

\begin{Korollar} \label{linRed}
$G$ ist genau dann linear reduktiv, wenn $H^{1}(G,V)=0$ f"ur jeden $G$-Modul $V$ gilt.
\end{Korollar}

\Bew Sei $G$ linear reduktiv. Dann hat jeder Untermodul eines jeden Moduls $V$ ein Komplement, also zerf"allt nach Satz \ref{zerfall} (b) insbesondere auch jede von einem Kozyklus induzierte kurze exakte Sequenz, was nach Bemerkung \ref{Kozykl} bedeutet, dass jeder Kozyklus ein Korand ist, also $H^{1}(G,V)=0$.

Falls umgekehrt $H^{1}(G,V)=0$ f"ur jeden $G$-Modul $V$, so ist jeder Kozyklus ein Korand, was nach voriger Proposition aber hei\ss t, dass jeder Untermodul eines jeden $G$-Moduls ein Komplement hat, also ist $G$ linear reduktiv (einen $K$-Homomorphismus $\iota$ wie in Proposition \ref{kompl} gibt es immer). \qed\\

\noindent Sei $g \in Z^{1}(G,V)$, wobei wir die zugeh"orige Restklasse von $g$ in $H^{1}(G,V)$ ebenfalls mit $g$ bezeichnen. Ist $W$ ein weiterer $G$-Modul, $w \in W^{G}$, so ist $\sigma \mapsto w \otimes g_{\sigma}$ wegen $\sigma w =w \quad \forall \sigma \in G$ aus $Z^{1}(G,W \otimes V)$, denn
\[
w \otimes g_{\sigma \tau}=w \otimes (\sigma g_{\tau} + g_{\sigma}) \stackrel{w \in W^{G}}=\sigma (w \otimes g_{\tau}) + w \otimes g_{\sigma}. 
\]
Die zugeh"orige Restklasse in $H^{1}(G, W \otimes V)$ bezeichnen wir mit $w \otimes g$. Sie ist unabh"angig von dem Repr"asentanten von $g$ in $Z^{1}(G,V)$, da
\[
w \otimes (g_{\sigma}+\sigma (v)-v)=w \otimes g_{\sigma}+\sigma(w \otimes v)- w \otimes v 
\]
wegen $w \in W^{G}$.\\

\begin{Prop} \label{anni}
Sei $V$ ein $G$-Modul, $g \in H^{1}(G,V)$ und 
\[
0 \rightarrow V  \hookrightarrow  \tilde{V}  \stackrel{\pi}{\rightarrow} K \rightarrow 0
\]
die von einem Repr"asentanten von $g$ in $Z^{1}(G,V)$  induzierte kurze exakte Sequenz. Dann gibt es einen $G$-Modul $W$ (n"amlich $W=\tilde{V}^{*}$) und ein $0 \ne w \in W^{G}$ (n"amlich $w=\pi$) mit $w\otimes g=0$ in $H^{1}(G, W \otimes V)$. (Man sagt, $g$ wird von $w$ \emph{annulliert}).
\end{Prop}

\Bew Wir m"ussen zeigen, dass $\sigma \mapsto \pi \otimes g_{\sigma}$ ein Korand in $\tilde{V}^{*} \otimes V$ ist, dass es also ein $x \in \tilde{V}^{*} \otimes V$ mit $\pi \otimes g_{\sigma}=(\sigma-1)x \forall  \sigma \in G$ gibt. Sei dazu $\{v_{1},\ldots,v_{n}\}$ eine Basis von $V$ und  $\mathcal{B}=\{v_{1},\ldots,v_{n},v_{n+1}\}$ (mit $v_{n+1} \in \pi^{-1}(1)$)  eine Basis von $\tilde{V}$, so dass $\tilde{V}$ die Darstellung
\[
\sigma \mapsto
\left( \begin{array}{cc}
A_{\sigma} & g_{\sigma} \\
0 & 1
\end{array} \right) \textrm{ bzgl. } \{v_{1},\ldots,v_{n},v_{n+1}\}
\]
hat (wir identifizieren $g_{\sigma}$ wieder mit seinem Koordinaten(spalten)vektor). Die duale Basis von $\mathcal{B}$ ist gegeben durch $\mathcal{B'}:=\{v_{1}^{*},\ldots,v_{n}^{*},v_{n+1}^{*}=\pi\}$ (denn $\pi(V)=0$), und die Darstellung auf $\tilde{V}^{*}$ ist nach Lemma \ref{dual} gegeben durch
\[
\sigma \mapsto
\left( \begin{array}{cc}
A_{\sigma^{-1}}^{T} & 0 \\
g_{\sigma^{-1}}^{T} & 1
\end{array} \right) \textrm{ bzgl. } \{v_{1}^{*},\ldots,v_{n}^{*},v_{n+1}^{*}=\pi\}.
\]
An dieser Darstellung sehen wir auch nochmals $\pi \in (\tilde{V}^{*})^{G}$ (was auch folgt, weil $\pi$ ein $G$-Homomorphismus ist).
Gem"a\ss{} Lemma \ref{Tensor} ordnen wir einem Vektor
\[
x=\sum_{i=1}^{n+1}\sum_{j=1}^{n} x_{ij} (v_{i}^{*} \otimes v_{j})=:\sum_{i,j=1}^{n} x_{ij}(v_{i}^{*} \otimes v_{j}) + \sum_{j=1}^{n} y_{j} (\pi \otimes v_{j}) \in \tilde{V}^{*} \otimes V
\]
(mit $x_{ij},y_{j} \in K$) die Koordinaten-Matrix
\[
\left( \begin{array}{ccc}
x_{11} & \cdots &x_{1n}\\
\vdots&&\vdots\\
x_{n1} & \cdots &x_{nn}\\
y_{1} &\cdots & y_{n}
\end{array} \right) =:
\left( \begin{array}{c}
X\\
y
\end{array} \right)
\]
zu. Nun setzen wir
\begin{equation} \label{xanni}
\fbox{$ \displaystyle
x:=-\left(v_{1}^{*} \otimes v_{1} + \ldots + v_{n}^{*} \otimes v_{n}\right)
$}
\end{equation}
mit der Koordinatenmatrix 
\[
\left( \begin{array}{c}
-I_{n}\\
0_{1 \times n}
\end{array} \right),
\]
und behaupten, dass dies das gesuchte $x$ ist - wir haben also
\begin{equation} \label{piganni}
\fbox{$ \displaystyle
\pi \otimes g_{\sigma}=(\sigma-1)x \quad \forall  \sigma \in G
$}
\end{equation}
zu zeigen. W"ahrend $\pi \otimes g_{\sigma}$ die Koordinatenmatrix
\[
\left( \begin{array}{c}
0_{n\times n}\\
g_{\sigma}^{T}
\end{array} \right)
\]
besitzt, hat $(\sigma-1)x$ nach Lemma \ref{Tensor} die Koordinatenmatrix
\[
\left( \begin{array}{cc}
A_{\sigma^{-1}}^{T} & 0 \\
g_{\sigma^{-1}}^{T} & 1
\end{array} \right)
\left( \begin{array}{c}
-I_{n}\\
0_{1 \times n}
\end{array} \right)
A_{\sigma}^{T}-
\left( \begin{array}{c}
-I_{n}\\
0_{1 \times n}
\end{array} \right)=
\left( \begin{array}{c}
0_{n\times n}\\
-g_{\sigma^{-1}}^{T}A_{\sigma}^{T}
\end{array} \right),
\]
welche zu dem Vektor $\pi \otimes (-\sigma g_{\sigma^{-1}})$ geh"ort. Wir m"ussen also nur noch $-\sigma g_{\sigma^{-1}}=g_{\sigma}$ zeigen: Wegen der Kozyklus-Eigenschaft $g_{\sigma \tau} =\sigma g_{\tau} + g_{\sigma}$ gilt zun"achst f"ur das neutrale Element $\iota$ von $G$
\[
g_{\iota}=g_{\iota \iota}=\iota g_{\iota} + g_{\iota} =g_{\iota} + g_{\iota},
\]
also $g_{\iota}=0$, und daher
\[
0=g_{\iota}=g_{\sigma \sigma^{-1}}=\sigma g_{\sigma^{-1}} + g_{\sigma},
\]
also die Behauptung. \qed\\

\noindent Sch"arfere Varianten dieser Proposition findet man in den Originalarbeiten \cite{Kem2} und \cite{Kem3}; die Zusammenfassung der dortigen Beweise findet sich auch in meinem Projekt \cite{Projekt}.\\

\noindent Wir haben also zu jedem (nichttrivialen) Kozyklus eine annullierende In\-variante. Kann man umgekehrt zu einer Invarianten einen nichttrivialen Ko\-zyk\-lus finden, der von dieser annulliert wird? Wir lesen dazu einfach die bisherige Kon\-struktions\-kette r"uckw"arts: Bei gegebenem Ko\-zyk\-lus $g$ in $V$ haben wir zun"achst $\tilde{V}$ konstruiert, und $\tilde{V}^{*}$ enth"alt dann eine annullierende Invariante $\pi$. Sei nun umgekehrt ein Modul $\tilde{V}^{*}$ mit einer Invarianten $\pi$ gegeben (wobei wir den Modul wegen $W^{**}\cong W$ als Dual seines Duals schreiben k"onnen). In dem wir $\pi$ als letztes Element einer Basis nehmen, hat $\tilde{V}^{*}$ die Darstellung
\[
\sigma \mapsto
\left( \begin{array}{cc}
B_{\sigma} &0 \\
h_{\sigma}^{T} &1
\end{array} \right),
\]
wobei $h_{\sigma}$ als Spalte geschrieben sei. Die Invariante $\pi$ erzeugt den ein\-di\-men\-sion\-al\-en Untermodul (!) $K\pi$. Gem"a\ss{} Proposition \ref{kompl} betrachten wir 
\[
\iota \in \textrm{Hom}_{K}(\tilde{V}^{*},K\pi) \cong \tilde{V}^{**}\cong \tilde{V}
\]
mit der Darstellungsmatrix
\[
J:=(0,\ldots,0,1)
\]
bez"uglich der betrachteten Basis von $\tilde{V}^{*}$ ($\pi$ ist das letzte Element dieser Basis). Der nach Proposition \ref{kompl} zugeh"orige Kozyklus $g_{\sigma}=(\sigma-1)\iota$ in $\textrm{Hom}_{K}(\tilde{V}^{*},K\pi)_{0}$ hat dann die Darstellungsmatrix
\[
(1) \cdot (0,\ldots,0,1)
\left( \begin{array}{cc}
B_{\sigma^{-1}} &0 \\
h_{\sigma^{-1}}^{T} &1
\end{array} \right) -(0,\ldots,0,1)\\
=(h_{\sigma^{-1}}^{T},0). 
\]
Dabei sind die Koordinatenvektoren von $\textrm{Hom}_{K}(\tilde{V}^{*},K\pi)$ als Zeilen geschrieben. Schreiben wir sie als Spalte und beachten dass $\tilde{V}^{**}=\tilde{V}$ die Darstellung
\[
\sigma \mapsto
\left( \begin{array}{cc}
B_{\sigma} &0 \\
h_{\sigma}^{T} &1
\end{array} \right)^{-T}=
\left( \begin{array}{cc}
B_{\sigma^{-1}}^{T} & h_{\sigma^{-1}} \\
0 &1
\end{array} \right)
\]
hat, sehen wir dass dieser Kozyklus in offensichtlicher Weise einem Kozyklus in einem Untermodul $V$ von $\tilde{V}$ entspricht. Nach Proposition $\ref{kompl}$ ist dieser genau dann nichttrivial, wenn $K\pi$ kein Komplement in $\tilde{V}^{*}$ hat. Zusammen\-gefasst und griffig formuliert also:

\begin{Bemerkung} \label{InvKoz}
Der aus einer Invarianten konstruierte Kozyklus ist genau dann nichttrivial, wenn die Invariante kein Komplement besitzt. \qed
\end{Bemerkung}

\noindent Dieses Verfahren ist deshalb interessant, weil man f"ur die
Berechnung von Invarianten gute Algorithmen hat, f"ur die Berechnung von
Kozyklen jedoch nur unpraktikable. Um einen nichttrivialen Kozyklus zu finden kann man also so vorgehen, dass man einen beliebigen Modul nimmt, die Invarianten berechnet und dann pr"uft, ob die zugeh"origen Kozyklen nichttrivial sind.

Ein Beispiel, wo dieses Verfahren erfolgreich angewendet wird, findet sich in Abschnitt \ref{munupi}. Es war das erste neue Resultat dieser Arbeit und f"uhrte durch Lenkung meiner Aufmerksamkeit auf das dort betrachtete Ten\-sor\-pro\-dukt auch zum Hauptresultat dieser Arbeit.

\subsection{Der Hauptsatz}
Der folgende Satz von Kemper \cite{Kem2} wird sp"ater st"andig f"ur die Konstruktion von nicht Cohen-Macaulay Invariantenringen benutzt.
\begin{Hauptsatz} \label{HSatz}
Sei $G$ eine reduktive Gruppe und $V$ ein $G$-Modul, der ein $0 \ne g \in H^{1}(G,K[V])$ enth"alt, sowie $a_{1},a_{2},a_{3} \in K[V]^{G}$ ein phsop in $K[V]$, das  $g$ annulliert, also $a_{i}g=0 \in H^{1}(G,K[V])$ f"ur $i=1,2,3$. Dann ist $a_{1},a_{2},a_{3}$ ein phsop in $K[V]^{G}$, aber dort keine regul"are Sequenz. Insbesondere ist $K[V]^{G}$ nicht Cohen-Macaulay. 
\end{Hauptsatz}

\Bew Nach Lemma \ref{redphsop} ist $a_{1},a_{2},a_{3}$ auch ein phsop in $K[V]^{G}$. Wir zeigen, dass es dort keine regul"are Sequenz ist. 
Nach Voraussetzung sind die Ko\-zy\-klen $\sigma \mapsto a_{i} g_{\sigma}, i=1,2,3$ in $K[V]$ dort sogar Kor"ander, also gibt es $b_{i} \in K[V]$ mit 
\begin{equation} \label{bi}
(\sigma -1)b_{i}=a_{i} g_{\sigma} \quad \forall \sigma \in G,\quad i=1,2,3.
\end{equation}
Sei 
\begin{equation} \label{uij}
u_{ij}=a_{i}b_{j}-a_{j}b_{i} \textrm{ f"ur } 1 \le i < j \le 3.
\end{equation}
Offenbar ist $u_{ij} \in K[V]^{G}$ ($\sigma u_{ij}=a_{i}(b_{j}+a_{j}g_{\sigma})-a_{j}(b_{i}+a_{i}g_{\sigma})=u_{ij}$), und es gilt
\[
u_{23}a_{1}-u_{13}a_{2}+u_{12}a_{3}=
\left|
\begin{array}{ccc}
a_{1} & a_{2} & a_{3}\\
a_{1} & a_{2} & a_{3}\\
b_{1} & b_{2} & b_{3}\\
\end{array}
\right|=0.
\]
{\bf Annahme:} $a_{1},a_{2},a_{3}$ ist eine regul"are Sequenz in $K[V]^{G}.$\\
Dann liegt  $u_{12}$ in dem von $a_{1},a_{2}$ erzeugten Ideal in $K[V]^{G}$, d.h. es gibt $f_{1},f_{2} \in K[V]^{G}$ mit 
\begin{equation} \label{u12Ann}
u_{12}=a_{1}b_{2}-a_{2}b_{1}=f_{1}a_{1}+f_{2}a_{2}.
\end{equation}
Weiter sind $a_{1},a_{2}$ teilerfremde Polynome, denn $a_{1},a_{2}$ bilden ein phsop in dem \CM Ring K[V], also dort sogar eine regul"are Sequenz; Ist nun $d \in K[V]$ ein gemeinsamer Teiler, $a_{1}=dk$ so folgt $a_{2}k \in (a_{1})_{K[V]}$, also gilt mit der Regularit"at $k \in (a_{1})_{K[V]}$ , d.h. $d$ ist Einheit.  
Aus $a_{1}(b_{2}-f_{1})=a_{2}(f_{2}+b_{1})$ folgt dann, dass $a_{1}$ Teiler von $f_{2}+b_{1}$ ist, also $f_{2}+b_{1}=a_{1} \cdot h$ mit $h \in K[V]$. Nun ist
\[
a_{1} \cdot (\sigma -1)h=(\sigma -1)(a_{1}h)=(\sigma -1)(f_{2}+b_{1})=(\sigma -1)b_{1}=a_{1}g_{\sigma} \qquad \forall \sigma \in G,
\]
also $g_{\sigma}=(\sigma -1)h$. Damit ist die Restklasse von $g$ in $H^{1}(G,K[V])$ gleich 0, was im Widerspruch zur Voraussetzung steht. Also war die Annahme falsch, und $a_{1},a_{2},a_{3}$ ist nicht regul"ar. \qed\\

\noindent Falls \emph{IsNotCohenMacauly} dann $a_{1},a_{2},a_{3}$ in seiner phsop-Liste hat, so erkennt es $K[V]^{G}$ nach Satz \ref{Erkenn} als nicht Cohen-Macaulay, wenn man im Fall homo\-gener $a_{i}$ und $g_{\sigma}$ (also $g_{\sigma} \in K[V]^{G}_{d} \quad \forall \sigma \in G$ f"ur ein festes $d=:\deg(g)$) $d_{\max} \ge \deg (u_{12}a_{3})=\deg a_{1} + \deg b_{2} + \deg a_{3} $, also
\begin{equation} \label{dschaetz}
\fbox{$ \displaystyle d_{\max}\ge  \deg (a_{1}) + \deg (a_{2}) + \deg (a_{3}) + \deg (g) $}
\end{equation}
setzt. (Da $G$ graderhaltend operiert, ist $\deg b_{i}=\deg a_{i} + \deg g$).\\

\noindent Die einfachste M"oglichkeit, mit dem Hauptsatz einen nicht Cohen-Macaulay Invariantenring zu konstruieren, wird durch folgendes Korollar gegeben. Ko\-zyk\-lus und annullierendes phsop finden sich dabei im Grad 1, was zu einer relativ hohen Dimension f"ur $V$ f"uhrt. Ziel dieser Arbeit war es, die Summan\-den des n"achsten Korollars (f"ur konkrete Beispiele) in h"oheren Potenzen von (kleineren) Moduln wiederzufinden und so die Dimension zu reduzieren.

\begin{Korollar} \label{NCM}
Ist $0 \rightarrow U \rightarrow \tilde{U} \rightarrow K \rightarrow 0$ eine nicht zerfallende kurze exakte Sequenz von $G$-Moduln, so ist mit
\[
\begin{array}{rcl}
V & := & U^{*} \oplus \tilde{U} \oplus \tilde{U} \oplus \tilde{U}
\end{array}
\]
$K[V]^{G}$ nicht Cohen-Macaulay.
\end{Korollar}

\Bew Es ist
\begin{equation} \label{VStern}
V^{*}=U \oplus \tilde{U}^{*} \oplus \tilde{U}^{*} \oplus \tilde{U}^{*}
\end{equation}
und
\[
K[V]=S(V^{*})=V^{*} \oplus S^{2}(V^{*}) \oplus  S^{3}(V^{*}) \oplus \ldots,
\]
also sind $U$ und $\tilde{U}^{*}$ direkte Summanden von $K[V]$. Der Modul $U$ enth"alt einen Kozyklus $g$, der kein Korand in $U$ ist (Bemerkung \ref{Kozykl}). Mittels Einbettung ist $g$ dann auch Kozyklus in $K[V]$. Da $U$ ein direkter Summand von $K[V]$ ist, ist dann $g$ auch kein Korand in $K[V]$: Denn ist etwa $K[V]=U \oplus W$ und w"are 
\[
g_{\sigma}=(\sigma-1)(u+w) \quad \forall \sigma \in G \textrm{ mit } u \in U, w \in W
\]
so folgte wegen der direkten Summe und $g_{\sigma} \in U$ jedenfalls $(\sigma-1)w=0 \forall \sigma$, also doch $g_{\sigma}=(\sigma-1)u$ im Widerspruch dazu, dass $g$ kein Korand in $U$ ist. Die drei Kopien der nach Proposition \ref{anni} annullierenden Elemente $\pi \in \tilde{U}^{*G}$ in $K[V]^{G}$ bezeichnen wir mit $a_{1},a_{2},a_{3}$. Da sie als Elemente einer Basis und damit als unabh"angige Variablen im Polynomring dienen k"onnen, bilden sie ein phsop in $K[V]$. Da
\[
U \otimes \tilde{U}^{*} \le S^{2}\left(U \oplus \tilde{U}^{*}\right)
\]
ein Untermodul ist, sind die $a_{i}$ nach Proposition \ref{anni} Annullatoren von $g$, d.h. $a_{i}g=0$ in $H^{1}(G,K[V])$ -  die Voraussetzungen des Hauptsatzes sind also erf"ullt und $K[V]^{G}$ ist nicht Cohen-Macaulay. \qed\\

\noindent Damit k"onnen wir auch das Hauptresultat aus \cite{Kem2} beweisen, eine Umkehrung des Satzes von Hochster und Roberts:

\begin{Korollar}[Kemper \cite{Kem2}]
Sei $G$ reduktiv und $K[V]^{G}$ Cohen-Macaulay f"ur jeden $G$-Modul $V$. Dann ist $G$ sogar linear reduktiv.
\end{Korollar}

\Bew W"are $G$ nicht linear reduktiv, so g"abe es nach Korollar \ref{linRed} eine nicht zerfallende kurze exakte Sequenz $0 \rightarrow U \rightarrow \tilde{U} \rightarrow K \rightarrow 0$ von $G$-Moduln, und nach vorigem Korollar dann einen nicht Cohen-Macaulay Invariantenring - im Widerspruch zur Voraussetzung. \qed\\

\noindent Man kann diesen Satz auch so formulieren: Wenn $G$ reduktiv, aber nicht linear reduktiv ist, so gibt es einen $G$-Modul $V$ mit nicht Cohen-Macau\-lay In\-va\-ri\-an\-ten\-ring $K[V]^{G}$.

\newpage
\section{Nichttriviale Kozyklen f"ur $\textrm{SL}_{n}(K)$ und $\textrm{GL}_{n}(K)$}

\noindent Wir kommen nun zu den eigentlichen Resultaten meiner Arbeit, n"amlich der Konstruktion konkreter Beispiele f"ur die zusammen\-h"angen\-den (und reduk\-tiven) Gruppen $\textrm{SL}_{n}(K)$ und $\textrm{GL}_{n}(K)$. In meinem Projekt habe ich dazu lediglich Korollar \ref{NCM} verwandt, wozu es also gen"ugte, nichttriviale Kozyklen anzugeben. Da diese auch f"ur die Resultate dieser Diplomarbeit verwendet werden, sollen in diesem Abschnitt nochmal die Hauptergebnisse des Projekts samt Beweisen angegeben werden.

\noindent Hier also das Hauptresultat aus meinem Projekt \cite{Projekt}:
\begin{Satz} \label{PKoz}
Sei $K$ ein (algebraisch abgeschlossener) K"orper mit $\textrm{char }K=p>0$, $n \ge 2$ und $G$ eine Untergruppe von $\textrm{GL}_{n}(K)$ mit\\
\begin{enumerate}
\renewcommand{\labelenumi}{(\alph{enumi})}
\item  Falls $p=2$:  
$
\left(
\begin{array}{ccc}
 a & a+1 &\\
a+1&a&\\
&&I_{n-2}
\end{array} \right)
\in G
$
      f"ur wenigstens drei verschiedene Werte von $a\in K$ ($a=1$ ist stets ein solcher).\\
\item Falls $p \ge 3$: 
$
\left(
\begin{array}{ccc}
1  & 1 &\\
0&1&\\
&&I_{n-2}
\end{array} \right)
\in G.
$
\end{enumerate}
Sei weiter $V$ der Raum aller homogenen Polynome vom Grad $p$ in den Variablen $X_{1},...,X_{n}$,
\[
V:=K[X_{1},\ldots,X_{n}]_{p}.
\]
Ist $(a_{ij}) \in G$, so ist durch $(a_{ij}) \cdot X_{j} = \sum_{i=1}^{n}a_{ij}X_{i}$ die kanonische Operation von $G$ auf $V$ gegeben, wodurch $V$ zu einem $G$-Modul wird. Sei 
\[
W:=<X_{1}^{p},\ldots,X_{n}^{p}> \; \, \le \; V 
\]
(Untermodul nach Frobenius) und
\[
U:=\textrm{Hom}_{K}(V,W)_{0}=\{ f \in \textrm{Hom}_{K}(V,W): f|_{W}=0\}.
\]
$U$ ist $G$-Modul nach Satz \ref{Hom0}.\\
Sei ferner $\iota \in \textrm{Hom}_{K}(V,W)$ gegeben durch 
\[
\iota |_{W}=\textrm{id}_{W}
\]
und 
\[
\iota \textrm{ gleich } 0 \textrm{ auf allen Monomen, die nicht in }W \textrm{ liegen.}
\] 
Sei $\tilde{U}:=U \oplus K \iota$ und $\pi: \tilde{U} \rightarrow K$ gegeben durch $\pi(u + \lambda \cdot \iota):=\lambda \quad \textrm{f"ur } u \in U, \lambda \in K$. Dann ist durch
\[
0 \rightarrow U  \rightarrow  \tilde{U} \stackrel{\pi}{\rightarrow} K \rightarrow 0
\]
eine kurze exakte Sequenz gegeben, die nicht zerf"allt.

"Aquivalent dazu ist nach Bemerkung \ref{Kozykl}, dass durch
\[
g_{\sigma}:=(\sigma-1) \iota
\]
ein nichttrivialer Kozyklus in $U$ gegeben ist.

Ebenfalls "aquivalent ist nach Proposition \ref{kompl}, dass der Untermodul $W$ von $V$ kein $G$-invariantes Komplement hat - insbesondere sind also die betrachteten Gruppen nicht linear reduktiv.
\end{Satz}

\Bew Wir verwenden f"ur $V$ eine monomiale Basis $\mathcal{B}$, wobei wir die Reihen\-fol\-ge der ersten $n+1$ bzw. $n+2$ Monome in den F"allen (a) bzw. (b) vorge\-ben, und zwar\\
im Fall (a):
\[
\mathcal{B}=\{X_{1}^{2},...,X_{n}^{2},X_{1}X_{2},...\}
\]\\
im Fall (b):
\[
\mathcal{B}=\{X_{1}^{p},...,X_{n}^{p},X_{1}^{p-1}X_{2},X_{1}^{p-2}X_{2}^{2},...\}
\]\\
Als Basis von $W$ dienen die ersten $n$ Eintr"age von $\mathcal{B}$. Sei $N:=|\mathcal{B}|= {n+p-1 \choose p}$. Wir bezeichnen mit $f_{p}:\textrm{GL}_{n}(K) \rightarrow \textrm{GL}_{n}(K)$ den koeffizientenweisen Frobenius-Homomorphismus, also $f_{p}(a_{ij})=(a_{ij}^{p})$. Ist $A_{\sigma} \in K^{N \times N}$ die Dar\-stell\-ungs\-matrix von $\sigma \in G$ bzgl. der Basis $\mathcal{B}$, so hat diese die Form
\[
A_{\sigma}=\left( \begin{array}{cc}
f_{p}(\sigma) & *\\
0&*\\
\end{array} \right).
\]
Weiter haben wir bzgl. der Basis $\mathcal{B}$
\[
U \cong \left\{
\left( \begin{array}{cc}
0_{n \times n} & B\\
\end{array} \right) \in K^{n \times N} \textrm{  mit  } B \in K^{n \times (N-n)} \right\},
\]
und f"ur $U \ni f \cong   
\left( \begin{array}{cc}
0_{n \times n} & B\\
\end{array} \right)$ bzgl. $\mathcal{B}$ haben wir die Operation gegeben durch
\[
\sigma \cdot f=\sigma \circ f \circ \sigma^{-1} \cong f_{p}(\sigma) \cdot
\left( \begin{array}{cc}
0_{n \times n} & B\\
\end{array} \right)
\cdot A_{\sigma^{-1}}.
\]
Die Darstellungsmatrix von $\iota$ ist gegeben durch
\[
\iota \cong
\left( \begin{array}{cc}
I_{n} & 0\\
\end{array} \right)=:J \in K^{n \times N}.
\]
Da
\[
\sigma \cdot \iota=\sigma \circ \iota \circ \sigma^{-1} \cong f_{p}(\sigma)
\left( \begin{array}{cc}
I_{n} & 0\\
\end{array} \right) 
\underbrace{
\left( \begin{array}{cc}
f_{p}(\sigma^{-1}) & *\\
0&*\\
\end{array} \right)}_{A_{\sigma^{-1}}} 
=
\left( \begin{array}{cc}
I_{n} & *\\
\end{array} \right)
\]
ist $\sigma \iota - \iota \in U$. Damit ist $\pi$ als $G$-Homomorphismus wohldefiniert, und mit $g_{\sigma}:=(\sigma -1)\iota$ ist $g \in Z^{1}(G,U)$ (wir haben dies ohne Matrizen auch schon im Beweis von Proposition \ref{kompl} ein\-gesehen). Wir m"ussen zeigen, dass $g \not\in B^{1}(G,U)$. Wir nehmen das Gegenteil an, also die Existenz eines
\[
U \ni u \cong Z =
\left( \begin{array}{cc}
0_{n \times n} & \hat{Z}\\
\end{array} \right) \in K^{n \times N}, \qquad \hat{Z}=(z_{ij}) \in K^{n \times (N-n)}
\]
mit $g_{\sigma}=(\sigma -1) \iota \stackrel{!}{=}(\sigma -1)u \quad \textrm{f"ur alle } \sigma \in G$ bzw.
\begin{equation} \label{zerf}
f_{p}(\sigma)J A_{\sigma^{-1}}-J \stackrel{!}{=}f_{p}(\sigma)Z A_{\sigma^{-1}}-Z \quad \forall \sigma \in G.
\end{equation}
Diese letzte Gleichung f"uhren wir nun in beiden F"allen zum Widerspruch.\\

\noindent (a) Mit
\[
\sigma:=
\left(
\begin{array}{ccc}
 a & a+1 &\\
a+1&a&\\
&&I_{n-2}
\end{array} \right) =\sigma^{-1}
\]
berechnen wir die $(n+1)$te Spalte von $A_{\sigma^{-1}}$:
\[
\begin{array}{rcl}
\sigma^{-1} \cdot X_{1}X_{2}&=&(aX_{1}+(a+1)X_{2})((a+1)X_{1}+aX_{2})\\
&=&(a^{2}+a)X_{1}^{2}+(a^{2}+a)X_{2}^{2}+X_{1}X_{2}\\
&\cong&(a^{2}+a,a^{2}+a,0_{n-2},1,0)^{T}.\\
\end{array}
\]
Damit vergleichen wir nun auf beiden Seiten von (\ref{zerf}) den Eintrag in der ersten Zeile und $(n+1)$ten Spalte:\\
Links:
\[
(a^{2},a^{2}+1,0_{n-2})
\left( \begin{array}{cc}
I_{n} & 0_{n \times (N-n)}\\
\end{array} \right)
\left( \begin{array}{c}
a^{2}+a\\
a^{2}+a\\
0_{n-2}\\
1\\
0\\
\end{array} \right) =a^{2}+a
\]
Rechts:
\[
\begin{array}{rcl}
(a^{2},a^{2}+1,0_{n-2})
\left( \begin{array}{cc}
0_{n \times n} & \hat{Z}\\
\end{array} \right)
\left( \begin{array}{c}
a^{2}+a\\
a^{2}+a\\
0_{n-2}\\
1\\
0\\
\end{array} \right) -z_{11}&=&a^{2}z_{11}+(a^{2}+1)z_{21}-z_{11}\\
&=&(a^{2}+1)(z_{11}+z_{21})\\
\end{array}
\]
Setzen wir $c:=z_{11}+z_{21}$, so sehen wir, dass $ca^2+c=a^{2}+a$ bzw.
\[
(c+1)a^{2}+a+c=0
\]
f"ur wenigstens drei verschiedene Werte $a \in K$ erf"ullt sein muss. Dies ist ein Widerspruch.\\

\noindent (b) Wir betrachten
\[
\sigma=
\left(
\begin{array}{ccc}
1  & 1 &\\
0&1&\\
&&I_{n-2}
\end{array} \right),
\sigma^{-1}=
\left(
\begin{array}{ccc}
1  & -1 &\\
0&1&\\
&&I_{n-2}
\end{array} \right)
\in G
\]
und berechnen die $(n+1)$te und $(n+2)$te Spalte von $A_{\sigma^{-1}}$:\\
$(n+1)$te Spalte:
\[
\begin{array}{rcl}
\sigma^{-1} \cdot X_{1}^{p-1}X_{2} &=& X_{1}^{p-1}(-X_{1}+X_{2})\\
&=&-X_{1}^{p}+X_{1}^{p-1}X_{2}\\
&\cong &(-1,0_{n-1},1,0)^{T}
\end{array}
\]
$(n+2)$te Spalte:
\[
\begin{array}{rcl}
\sigma^{-1} \cdot X_{1}^{p-2}X_{2}^{2} &=& X_{1}^{p-2}(X_{1}^{2}-2X_{1}X_{2}+X_{2}^{2})\\
&=&X_{1}^{p}-2X_{1}^{p-1}X_{2}+X_{1}^{p-2}X_{2}^{2}\\
&\cong &(1,0_{n-1},-2,1,0)^{T}
\end{array}
\]
Wir vergleichen nun wieder beide Seiten von (\ref{zerf}):\\
(i) \underline{erste Zeile, $(n+1)$te Spalte}\\
Links:
\[
\left( \begin{array}{ccc} 1 & 1 & 0_{n-2}\\ \end{array} \right)
\left( \begin{array}{cc} I_{n} & 0_{n \times (N-n)}\\ \end{array} \right)
\left( \begin{array}{c}
-1\\
0_{n-1}\\
1\\
0\\
\end{array} \right)=-1
\]
Rechts:
\[
\left( \begin{array}{ccc} 1 & 1 & 0_{n-2}\\ \end{array} \right)
\left( \begin{array}{cc} 0_{n \times n} & \hat{Z}\\ \end{array} \right)
\left( \begin{array}{c}
-1\\
0_{n-1}\\
1\\
0\\
\end{array} \right)-z_{11}=z_{11}+z_{21}-z_{11}=z_{21}
\]
Da Gleichheit gelten soll, haben wir
\begin{equation} \label{wid1}
z_{21}=-1.
\end{equation}

\noindent (ii) \underline{zweite Zeile, $(n+2)$te Spalte}\\
Links:
\[
\left( \begin{array}{ccc} 0 & 1 & 0_{n-2}\\ \end{array} \right)
\left( \begin{array}{cc} I_{n} & 0_{n \times (N-n)}\\ \end{array} \right)
\left( \begin{array}{c}
1\\
0_{n-1}\\
-2\\
1\\
0\\
\end{array} \right)=0
\]
Rechts:
\[
\left( \begin{array}{ccc} 0 & 1 & 0_{n-2}\\ \end{array} \right)
\left( \begin{array}{cc} 0_{n \times n} & \hat{Z}\\ \end{array} \right)
\left( \begin{array}{c}
1\\
0_{n-1}\\
-2\\
1\\
0\\
\end{array} \right)-z_{22}=-2z_{21}+z_{22}-z_{22}=-2z_{21}
\]
Da $p \ge 3$ ist $2 \ne 0$, und der Vergleich beider Seiten liefert
\[
z_{21}=0,
\]
im Widerspruch zu (\ref{wid1}). \qed\\

\noindent Mit Korollar \ref{NCM} ist dann also mit obiger Notation mit 
\[
\begin{array}{crcl}
&X&:=&U^{*} \oplus \tilde{U} \oplus \tilde{U} \oplus \tilde{U}\\
\textrm{bzw. }& X^{*}&:=&U \oplus \tilde{U}^{*} \oplus \tilde{U}^{*} \oplus \tilde{U}^{*}
\end{array}
\]
durch $K[X]^{G}$ ein nicht Cohen-Macaulay Invariantenring gegeben, falls $G$ reduktiv ist. Insbesondere haben wir also Beispiele f"ur die Gruppen $\textrm{SL}_{n}(K)$ und $\textrm{GL}_{n}(K)$ mit $p>0$ und $n \ge 2$ beliebig. Wir interessieren uns ab jetzt nur noch f"ur den Fall $n=2$. F"ur die Dimension des Moduls gilt dann offenbar $\dim X=4\cdot\dim U+3=4 \cdot \left(2 (p+1-2) \right) +3$, d.h.
\[
\dim X=8p-5,
\] 
also $11$ f"ur $p=2$ und $19$ f"ur $p=3$. Diese Dimensionen wollen wir reduzieren.

\newpage
\section{Beispiele f"ur \SL2 und \GL2 sowie \SO}
Wir wollen nun Beispiele m"oglichst kleiner Dimension konstruieren, wobei wir das Konstruktionsverfahren aus Korollar \ref{NCM} leicht ab"andern. Nach dem Hauptsatz, den wir stets f"ur die Konstruktion verwenden, ben"otigen wir einen nichttrivialen Kozyklus, der sich in einem Modul $U$ befindet, und drei Annullatoren, die sich jeweils in Kopien von $\tilde{U}^{*}$ befinden. Damit sich diese Zutaten in $S(V^{*})$ wiederfinden, haben wir $V^{*}$ in \ref{NCM} einfach als die direkte Summe all dieser Zutaten definiert. F"ur den Beweis von \ref{NCM} reicht es jedoch aus, wenn $S(V^{*})$ diese Zutaten enth"alt, genau genommen brauchen wir Folgendes:
\begin{itemize}
\item Der Modul mit dem nichttrivialen Kozyklus $U$ ist ein Untermodul von $S(V^{*})=K[V]$ und hat ein $G$-invariantes Komplement $U'$, also $S(V^{*})=U \oplus U'$. Dieses Komplement wird gebraucht, damit der Kozyklus auch in $S(V^{*})$ nichttrivial bleibt (vgl. den Beweis von Korollar \ref{NCM}).
\item Es gibt drei (nicht notwendig disjunkte) Untermoduln $W_{1},W_{2},W_{3}$ von $S(V^{*})$, die isomorph zu $\tilde{U}^{*}$ sind, und deren Kopien $\pi_{1},\pi_{2},\pi_{3}$ des Annul\-lators $\pi \in \tilde{U}^{*}$ bilden ein phsop in $K[V]$.
\end{itemize}
In diesem Fall l"asst sich n"amlich der Hauptsatz genau wir im Beweis von Korollar \ref{NCM} anwenden, d.h. das phsop $\pi_{1},\pi_{2},\pi_{3}$ ist keine regul"are Sequenz in $K[V]^{G}$, und $K[V]^{G}=S(V^{*})^{G}$ ist damit nicht Cohen-Macaulay.\\

\noindent Die einfachste M"oglichkeit, diese beiden Punkte umzusetzen, gibt sich, wenn sich die Moduln $U$ oder $\tilde{U}^{*}$ als Tensorprodukt schreiben lassen, denn das Tensorprodukt zweier Moduln ist Untermodul (mit Komplement) der zweiten Potenz der direkten Summe. Dieses Verfahren werden wir in den folgenden Abschnitten durchf"uhren.

\subsection{Kozyklen im Grad 2}
Wir betrachten den Modul mit dem Kozyklus, $U=\textrm{Hom}_{K}(V,W)_{0}$ in Satz \ref{PKoz} und der Dimension $\dim W \cdot (\dim V - \dim W)=2(p+1-2)=2(p-1)$ (f"ur $n=2$).
Nach Satz \ref{Hom0} gilt nun aber
\[
\textrm{Hom}_{K}(V,W)_{0} \cong W \otimes (V/W)^{*} \le S^{2}\left(W \oplus (V/W)^{*} \right),
\]
und wir k"onnen $U$ in (\ref{VStern})  in Korollar \ref{NCM} durch diese direkte Summe der Dimension $2+(p+1-2)=p+1$ ersetzen - mit der Notation im Text nach dem Beweis von Satz \ref{PKoz} also $X^{*}$ durch
\[
X^{*}:=W \oplus (V/W)^{*} \oplus \tilde{U}^{*} \oplus \tilde{U}^{*}
\oplus \tilde{U}^{*}
\]
(und $K[X]^{G}$ dann nicht Cohen-Macaulay).
Dies ergibt eine Dimensions\-erspar\-nis von $2(p-1)-(p+1)=p-3$, oder wie schon angek"undigt die Dimensions\-reduk\-tion von $8p-5$ auf $(8p-5)-(p-3)=7p-2$. Ab $p \ge 5$ handelt es sich hier um eine echte Ersparnis. F"ur diese Primzahlen konnte ich ansonsten keine weiteren Beispiele finden. Wir notieren unser Ergebnis:

\begin{Satz} \label{p72}
F"ur einen algebraisch abgeschlossenen K"orper $K$ der Charakteristik $p>0$  und einer reduktiven Zwischengruppe $G$ von \SL2 und \GL2 gibt es einen $G$-Modul $V$ der Dimension $7p-2$ mit nicht Cohen-Macaulay Invariantenring $K[V]^{G}$. Dabei liegt der nichttriviale Kozyklus im Grad $2$ und das annullierende phsop im Grad 1. \qed
\end{Satz}

\subsection{Beispiele f"ur \SL2 in Charakteristik 2}
Der Ausgangspunkt f"ur die Konstruktion ist meist ein Modul mit nicht\-tri\-vialem Kozyklus. Von diesem ausgehend versucht man, m"oglichst kleine Moduln mit Annullatoren dazu zu addieren. Daher ist dieser Abschnitt nach den verwendeten Kozyklen gegliedert.
\subsubsection{Beispiele mit dem Kozyklus in $\langle X^{2},Y^{2} \rangle$}
Wir wollen zun"achst die in Satz \ref{PKoz} und im vorigen Abschnitt auftretenden Moduln f"ur den Fall $p=n=2$ etwas expliziter angeben. Setzen wir $X:=X_{1}, Y:=Y_{2}$ so ist mit den Bezeichnungen aus \ref{PKoz}
\[
\begin{array}{rcl}
V&=&\langle X^{2},Y^{2},XY \rangle\\
W&=&\langle X^{2},Y^{2} \rangle.
\end{array}
\]
Die Darstellung des Moduls $U=\textrm{Hom}_{K}(V,W)_{0}$, welcher durch $2 \times 3$ Matrizen mit Eintr"agen nur in der letzten Spalte beschrieben wird, erh"alt man aus
\[
\left( \begin{array}{cc}
a^{2} & b^{2}\\
c^{2} &  d^{2}
\end{array} \right)
\left( \begin{array}{ccc}
0 & 0 & x_{1}\\
0 & 0 & x_{2}
\end{array} \right)
\left( \begin{array}{ccc}
a^{2} & b^{2} & ab\\
c^{2} & d^{2} & cd\\
0 & 0 & 1
\end{array} \right)^{-1}
\]
zu
\[
\sigma=
\left( \begin{array}{cc}
a & b\\
c &  d
\end{array} \right) \mapsto
\left( \begin{array}{cc}
a^{2} & b^{2}\\
c^{2} &  d^{2}
\end{array} \right),
\] 
denn von der rechten Matrix wirkt blos die letzte Zeile auf die Koeffizienten $x_{1},x_{2}$. Damit ist also 
\[
U= \langle X^{2},Y^{2} \rangle,
\]
vgl. Tabelle \ref{Darstellungen} (S. \pageref{Darstellungen}). F"ur $\tilde{U}$ m"ussen wir noch die Operation auf $\iota \cong \left( \begin{array}{ccc} 1&0&0\\ 0&1 &0 \end{array} \right)$ bestimmen. Diese erhalten wir mit
$
\left( \begin{array}{cc}
a & b\\
c &  d
\end{array} \right)^{-1}
=
\left( \begin{array}{cc}
d & b\\
c &  a
\end{array} \right)
$ aus
\[
\begin{array}{cl}
&\left( \begin{array}{cc}
a^{2} & b^{2}\\
c^{2} &  d^{2}
\end{array} \right)
\left( \begin{array}{ccc}
1 & 0 & 0\\
0 & 1 & 0\\
\end{array} \right)
\left( \begin{array}{ccc}
a^{2} & b^{2} & ab\\
c^{2} & d^{2} & cd\\
0 & 0 & 1
\end{array} \right)^{-1}\\
=&
\left( \begin{array}{ccc}
a^{2} & b^{2} & 0\\
c^{2} & d^{2} & 0\\
\end{array} \right)
\left( \begin{array}{ccc}
d^{2} & b^{2} & bd\\
c^{2} & a^{2} & ac\\
0 & 0 & 1
\end{array} \right)\\
=&
\left( \begin{array}{ccc}
1 & 0 & ab(ad+bc)\\
0 & 1 & cd(ad+bc)
\end{array} \right)
\end{array}
\]
und $ad+bc=1$ zu $\sigma \iota=abX^{2}+cdY^{2}+\iota$. Wieder anhand von Tabelle \ref{Darstellungen} erkennen wir 
\[
\tilde{U}=\langle X^{2},Y^{2},XY \rangle,
\]
und der nichttriviale Kozyklus in $U$ ist gegeben durch
\[
g_{\sigma}=(\sigma-1)XY \cong (\sigma-1)\iota.
\]
Weiter ist dann
\[
\tilde{U}^{*} =\langle \mu,\nu,\pi \rangle
\]
(siehe wieder Tabelle \ref{Darstellungen}), und $\pi$ ist nach Proposition \ref{anni} die annullierende Invariante.\\ 

\noindent Wie steht es mit der im vorigen Abschnitt genannten Darstellung von $U$ als Tensorprodukt? Anhand der Darstellungen von $V$ und $W$ sieht man sofort, dass hier $V/W \cong K$ der triviale Modul ist. Also ist
\[
U=\textrm{Hom}_{K}(V,W)_{0} \cong W \otimes (V/W)^{*} \cong W \otimes K \cong W.
\]
Zum einen haben wir hier nochmals ohne Rechnung $U=W$ gesehen, zum anderen sehen wir, dass es hier keinen Sinn macht $U$ durch die direkte Summe $U  \oplus K$ zu ersetzen, weil dies offenbar keine Dimensionsersparnis gibt.\\
Insgesamt erhalten wir jedenfalls mit Korollar \ref{NCM} 

\begin{samepage}
\begin{Bsp}[Das Beispiel nach Satz \ref{PKoz}]
F"ur Charakteristik p=2 ist mit
\[
\begin{array}{cl}
&V^{*} = \langle X^2,Y^2 \rangle \bigoplus_{i=1}^{3} \langle \mu, \nu, \pi 
\rangle\\
\textrm{bzw.} & V=\langle X^2,Y^2 \rangle \bigoplus_{i=1}^{3} \langle X^{2}, Y^{2}, XY \rangle
\end{array}
\]
der Invariantenring \INVSL2 nicht Cohen-Macaulay, $\dim V=11$. Der Kozyklus im Grad $1$ liegt in $\langle X^2,Y^2 \rangle$, und die drei Kopien von $\pi$ bilden ein annullierendes phsop im Grad 1, welches nach dem Hauptsatz keine regul"are Sequenz ist. \qed
\end{Bsp}
\end{samepage}

Dieses Beispiel kann nat"urlich auch mit {\tt IsNotCohenMacaulay} untersucht werden (Datei {\tt Beispiel-6.02.txt}.) Nach der Absch"atzung (\ref{dschaetz})  (S. \pageref{dschaetz}) gen"ugt es, {\tt dmax=1+1+1+1=4} zu setzen. Der Aufruf von 
{\tt IsNotCohenMacaulay} zeigt nach wenigen Sekunden tats"achlich, dass die Hilbertreihe bis zum grad $3$ richtig gesch"atzt wird, aber im Grad $4$ werden $33$ Invarianten gesch"atzt, w"ahrend es nur $32$ sind (insbesondere wird $K[V]^{G}$ also wie nach Satz \ref{Erkenn} erwartet als nicht Cohen-Macaulay erkannt). Dies zeigt, dass die im Hauptsatz angegebene nichttriviale Relation in diesem Fall auch (im wesentlichen) die einzige (bis zum Grad 4) ist.\\

\noindent Wir wollen dieses Beispiel modifizieren. Zun"achst suchen wir den Modul $\langle \mu,\nu,\pi \rangle$ mit dem Annullator des Kozyklus als Untermodul eines einfacheren Moduls. Wir betrachten dazu das Tensorprodukt von $\langle X,Y \rangle$ mit sich selbst, welches nach Gleichung (\ref{XXYYbasis}) auf S. \pageref{XXYYbasis} bzgl. der Basis
\[
\mathcal{B}=\left( X \otimes X, Y \otimes Y, X\otimes Y-Y\otimes X,Y\otimes X \right)
\]
die Darstellung
\begin{equation} \label{TensBasis}
\sigma \mapsto \left( \begin{array}{cccc}
a^{2} & b^{2} & 0 & ab\\
c^{2} & d^{2} & 0 & cd\\
ac & bd & 1 & bc\\
0 & 0 & 0 & 1\\
\end{array} \right)
\end{equation}
besitzt (f"ur die Gruppe \SL2 mit $p=2$) - insbesondere erkennen wir $X\otimes Y-Y\otimes X$ als eine Invariante.
In Tabelle \ref{Darstellungen} (S. \pageref{Darstellungen}) finden wir die Darstellung von $\langle \mu,\nu,\pi \rangle$ bzgl. der angedeuteten Basis. Vertauschen der ersten beiden Zeilen und Spalten liefert die Darstellung bzgl. der Basis $(\nu,\mu,\pi)$, n"amlich
\[
\sigma \mapsto \left( \begin{array}{ccc}
a^{2} & b^{2} & 0\\
c^{2} & d^{2} & 0\\
ac & bd & 1
\end{array} \right).
\]
Damit ist $\langle \mu,\nu,\pi \rangle$ Untermodul des Tensorprodukts $\langle X\otimes X, \ldots, Y\otimes Y \rangle$. Die Invariante $\pi$ entspricht dabei der Invarianten $X\otimes Y-Y\otimes X$, d.h. dies ist der neue Annullator. Wie wir in Abschnitt \ref{BspMod} gesehen haben, ist das Tensorprodukt selbstdual. Damit haben wir also

\begin{Bsp}
Mit
\[
V^{*} = \langle X^2,Y^2 \rangle \bigoplus_{i=1}^{3} \langle X \otimes X,X \otimes Y, Y\otimes X,Y\otimes Y \rangle
\]
ist der der Invariantenring \INVSL2 nicht Cohen-Macaulay. $V$ ist selbstdual. Es ist $\dim V=14$ und phsop und Kozyklus liegen im Grad 1, wobei die annullierenden Invarianten die drei Kopien von $X \otimes Y+Y\otimes X$ sind. \qed
\end{Bsp}

\noindent (siehe auch Datei {\tt Beispiel-6.03.txt}).\\
Man kann nun nat"urlich hier ein bisschen rumspielen, und ein Tensorprodukt durch den Modul $\langle \mu,\nu,\pi \rangle$ austauschen, so dass man verschiedene, jedoch  sehr "ahnliche Beispiele hat (f"ur ein Beispiel dieser Art siehe die Bemerkungen am Ende dieses Abschnitts). Interessanter ist jedoch folgende "Uberlegung: Das Tensorprodukt ist ja Unter\-modul der zweiten symmetrischen Potenz $S^{2} \left( \langle X,Y \rangle \oplus \langle X,Y \rangle \right)$. Wir k"onnen das Tensorprodukt also durch $\langle X,Y \rangle \oplus \langle X,Y \rangle$ ersetzen, bekommen dann also annullierende Invarianten im Grad $2$. Wir m"ussen nun aber nicht jedes der drei Tensorprodukte durch solch eine direkte Summe ersetzen, denn es reicht ja, wenn wir ein phsop an annullierenden Invarianten haben. Die zweite Potenz von $\langle X,Y \rangle \oplus \langle X,Y \rangle \oplus \langle X,Y \rangle$  enth"alt z.B. schon drei Kopien des Tensor\-produkts (jeweils f"ur die zweite Potenz von zwei Summanden), aber leider bilden die Invarianten hier kein phsop, wie man sich z.B. leicht mit \Magma "uberzeugt. Erst bei vier Summanden von $\langle X,Y \rangle$ (wo die zweite Potenz dann $6$ Kopien des Tensorpro\-dukts enth"alt), bilden drei der Invarianten der Form $X_{1}Y_{2}+X_{2}Y_{1}$ ein annul\-lier\-en\-des phsop (wie man wieder mit \Magma best"atigen kann). Aus dem Hauptsatz folgt dann eines der Hauptresultate dieser Arbeit:

\begin{Bsp}[Das Hauptbeispiel in Dimension 10]
Mit
\[
V^{*} = \langle X^2,Y^2 \rangle \bigoplus_{i=1}^{4} \langle X_{i},Y_{i} \rangle 
\]
ist der der Invariantenring \INVSL2 nicht Cohen-Macaulay. $V$ ist selbst\-dual. Es ist $\dim V=10$ und der Kozyklus liegt im Grad 1 im Untermodul $\langle X^2,Y^2 \rangle$. Ein phsop an annul\-lierenden Invarianten im Grad $2$ ist bei\-spiels\-weise gegeben durch $f_{12},f_{23},f_{34}$ mit $f_{ij}=X_{i}Y_{j}+X_{j}Y_{i}$, und dieses ist nach dem Hauptsatz nicht regul"ar. \qed
\end{Bsp}

Wenn man dieses Beispiel mit {\tt IsNotCohenMacaulay} untersucht (Datei {\tt Beispiel-6.04.txt}), muss man nach (\ref{dschaetz}) mindestens {\tt dmax:=1+2+2+2=7} setzen, denn das phsop liegt nun im Grad $2$. Man muss also Invarianten bis zum Grad $7$ berechnen! Mit dem nicht modifizierten Bayer-Algorithmus dauert das sehr lange (\emph{IsNotCohenMacaulay} rechnet dann etwa $3$ Stunden). "Ubergibt man jedoch den Gewichtsvektor $w$ und berechnet also vorab die Invarianten des Torus, so erkennt \emph{IsNotCohenMacaulay} dieses Beispiel in etwa einer halben Minute als nicht Cohen-Macaulay. Bis zum Grad $6$ stimmt dann die gesch"atzte mit der tats"achlichen Hilbertreihe "uberein, f"ur den Grad $7$ dagegen werden 77 Basis-Invarianten gesch"atzt, wobei es nur 76 sind. Wieder ein Indiz, dass die im Hauptsatz angegebene Relation bis zum unter\-suchten Grad die einzige ist.\\

\noindent Hier noch eine Interpretation der verwendeten Invariante als Determinante: Wir betrachten die zweifache Summe des nat"urlichen Moduls, $\langle X_{1},Y_{1} \rangle \oplus \langle X_{2},Y_{2} \rangle$, und bilden "uber dem Polynomring $K[X_{1},Y_{1},X_{2},Y_{2}]$ die $2 \times 2$-Matrix $\left( \begin{array}{cc} X_{1}&X_{2}\\ Y_{1} & Y_{2} \end{array} \right)$ mit Determinante $X_{1}Y_{2}-X_{2}Y_{1}$. Linksmultiplikation dieser Matrix mit $\sigma^{T} = \left( \begin{array}{cc} a&c\\ b & d \end{array} \right) \in \textrm{SL}_{2}(K)$ liefert die Matrix $\left( \begin{array}{cc} \sigma X_{1}&\sigma X_{2}\\ \sigma Y_{1} & \sigma Y_{2} \end{array} \right)$ mit Determinante $(\sigma X_{1})(\sigma Y_{2})-(\sigma X_{2})(\sigma Y_{1})=\sigma(X_{1}Y_{2}-X_{2}Y_{1})$. Da jedoch $\det \sigma =1$, hat sich nach dem Deter\-minanten\-multi\-pli\-kations\-satz die Determinan\-te nicht ver"andert, d.h.  $X_{1}Y_{2}-X_{2}Y_{1}$ ist invariant.\\

\noindent Wir geben noch eine der bereits erw"ahnten Spielereien an: Mit
\[
V^{*} = \langle X^2,Y^2 \rangle \bigoplus_{i=1}^{3} \langle X,Y \rangle \oplus \langle \mu, \nu, \pi  \rangle\\
\]
ist der der Invariantenring \INVSL2 nicht Cohen-Macaulay - zwei der annul\-lierenden Invarianten liegen in zweiter Potenz und eine in der Ersten. Die Dimension dieses Beispiels ist $11$.\\

\subsubsection{Beispiele mit dem Kozyklus in $\langle \mu, \nu, \pi \rangle$} \label{munupi}
Wir konstruieren aus dem zu Bemerkung \ref{InvKoz} geh"origen Verfahren einen Kozyklus aus der Invarianten $X\otimes Y-Y\otimes X$ des Tensorprodukts. Zun"achst sei die Charakteristik $p$ hier noch beliebig. Wir beschreiben das Tensorprodukt $\langle X,Y \rangle \otimes \langle -Y,X \rangle$ gem"a\ss{} Lemma \ref{Tensor} durch Matrizen, wobei wir f"ur den ersten Faktor die Basis $(X,Y)$ und f"ur den zweiten die Basis $(-Y,X)$ (entspricht der Dualbasis $(X^{*},Y^{*})$) und damit die Darstellung
\[
\sigma \mapsto
\left( \begin{array}{cc}
d&-c\\
-b&a
\end{array} \right)
=
\left( \begin{array}{cc}
a&b\\
c&d
\end{array} \right)^{-T}
\]
verwenden wollen. Auf den Koordinatenmatrizen des Tensorprodukts haben wir dann nach Lemma \ref{Tensor} die Operation
\[
\begin{array}{rcl}
\sigma \cdot X &=&
\left( \begin{array}{cc}
a&b\\
c&d
\end{array} \right)
\left( \begin{array}{cc}
x_{11}&x_{12}\\
x_{21}&x_{22}
\end{array} \right)
\left( \begin{array}{cc}
a&b\\
c&d
\end{array} \right)^{-TT}\\
&=&
\left( \begin{array}{cc}
a&b\\
c&d
\end{array} \right)
\left( \begin{array}{cc}
x_{11}& x_{12}\\
x_{21}& x_{22}
\end{array} \right)
\left( \begin{array}{cc}
a&b\\
c&d
\end{array} \right)^{-1}.
\end{array}
\]
Der Invarianten $X\otimes Y-Y\otimes X=-X\otimes(-Y) - Y \otimes X$ entspricht die Koordinatenmatrix
$
\left( \begin{array}{cc}
-1&0\\
0&-1
\end{array} \right)=-I_{2}
$. Da die Einheitsmatrix unter Kon\-jugation mit Matrizen invariant ist, sieht man hier auch nochmal sch"on, dass $X\otimes Y-Y\otimes X$ invariant ist. Die Spur einer Matrix $X$ ist unter Konjugation ebenfalls invariant, daher bilden die Matrizen mit Spur $0$ einen \SL2 - invarianten Unterraum der Kodimension $1$. Im Falle $p \ne 2$ ist $\textrm{Spur}(-I_{2})=-2\ne 0$, und die spurlosen Matrizen bilden ein \SL2-invariantes Komplement zu $K\cdot I_{2}$. Dies bedeutet nach Bemerkung \ref{InvKoz}, dass der aus der Invarianten $X\otimes Y-Y\otimes X$ konstruierte Kozyklus trivial ist. Im Fall $p=2$ hat die zu dieser Invarianten geh"orende Koordinatenmatrix ebenfalls Spur $0$ - es gibt also eine Chance, dass sie kein Komplement hat. Wir zeigen nun einfach direkt, dass der entstehende Kozyklus nichttrivial ist: Wie wir gesehen haben, hat das Tensorprodukt 
\[
\tilde{V}^{*}:=\left( X \otimes X, Y \otimes Y, X\otimes Y-Y\otimes X,Y\otimes X \right)
\]
bzgl. der angedeuteten Basis die Darstellung
\[
\sigma \mapsto \left( \begin{array}{cccc}
a^{2} & b^{2} & 0 & ab\\
c^{2} & d^{2} & 0 & cd\\
ac & bd & 1 & bc\\
0 & 0 & 0 & 1\\
\end{array} \right).
\]
Damit die Invariante das letzte Basiselement ist, vertauschen wir die letzten beiden Zeilen und Spalten und erhalten
\begin{equation} \label{TensInv}
\left( \begin{array}{cccc}
a^{2} & b^{2} & ab & 0\\
c^{2} & d^{2} & cd & 0\\
0 &0&1&0\\
ac & bd & bc & 1
\end{array} \right).
\end{equation}
Wir gehen zum Dual $\tilde{V}$ "uber, indem wir invertieren (d.h. an der Stelle $\left( \begin{array}{cc}
d&b\\
c&a
\end{array} \right)$ auswerten) und transponieren; das ergibt
\begin{equation} \label{TensInv0}
\left( \begin{array}{cccc}
d^{2} & c^{2} & 0 & cd\\
b^{2} & a^{2} & 0 & ab\\
bd & ac & 1 & bc\\
0&0&0&1
\end{array} \right).
\end{equation}
Vertauschen der ersten beiden Zeilen und Spalten liefert wieder die Ausgangs\-matrix. Die letzte Spalte beschreibt dann jedenfalls einen nichttrivialen Ko\-zyklus im Untermodul $V=\langle \mu,\nu,\pi \rangle$ des (selbstdualen) Tensorprodukts: Sonst g"abe es n"amlich $x=(x_{1},x_{2},x_{3})^{T} \in K^{3}$ mit
\[
\left( \left( \begin{array}{ccc}
a^{2} & b^{2} & 0 \\
c^{2} & d^{2} & 0 \\
ac & bd & 1 
\end{array} \right)-I_{3} \right) 
\left(\begin{array}{c}
x_{1} \\ x_{2} \\ x_{3}
\end{array} \right)
=
\left( \begin{array}{c}
ab\\
cd\\
bc
\end{array} \right) \quad \forall
\left( \begin{array}{cc}
a&b\\
c&d
\end{array} \right) \in \textrm{SL}_{2}(K).
\]
Betrachtet man hiervon nur die ersten beiden Zeilen, so w"are auch der Ko\-zyklus $\sigma \mapsto (\sigma-1)XY$  in $U:=\langle X^{2},Y^{2} \rangle$ trivial - siehe die Darstellung von $\tilde{U}=\langle X^{2},Y^{2},XY \rangle$ in Tabelle \ref{Darstellungen}. Wie wir im vorigen Abschnitt gesehen haben, ist dies jedoch nicht der Fall, also ein Widerspruch.\\

\noindent Es gibt noch zwei andere M"oglichkeiten einzusehen, dass der betrachtete Kozyklus nichttrivial ist (die aber letztendlich auf dieselbe Rechnung hinaus\-f"uhren):

Zum einen kann man den Kozyklus auch im Faktormodul $\langle X^{2},Y^{2} \rangle= \langle \mu, \nu, \pi \rangle / \langle \pi \rangle$ von $\langle \mu,\nu,\pi \rangle$ betrachten. Da er hier nichttrivial ist (wieder auf\-grund des vorigen Abschnitts), ist er es erst recht in $\langle \mu,\nu,\pi \rangle$.

Oder man stellt fest, dass der betrachtete Kozyklus gem"a\ss{} Lemma \ref{InvKoz} zu der Invarianten $\pi \cong X\otimes Y -Y \otimes X$ im Tensorprodukt geh"ort. Diese hat jedoch bereits im Untermodul $\langle \mu,\nu,\pi \rangle$ des Tensorprodukts kein Komplement, denn der der Invarianten $\pi$ in $\langle \mu,\nu,\pi \rangle$ nach Lemma \ref{InvKoz} zugeh"orige Kozyklus in $\langle X^{2},Y^{2} \rangle$ ist ja (wieder nach vorigem Abschnitt) nichttrivial. Also hat $\pi$ nach Bemerkung \ref{kleinKompl} erst recht kein Komplement im gr"o\ss eren Modul, dem Tensorprodukt, und wieder nach Lemma \ref{InvKoz} ist dann der zugeh"orige Kozyklus in $\langle \mu,\nu,\pi \rangle$ nichttrivial.\\

{\noindent \bf Zusammengefasst: }Der Modul $V=\langle \mu,\nu,\pi \rangle$ enth"alt einen nichttrivialen Kozyklus, der $V$ erweitert zum Tensorprodukt $\tilde{V}=\langle X\otimes X,\ldots, Y\otimes Y \rangle$. Der Annullator im Dual $\tilde{V}^{*}=\tilde{V}$ ist gegeben durch $X\otimes Y-Y \otimes X$.\\

\noindent Aus Korollar \ref{NCM} erhalten wir also sofort

\begin{Bsp}
F"ur Charakteristik p=2 ist mit
\[
\begin{array}{cl}
&V^{*} = \langle \mu, \nu, \pi 
 \rangle \bigoplus_{i=1}^{3} \langle X\otimes X,X\otimes Y,Y\otimes X, Y\otimes Y \rangle\\
\textrm{bzw.} & V=\langle X^2,Y^2,XY \rangle \bigoplus_{i=1}^{3} \langle X\otimes X,X\otimes Y,Y\otimes X, Y\otimes Y \rangle
\end{array}
\]
der Invariantenring \INVSL2 nicht Cohen-Macaulay, wobei $\dim V=15$. 
Der Kozyklus liegt im Grad 1 im Untermodul $\langle \mu, \nu, \pi  \rangle$,
und ein phsop an annullierenden Invarianten im Grad $1$, (welches keine regul"are Sequenz ist), bilden die drei Kopien der $X\otimes Y-Y\otimes X$.  \qed
\end{Bsp}

F"ur \emph{IsNotCohenMacaulay} (siehe Datei {\tt Beispiel-6.05.txt}) muss man (nach (\ref{dschaetz} auf S. \pageref{dschaetz}) also wieder {\tt dmax:=4} setzen. Der Koeffizient der Hilbertreihe zum Grad $4$ wird gesch"atzt zu $143$, w"ahrend er nur $142$ ist.\\

\noindent Mit der Diskussion im vorigen Abschnitt, wo wir das Tensorprodukt in einer zweiten Potenz wiederfanden, erhalten wir aus dem Hauptsatz
\begin{Bsp}
F"ur Charakteristik p=2 ist mit
\[
\begin{array}{cl}
&V^{*} = \langle \mu, \nu, \pi 
 \rangle \bigoplus_{i=1}^{4} \langle X_{i},Y_{i} \rangle\\
\textrm{bzw.} & V=\langle X^2,Y^2,XY \rangle \bigoplus_{i=1}^{4} \langle X_{i},Y_{i} \rangle
\end{array}
\]
der Invariantenring \INVSL2 nicht Cohen-Macaulay. Es ist $\dim V=11$ und der Kozyklus liegt im Untermodul $\langle \mu, \nu, \pi 
 \rangle$ im Grad 1. Ein phsop an annul\-lierenden Invarianten im Grad $2$, welches keine regul"are Sequenz ist, ist beispielsweise gegeben durch $f_{12},f_{23},f_{34}$ mit $f_{ij}=X_{i}Y_{j}+X_{j}Y_{i}$. \qed
\end{Bsp}

Da das phsop im Grad $2$ liegt, muss man bei \emph{IsNotCohenMacaulay} (siehe Datei {\tt Beispiel-6.06.txt}) nun also {\tt dmax:=7} setzen. Auch hier f"uhrt die Angabe des Gewichtsvektors $w$ wieder zu einer enormen Beschleunigung. Der Koeffizient der Hilbertreihe zum Grad $7$ wird dann gesch"atzt zu 193, w"ahrend er nur 192 ist.\\

\subsubsection{Kozyklus und Invariante in einem}
Wir wollen noch einmal das Verfahren nach Bemerkung \ref{InvKoz} anwenden, aus\-ge\-hend von einer Invarianten einen Modul mit nichttrivialem Kozyklus zu konstruieren. Wir betrachten die zweite symmetrische Potenz $\tilde{U}^{*}$ des Moduls $\langle \mu,\nu,\pi \rangle$, wobei $\pi^{2}$ die annullierende Invariante werden soll. Ausgehend von der Dar\-stel\-lung von $\langle \mu,\nu,\pi \rangle$,
\[
\sigma \mapsto
\left( \begin{array}{ccc}
d^{2} & c^{2} & 0\\
b^{2} & a^{2} & 0\\
bd & ac & 1
\end{array} \right)
\]
berechnen wir die Darstellung von $\langle \mu\nu,\mu^{2},\nu^{2},\mu\pi,\nu\pi,\pi^{2} \rangle$:
\[
\begin{array}{rcl}
\sigma( \mu\nu)&=&(d^{2} \mu+b^{2} \nu + bd \pi)(c^{2} \mu+a^{2}\nu+ac\pi)\\
&=&\mu\nu+c^{2}d^{2}\mu^{2}+a^{2}b^{2}\nu^{2}+(acd^{2}+bc^{2}d)\mu\pi+(acb^{2}+a^{2}bd)\nu\pi+abcd\pi^{2}\\
&\cong&(1,c^{2}d^{2},a^{2}b^{2},cd,ab,abcd)^{T}\\
&&\\
\sigma(\mu^{2})&=&(d^{2}\mu+b^{2}\nu+bd\pi)^{2}\\
&=&d^{4}\mu^{2}+b^{4}\nu^{2}+b^{2}d^{2}\pi^{2}\\
&\cong&(0,d^{4},b^{4},0,0,b^{2}d^{2})^{T}\\
&&\\
\sigma(\nu^{2})&=&(c^{2}\mu+a^{2}\nu+ac\pi)^{2}\\
&=&c^{4}\mu^{2}+a^{4}\nu^{2}+a^{2}c^{2}\pi^{2}\\
&\cong&(0,c^{4},a^{4},0,0,a^{2}c^{2})^{T}
\end{array}
\]
Die Berechnung der letzten drei Spalten ist einfach, denn $\langle \mu,\nu,\pi \rangle \cong \langle \mu\pi,\nu\pi,\pi^{2} \rangle$. Also hat $\tilde{U}^{*}=\langle \mu\nu,\mu^{2},\nu^{2},\mu\pi,\nu\pi,\pi^{2} \rangle$ die Darstellung
\[
\sigma \mapsto
\left( \begin{array}{cccccc}
1 &0&0&0&0&0\\
c^{2}d^{2} & d^{4} & c^{4} &0&0&0\\
a^{2}b^{2} &b^{4} & a^{4} &0&0&0\\
cd &0&0 &d^{2} &c^{2} &0\\
ab&0&0&b^{2} &a^{2}&0\\
abcd&b^{2}d^{2}&a^{2}c^{2} &bd&ac&1
\end{array} \right).
\]
Durch Invertieren (auswerten an $\left( \begin{array}{cc}
d&b\\
c&a
\end{array} \right)$) und Transponieren gehen wir zum Dual $\tilde{U}$ "uber, mit der Darstellung
\[
\sigma \mapsto
\left( \begin{array}{cccccc}
1&a^{2}c^{2}&b^{2}d^{2}&ac&bd&abcd\\
0&a^{4}&b^{4}&0&0&a^{2}b^{2}\\
0&c^{4}&d^{4}&0&0&c^{2}d^{2}\\
0&0&0&a^{2}&b^{2}&ab\\
0&0&0&c^{2}&d^{2}&cd\\
0&0&0&0&0&1
\end{array} \right).
\]
Die letzte Spalte geh"ort dabei wieder zu einem nichttrivialen Kozyklus. Dies sieht man mit dem selben Argument wie im letzten Abschnitt durch be\-trach\-ten der letzten drei Zeilen: W"are der Kozyklus ein Korand, so w"are erst recht der Kozyklus $\sigma \mapsto (\sigma-1)XY$ in $\langle X^{2},Y^{2} \rangle$ ein Korand, was dieser aber nicht ist. 

Vertauschen wir in der Darstellung von $\tilde{U}$ jeweils die Zeilen/Spalten $1 \leftrightarrow 6,2\leftrightarrow 3,4 \leftrightarrow 5$, so erhalten wir wieder die Darstellung von $\tilde{U}^{*}$ - dieser Modul ist also selbstdual.

Damit ist durch $\sigma \mapsto (\sigma-1)\mu\nu$ ein nichttrivialer Kozyklus in $U:=\langle \mu^{2},\nu^{2},\mu\pi,\nu\pi,\pi^{2} \rangle$ gegeben, der von der Invarianten $\pi^{2}$ in $\tilde{U}^{*}=S^{2}\left(\langle \mu,\nu,\pi \rangle \right)$ annulliert wird. Da die zweite Potenz von $U$ wegen $\langle \mu,\nu,\pi \rangle \cong \langle \mu\pi,\nu\pi,\pi^{2} \rangle$ einen zu $\tilde{U}^{*}$ isomorphen Untermodul enth"alt, gen"ugt es also diesmal, blos zwei Annullatoren von aussen dazu zu addieren. Insgesamt folgt also mit dem Hauptsatz:

\begin{Bsp}
F"ur Charakteristik p=2 ist mit
\[
V^{*} = \langle \mu_{1}^{2},\nu_{1}^{2},\mu_{1}\pi_{1},\nu_{1}\pi_{1},\pi_{1}^{2} \rangle \bigoplus_{i=2}^{3} \langle \mu_{i},\nu_{i},\pi_{i} \rangle\\
\]
der Invariantenring \INVSL2 nicht Cohen-Macaulay, wobei $\dim V=11$. 
Der Kozyklus liegt im Grad 1 im Untermodul $\langle \mu_{1}^{2},\nu_{1}^{2},\mu_{1}\pi_{1},\nu_{1}\pi_{1},\pi_{1}^{2} \rangle$,
und ein phsop an annullierenden Invarianten im Grad $2$ (welches keine regul"are Sequenz ist), ist gegeben durch $(\pi_{1}^{2})^{2},\pi_{2}^{2},\pi_{3}^{2}$.
\qed
\end{Bsp}

Wir "uberpr"ufen dieses Beispiel wieder mit \emph{IsNotCohenMacaulay} (Datei {\tt Beispiel-6.07.txt}). Hier ist Folgendes zu beachten: Das phsop, das keine regul"are Sequenz ist, ist gegeben durch $(\pi_{1}^{2})^{2},\pi_{2}^{2},\pi_{3}^{2}$. Da jedoch auch $(\pi_{1}^{2}),\pi_{2},\pi_{3}$ invariant ist, wird \emph{IsNotCohenMacaulay} nun \emph{diese} Elemente in sein phsop mit aufnehmen (und die zweiten Potenzen sind damit f"urs phsop blockiert). Eine M"oglichkeit, dies zu umgehen, w"are den Parameter {\tt mdp} (den ich urspr"unglich aufgrund dieses Beispiels eingef"uhrt hatte)  von \emph{IsNotCohenMacaulay} auf $2$ zu setzen, damit gegebenenfalls erst das phsop ab Grad 2 gesucht wird. Allerdings ist dies nicht n"otig: Nach Lemma \ref{regPot} ist n"amlich auch das phsop in erster Potenz keine regul"are Sequenz (was von \emph{IsNotCohenMacaulay} er\-kannt wird, sofern es das phsop findet - und das tut es). Es ist hier nicht von vorneherein klar, wir gro\ss{} man nun {\tt dmax} setzen muss, denn die im Hauptsatz gegebene Relation war ja f"ur das  phsop im Grad 2 nichttrivial, jedoch zeigt sich, dass auch hier Grad $7=1+2+2+2$ gen"ugt. Der entsprechende Koeffizient der Hilbertreihe wird dann zu $120$ gesch"atzt, w"ahrend er nur $119$ ist.

\subsubsection{Verheftung an Invarianten}
Anstatt ben"otigte Untermoduln einfach als direkte Summanden dazu zu nehmen, kann man sie auch durch \emph{verheften} (siehe Abschnitt \ref{Verheft}) hinzu\-nehmen. Dabei muss man jedoch aufpassen, dass man die phsop-Eigenschaft der annul\-lierenden In\-varian\-ten nicht zerst"ort (wenn man zwei annullierende Invarianten durch Verheften zu einer verschmilzt, hat man offenbar eine weniger!) bzw. dass nichttriviale Kozyklen dies auch nach der Verheftung bleiben. Wir geben ein Beispiel f"ur dieses Verfahren:\\
Wir starten mit dem $13$-dimensionalen Modul
\[
V^{*}=\langle \mu,\nu,\pi \rangle \oplus \langle X\otimes X,...,Y\otimes Y \rangle \bigoplus_{i=1}^{3} \langle X_{i},Y_{i} \rangle
\]
mit \INVSL2 nicht Cohen-Macaulay; der nichttriviale Kozyklus liegt dabei im ersten Summanden, und die drei Annullatoren sind gegeben durch $X\otimes Y+Y\otimes X,X_{1}Y_{2}+X_{2}Y_{1},X_{2}Y_{3}+X_{3}Y_{2}$. Die Invariante $\pi$ des ersten Summanden wird jedoch nicht zur Annullation ben"otigt. Wir k"onnen daher die ersten beiden Summanden an ihren Invarianten (genauer: an dem von den In\-varian\-ten erzeugten eindimensionalen Untermoduln) zu einem einzigen Modul ver\-heften. Die Verheftung der beiden Darstellungen
\[
\sigma \mapsto
\left( \begin{array}{ccc}
d^{2} & c^{2} & 0\\
b^{2} & a^{2} & 0\\
bd & ac & 1
\end{array} \right)
\textrm{ bzw. }
\left( \begin{array}{cccc}
a^{2} & b^{2} & ab & 0\\
c^{2} & d^{2} & cd & 0\\
0 &0&1&0\\
ac & bd & bc & 1
\end{array} \right)
\]
(siehe Tabelle \ref{Darstellungen} bzw. (\ref{TensInv}) auf S. \pageref{TensInv}) an der gemeinsamen Invarianten f"uhrt dann zu einem Modul $U$ mit der Darstellung
\[
\sigma \mapsto
\left( \begin{array}{cccccc}
d^{2} & c^{2} & 0&0&0&0\\
b^{2} & a^{2} & 0&0&0&0\\
0&0&a^{2} & b^{2} & ab & 0\\
0&0&c^{2} & d^{2} & cd & 0\\
0&0&0 &0&1&0\\
bd & ac &ac & bd & bc & 1
\end{array} \right).
\]
Um zu sehen, dass der Kozyklus nichttrivial bleibt, erweitern wir $U$ mit Hilfe des Kozyklus (siehe (\ref{TensInv0}) auf S. \pageref{TensInv0}) zu $\tilde{U}$ mit der Darstellung
\[
\sigma \mapsto
\left( \begin{array}{ccccccc}
d^{2} & c^{2} & 0&0&0&0&cd\\
b^{2} & a^{2} & 0&0&0&0&ab\\
0&0&a^{2} & b^{2} & ab & 0&0\\
0&0&c^{2} & d^{2} & cd & 0&0\\
0&0&0 &0&1&0&0\\
bd & ac &ac & bd & bc & 1&bc\\
0&0&0&0&0&0&1
\end{array} \right).
\]
Dass dieser Kozyklus nichttrivial ist, sehen wir wieder sofort durch Betrachtung der ersten beiden Zeilen, denn bereits $\sigma \mapsto (\sigma-1)XY$ ist nichttrivial in $\langle X^{2},Y^{2} \rangle$. Wir erhalten also

\begin{Bsp}
F"ur Charakteristik p=2 ist mit
\[
V^{*} = U \bigoplus_{i=1}^{3} \langle X_{i},Y_{i} \rangle
\]
der Invariantenring \INVSL2 nicht Cohen-Macaulay, wobei $\dim
V=12$. Dabei ist $U$ die Verheftung von $\langle \mu,\nu,\pi \rangle$
und $\langle X\otimes X,...,Y\otimes Y \rangle$ an den Invarianten
$\pi$ bzw. $X\otimes Y+Y\otimes X$.
Der Kozyklus liegt im Grad 1 im Untermodul $U$, und dort liegt auch
eine der annullierenden Invarianten. Zwei weitere annullierende
Invarianten im Grad $2$ sind gegeben durch
$X_{1}Y_{2}+X_{2}Y_{1},X_{2}Y_{3}+X_{3}Y_{2}$.

\qed
\end{Bsp}

\noindent Gegen"uber dem Ausgangsmodul haben wir also eine Dimension gespart.

Aufgrund des phsop in den Graden $1,2,2$ muss man {\tt dmax:=1+1+2+2=6} f"ur \emph{IsNotCohenMacaulay} setzen (Datei {\tt Beispiel-6.08.txt}). Im Grad $6$ werden dann 137 Invarianten erwartet, w"ahrend es nur 136 sind.

\subsection{Beispiele f"ur \SL2 in Charakteristik 3}
Wir beginnen wieder damit, dass wir den nach Satz \ref{PKoz} gegebenen Modul mit nichttrivialem Kozyklus explizit berechnen. Wir verwenden die dortige Notation bis auf $X:=X_{1}$ und $Y:=X_{2}$. Dann ist
\[
\begin{array}{rcl}
V&:=&S^{3}(\langle X,Y \rangle)=\langle X^{3},Y^{3},X^{2}Y,XY^{2} \rangle\\
\textrm{und }W&:=&\langle X^{3},Y^{3} \rangle.
\end{array}
\]
Als erstes berechnen wir die Darstellung auf $V$. Mit $\sigma =\left( \begin{array}{cc} a&b\\ c&d \end{array} \right)$ ist

\begin{equation} \label{SL2p3}
\begin{array}{rcl}
\sigma\cdot X^{3} &=&(aX+cY)^{3}=a^{3}X^{3}+c^{3}Y^{3}\\
&\cong&(a^{3},c^{3},0,0)^{T}\\
&&\\
\sigma\cdot Y^{3} &=&(bX+dY)^{3}=b^{3}X^{3}+d^{3}Y^{3}\\
&\cong&(b^{3},d^{3},0,0)^{T}\\
&&\\
\sigma \cdot X^{2}Y&=&(aX+cY)^{2}(bX+dY)=(a^{2}X^{2}+c^{2}Y^{2}-acXY)(bX+dY)\\
&=&(a^{2}bX^{3}+c^{2}dY^{3}+a(ad-bc)X^{2}Y-c(ad-bc)XY^{2})\\
&\cong&(a^{2}b,c^{2}d,a,-c)^{T}\\
&&\\
\sigma \cdot XY^{2}&=&(aX+cY)(bX+dY)^{2}=(aX+cY)(b^{2}X^{2}+d^{2}Y^{2}-bdXY)\\
&=&(ab^{2}X^{3}+cd^{2}Y^{3}-b(ad-bc)X^{2}Y+d(ad-bc)XY^{2})\\
&\cong&(ab^{2},cd^{2},-b,d)^{T},\\
\end{array}
\end{equation}
so dass die Darstellung auf $V$ gegeben ist durch
\[
\sigma \mapsto A_{\sigma}=
\left( \begin{array}{cccc}
a^{3}&b^{3}&a^{2}b&ab^{2}\\
c^{3}&d^{3}&c^{2}d&cd^{2}\\
0&0&a&-b\\
0&0&-c&d
\end{array} \right).
\]
Die linke obere Teilmatrix beschreibt dabei die Darstellung auf $W$, die rechte untere Teilmatrix die auf $V/W$. Durch auswerten an 
\[
\left( \begin{array}{cc} a&b\\ c&d \end{array} \right)^{-1}=\left( \begin{array}{cc} d&-b\\ -c&a \end{array} \right)
\]
 und transponieren erhalten wir f"ur $(V/W)^{*}$ die Darstellung $\sigma \mapsto \left( \begin{array}{cc} d&c\\ b&a \end{array} \right)$, was die Darstellung von $\langle Y,X \rangle$ ist. Damit haben wir also  $(V/W)^{*}\cong \langle X,Y \rangle$, und f"ur den Modul mit dem Kozyklus $U$ nach Satz \ref{PKoz} gilt dann mit Satz \ref{Hom0}
\[
U=\textrm{Hom}_{K}(V,W)_{0} \cong W \otimes (V/W)^{*}=\langle X^{3},Y^{3} \rangle \otimes \langle X,Y \rangle.
\]
Da wir das Tensorprodukt als Untermodul der zweiten Potenz der direkten Summe wiederfinden, also sp"ater $U$ durch $\langle X^{3},Y^{3} \rangle \oplus \langle X,Y \rangle$ ersetzen k"onnen, erhalten wir so einen Kozyklus in Grad $2$ - leider hat die Summe hier jedoch dieselbe Dimension wie das Tensorprodukt, n"amlich $4$. Dies f"uhrt hier also nur zu einer Strukturvereinfachung und noch nicht zu einer niedrigeren Dimension. Wir berechnen nun die Darstellungen von $U$ und $\tilde{U}$. Dazu identifizieren wir
\[
U:=\textrm{Hom}_{K}(V,W)_{0}=\{ f \in \textrm{Hom}_{K}(V,W): f|_{W}=0\}
\]
mit Darstellungsmatrizen der Form
\[
\left( \begin{array}{cccc}
0 & 0 & x_{1} & x_{2}\\
0 & 0 & x_{3} & x_{4}
\end{array} \right).
\]
Die entsprechende Operation hierauf ist dann gegeben durch
\[
\sigma \cdot
\left( \begin{array}{cccc}
0 & 0 & x_{1} & x_{2}\\
0 & 0 & x_{3} & x_{4}
\end{array} \right)=
\left( \begin{array}{cc} a^{3}&b^{3}\\ c^{3}&d^{3} \end{array} \right)
\left( \begin{array}{cccc}
0 & 0 & x_{1} & x_{2}\\
0 & 0 & x_{3} & x_{4}
\end{array} \right) A_{\sigma^{-1}}.
\]
Von $A_{\sigma^{-1}}$ ist hierf"ur nur der rechte untere Block interessant, so dass wir letztendlich die Operation
\[
\begin{array}{rcl}
\sigma \cdot
\left( \begin{array}{cc}
x_{1} & x_{2}\\
x_{3} & x_{4}
\end{array} \right)&=&
\left( \begin{array}{cc} a^{3}&b^{3}\\ c^{3}&d^{3} \end{array} \right)
\left( \begin{array}{cc}
x_{1} & x_{2}\\
x_{3} & x_{4}
\end{array} \right) 
\left( \begin{array}{cc} d&b\\ c&a \end{array} \right)\\
&=&\left( \begin{array}{cc} a^{3}&b^{3}\\ c^{3}&d^{3} \end{array} \right)
\left( \begin{array}{cc}
x_{1} & x_{2}\\
x_{3} & x_{4}
\end{array} \right) 
\left( \begin{array}{cc} d&c\\ b&a \end{array} \right)^{T}
\end{array}
\]
erhalten. Wir haben also nach Lemma \ref{Tensor} eine Darstellung von $\langle X^{3},Y^{3} \rangle \otimes \langle Y,X \rangle$, wobei die Koordinaten $(x_{1},x_{2},x_{3},x_{4})^{T}$ zu $x_{1}X^{3}\otimes Y+x_{2}X^{3}\otimes X+x_{3}Y^{3}\otimes Y+x_{4}Y^{3}\otimes X$ geh"oren. Bez"uglich dieser Koordinaten hat dann jedenfalls $U$ nach Lemma \ref{TensDarst} die Darstellung
\begin{equation} \label{S4XY}
\sigma \mapsto
\left( \begin{array}{cc} a^{3}&b^{3}\\ c^{3}&d^{3} \end{array} \right)
\otimes
\left( \begin{array}{cc} d&c\\ b&a \end{array} \right)
=
\left( \begin{array}{cccc}
a^{3}d&a^{3}c&b^{3}d&b^{3}c\\
a^{3}b&a^{4}&b^{4}&ab^{3}\\
c^{3}d&c^{4}&d^{4}&cd^{3}\\
bc^{3}&ac^{3}&bd^{3}&ad^{3}
\end{array} \right).
\end{equation}
Dies ist auch die Darstellungsmatrix des Untermoduls $\langle X^{3}Y,X^{4},Y^{4},XY^{3} \rangle$ von $S^{4}(\langle X,Y \rangle)$, wie man anhand des durch $X^{3} \otimes Y \mapsto X^{3}Y, \; X^{3} \otimes X \mapsto X^{4}, \; Y^{3} \otimes Y \mapsto Y^{4}, \; Y^{3} \otimes X \mapsto XY^{3}$ gegebenen Isomorphismus erkennt.\\
Wir ben"otigen noch die Darstellung des Kozyklus, die wir gem"a\ss{} Satz \ref{PKoz} aus
\[
\begin{array}{rcl}
\sigma \cdot  
\left( \begin{array}{cccc}
1 & 0 & 0 & 0\\
0 & 1 & 0 & 0
\end{array} \right)&=&
\left( \begin{array}{cc} a^{3}&b^{3}\\ c^{3}&d^{3} \end{array} \right)
\left( \begin{array}{cccc}
1 & 0 & 0 & 0\\
0 & 1 & 0 & 0
\end{array} \right) A_{\sigma^{-1}}\\
&=&
\left( \begin{array}{cccc}
a^{3} & b^{3} & 0 & 0\\
c^{3} & d^{3} & 0 & 0
\end{array} \right)
\left( \begin{array}{cccc}
d^{3}&-b^{3}&-bd^{2}&b^{2}d\\
-c^{3}&a^{3}&ac^{2}&-a^{2}c\\
0&0&d&b\\
0&0&c&a
\end{array} \right)\\
&=&
\left( \begin{array}{cccc}
1 & 0 & ab^{3}c^{2}-a^{3}bd^{2} & a^{3}b^{2}d-a^{2}b^{3}c\\
0 & 1 & ac^{2}d^{3}-bc^{3}d^{2} & b^{2}c^{3}d-a^{2}cd^{3}
\end{array} \right)\\
&=&
\left( \begin{array}{cccc}
1 & 0 & -ab(ad+bc) & a^{2}b^{2}\\
0 & 1 & c^{2}d^{2} & -cd(ad+bc)
\end{array} \right)
\end{array}
\]
zu
\[
g_{\sigma} \cong (-ab(ad+bc),a^{2}b^{2},c^{2}d^{2}, -cd(ad+bc))^{T}
\]
erhalten. F"ur die Darstellung auf $\tilde{U}$ erhalten wir damit
\begin{equation} \label{US3}
\sigma \mapsto
\left( \begin{array}{ccccc}
a^{3}d&a^{3}c&b^{3}d&b^{3}c&-ab(ad+bc)\\
a^{3}b&a^{4}&b^{4}&ab^{3}&a^{2}b^{2}\\
c^{3}d&c^{4}&d^{4}&cd^{3}&c^{2}d^{2}\\
bc^{3}&ac^{3}&bd^{3}&ad^{3}&-cd(ad+bc)\\
0&0&0&0&1
\end{array} \right).
\end{equation}
Um diesen Kozyklus in $\langle X^{3}Y,X^{4},Y^{4},XY^{3} \rangle$ zu interpretieren, berechnen wir in $S^{4}(\langle X,Y \rangle)$
\[
\begin{array}{rcl}
\sigma \cdot X^{2}Y^{2}&=&(aX+cY)^{2}(bX+dY)^{2}\\
&=&
(a^{2}X^{2}+c^{2}Y^{2}-acXY)(b^{2}X^{2}+d^{2}Y^{2}-bdXY)\\
&=&-(ab^{2}c+a^{2}bd)X^{3}Y+a^{2}b^{2}X^{4}+c^{2}d^{2}Y^{4}+\\
&&-(acd^{2}+bc^{2}d)XY^{3}+(a^{2}d^{2}+b^{2}c^{2}-2abcd)X^{2}Y^{2}\\
&=&-ab(ad+bc)X^{3}Y+a^{2}b^{2}X^{4}+c^{2}d^{2}Y^{4}-cd(ad+bc)XY^{3}+X^{2}Y^{2}.
\end{array}
\]
Der Kozyklus entspricht also $\sigma \mapsto (\sigma-1)X^{2}Y^{2}$.\\

\begin{Zus}
Der Untermodul $U:=\langle X^{4},X^{3}Y,XY^{3},Y^{4} \rangle$ von $\tilde{U}=\langle X^{4},X^{3}Y,X^{2}Y^{2},XY^{3},Y^{4} \rangle =S^{4}\left(\langle X,Y \rangle \right)$ hat einen nicht\-trivial\-en Ko\-zyk\-lus, der durch $g_{\sigma}:=(\sigma-1)X^{2}Y^{2}$ gegeben ist. Es ist \[
U \cong \langle X^{3},Y^{3} \rangle \otimes \langle X,Y \rangle \le S^{2}\left(\langle X^{3},Y^{3} \rangle \oplus \langle X,Y \rangle \right).
\]
$U$ ist selbstdual, da es die Faktoren $\langle X^{3},Y^{3} \rangle$ und $\langle X,Y \rangle$ des Tensorprodukts sind.
Eine Darstellung von $\tilde{U}$ ist gegeben durch (\ref{US3}).
\end{Zus}

\noindent Nach Satz \ref{NCM} ist dann also mit
\[
\begin{array}{rcl}
V^{*}&:=&\langle X^{4},X^{3}Y,XY^{3},Y^{4} \rangle \oplus_{i=1}^{3} \tilde{U}^{*}\\
\textrm{bzw. } V^{*}&:=&\langle X^{3},Y^{3} \rangle \oplus \langle X,Y \rangle \oplus_{i=1}^{3} \tilde{U}^{*}
\end{array}
\]
\INVSL2 nicht Cohen-Macaulay, wobei der Kozyklus im Grad $1$ bzw. $2$
liegt und die Invarianten beide Male im Grad $1$. Die Dimension ist
ebenfalls in beiden F"allen gleich $19$. Dies alles stellt bis jetzt
blos das Beispiel aus Satz \ref{p72} in expliziter Form dar. Wir
testen wieder mit \emph{IsNotCohenMacaulay} (Datei {\tt
  Beispiel-6.09.txt}). Im Fall der Kozyklen im Grad $1$ ({\tt
  dmax:=4}) gibt es keine gro\ss e "Uberraschung - im Grad $4$ werden
$60$ Invarianten gesch"atzt, w"ahrend es nur $59$ sind. Bei dem
Beispiel mit Kozyklus im Grad $2$ (also {\tt dmax:=5} nach
(\ref{dschaetz}), S. \pageref{dschaetz}) werden im Grad $5$ dagegen $102$ (Basis-)Invarianten ge\-sch"atzt, aber es sind nur $98$! Es muss hier also weitere Relationen geben. F"ur dieses Beispiel ben"otigt {\tt IsNotCohenMacaulay} "ubrigens bereits gut eine viertel Stun\-de.\\

\noindent Wir wollen nun die Dimension weiter reduzieren, indem wir $\tilde{U}^{*}$ in einer zweiten Potenz suchen. Zun"achst berechnen wir die Darstellung von $\tilde{U}^{*}$, indem wir (\ref{US3}) invertieren (an der Stelle $\left( \begin{array}{cc} d&-b\\ -c&a \end{array} \right)$ auswerten) und trans\-po\-nier\-en zu
\begin{equation} \label{USS3}
\sigma \mapsto
\left( \begin{array}{ccccc}
ad^{3} & -bd^{3} & -ac^{3} & bc^{3} & 0\\
-cd^{3} & d^{4} & c^{4} & -c^{3}d & 0\\
-ab^{3} & b^{4} & a^{4} & -a^{3}b & 0\\
b^{3}c & -b^{3}d & -a^{3}c & a^{3}d &0\\
bd(ad+bc) & b^{2}d^{2} & a^{2}c^{2} & ac(ad+bc) & 1
\end{array} \right).
\end{equation}
Als n"achste betrachten wir den Modul $M:=\langle X^{2},Y^{2},XY \rangle$. Dieser hat nach Tabelle \ref{Darstellungen} die Darstellung
\[
\sigma \mapsto
\left( \begin{array}{ccc}
a^{2} & b^{2} & ab\\
c^{2} & d^{2} & cd\\
-ac & -bd & ad+bc
\end{array} \right).
\]
Die Darstellung des Duals ist damit gegeben durch 
\[
\sigma \mapsto
\left( \begin{array}{ccc}
d^{2} & c^{2} & cd\\
b^{2} & a^{2} & ab\\
-bd & -ac & ad+bc
\end{array} \right).
\]
Durch Vertauschen der ersten beiden Zeilen und Spalten erh"alt man wieder die Darstellung von $M$, so dass dieser Modul selbstdual ist. Wir berechnen nun
\[
S^{2}\left( M \right) = \langle (XY)Y^{2},-(Y^{2})^{2},-(X^{2})^{2},(XY)X^{2},X^{2}Y^{2}-(XY)^{2},(XY)^{2} \rangle
\]
bez"uglich der angedeuteten Basis (man unterscheide hier $(XY)^{2}$ und $X^{2}Y^{2}$!). Wir erhalten \label{S2M}
\[
\begin{array}{rcl}
\sigma \cdot (XY)Y^{2}&=&(abX^{2}+cdY^{2}+(ad+bc)XY)(b^{2}X^{2}+d^{2}Y^{2}-bdXY)\\
&=&(XY)Y^{2}(ad^{3}+bcd^{2}-bcd^{2})-(Y^{2})^{2}(-cd^{3})-(X^{2})^{2}(-ab^{3})+\\
&&+(XY)X^{2}(ab^{2}d+b^{3}c-ab^{2}d)
+(X^{2}Y^{2})(abd^{2}+b^{2}cd)+\\
&&-(XY)^{2}(abd^{2}+b^{2}cd)\\
&\cong&(ad^{3},-cd^{3},-ab^{3},b^{3}c,bd(ad+bc),0)^{T}\\
\end{array}
\]
\[
\begin{array}{rcl}
\sigma \cdot -(Y^{2})^{2}&=&-(b^{2}X^{2}+d^{2}Y^{2}-bdXY)^{2}\\
&=&(XY)Y^{2}(-bd^{3})-(Y^{2})^{2}(d^{4})-(X^{2})^{2}(b^{4})+(XY)X^{2}(-b^{3}d)+\\
&&+X^{2}Y^{2}(b^{2}d^{2})-(XY)^{2}(b^{2}d^{2})\\
&\cong&(-bd^{3},d^{4},b^{4},-b^{3}d,b^{2}d^{2},0)^{T}\\
\end{array}
\]
\[
\begin{array}{rcl}
\sigma \cdot -(X^{2})^{2}&=&-(a^{2}X^{2}+c^{2}Y^{2}-acXY)^{2}\\
&=&(XY)Y^{2}(-ac^{3})-(Y^{2})^{2}(c^{4})-(X^{2})^{2}(a^{4})+(XY)X^{2}(-a^{3}c)+\\
&&+X^{2}Y^{2}(a^{2}c^{2})-(XY)^{2}(a^{2}c^{2})\\
&\cong&(-ac^{3},c^{4},a^{4},-a^{3}c,a^{2}c^{2},0)^{T}\\
\end{array}
\]
\[
\begin{array}{rcl}
\sigma \cdot (XY)X^{2}&=&(abX^{2}+cdY^{2}+(ad+bc)XY)(a^{2}X^{2}+c^{2}Y^{2}-acXY)\\
&=&(XY)Y^{2}(ac^{2}d+bc^{3}-ac^{2}d)-(Y^{2})^{2}(-c^{3}d)-(X^{2})^{2}(-a^{3}b)+\\
&&+(XY)X^{2}(a^{3}d+a^{2}bc-a^{2}bc)
+(X^{2}Y^{2})(abc^{2}+a^{2}cd)+\\
&&-(XY)^{2}(a^{2}cd+abc^{2})\\
&\cong&(bc^{3},-c^{3}d,-a^{3}b,a^{3}d,ac(ad+bc),0)^{T}\\
\end{array}
\]
\[
\begin{array}{rcl}
\sigma \cdot (X^{2}Y^{2}-(XY)^{2})&=&(a^{2}X^{2}+c^{2}Y^{2}-acXY)(b^{2}X^{2}+d^{2}Y^{2}-bdXY)+\\
&&-(abX^{2}+cdY^{2}+(ad+bc)XY)^{2}\\
&=&(X^{2})^{2}(a^{2}b^{2}-(ab)^{2})+(Y^{2})^{2}(c^{2}d^{2}-(cd)^{2})+\\
&&+(XY)^{2}(4abcd-(ad+bc)^{2})+\\
&&+X^{2}Y^{2}(a^{2}d^{2}+b^{2}c^{2}-2abcd)\\
&&+(XY)X^{2}(-a^{2}bd-ab^{2}c+a^{2}bd+ab^{2}c)+\\
&&+(XY)Y^{2}(-acd^{2}-bc^{2}d+acd^{2}+bc^{2}d)+\\
&\cong&(0,0,0,0,1,0)^{T}
\end{array}
\]
\[
\begin{array}{rcl}
\sigma \cdot (XY)^{2}&=&(abX^{2}+cdY^{2}+(ad+bc)XY)^{2}\\
&=&(XY)Y^{2}(-cd(ad+bc))-(Y^{2})^{2}(-c^{2}d^{2})-(X^{2})^{2}(-a^{2}b^{2})+\\
&&+(XY)X^{2}(-ab(ad+bc))+\\
&&+(X^{2}Y^{2}-(XY)^{2})(-abcd)+(XY)^{2}((ad+bc)^{2}-4abcd)\\
&\cong&(-cd(ad+bc),-c^{2}d^{2},-a^{2}b^{2},-ab(ad+bc),-abcd,1)^{T}.
\end{array}
\]
Eine Darstellung auf $S^{2}(M)$ ist damit gegeben durch
\begin{equation} \label{DarstS2M}
\sigma \mapsto
\left( \begin{array}{cccccc}
ad^{3} & -bd^{3} & -ac^{3} & bc^{3} & 0&-cd(ad+bc)\\
-cd^{3} & d^{4} & c^{4} & -c^{3}d & 0&-c^{2}d^{2}\\
-ab^{3} & b^{4} & a^{4} & -a^{3}b & 0&-a^{2}b^{2}\\
b^{3}c & -b^{3}d & -a^{3}c & a^{3}d &0&-ab(ad+bc)\\
bd(ad+bc) & b^{2}d^{2} & a^{2}c^{2} & ac(ad+bc) & 1 &-abcd\\
0&0&0&0&0&1
\end{array} \right).
\end{equation}
Durch Vergleich mit (\ref{USS3}) erkennen wir $\tilde{U}^{*}$ als Untermodul von $S^{2}(\langle X^{2},Y^{2},XY \rangle)$.\\ 

{\small
\noindent Au\ss erdem ist mit der letzten Spalte ein neuer Kozyklus entstanden. Er ist nicht\-trivial, wie man durch betrachten der ersten vier Zeilen sieht: Vertauscht  man n"amlich in (\ref{US3}) (mit nichttrivialen Kozyklus) die Zeilen/Spalten $1\leftrightarrow 4, 2 \leftrightarrow 3$ und multipliziert dann die Zeilen/Spalten $2,3$ mit $-1$, so erh"alt man dieselben vier Zeilen wie hier - und da dort der Kozyklus nichttrivial ist, muss er es auch hier sein. Aufgrund der hohen Dimension dieses Moduls wird dies jedoch hier nicht weiter verfolgt.\\

\noindent "Ubrigens ist $S^{2}(M)$ selbstdual. Denn vertauscht man in der invertierten und trans\-ponierten Darstellungsmatrix die Zeilen/Spalten $1\leftrightarrow 4, 2 \leftrightarrow 3, 5 \leftrightarrow 6 $ und Multi\-pliziert die Zeilen/Spalten $2,3$ mit $-1$, so erh"alt man wieder die Ausgangsmatrix.\\

\noindent Man kann die ganze Rechnung auch f"ur allgemeines $p$ durchf"uhren: Dabei ist stets $X^{2}Y^{2}-(XY)^{2}$ eine Invariante. F"ur $p \ne 3$ ist der zugeh"orige Kozyklus jedoch trivial. Rechts unten in der obigen Matrix steht dann au\ss erdem keine $1$ mehr, so dass auch kein neuer Kozyklus entsteht.\\
}

\noindent Nach all diesen negativ-Berichten nun aber zur"uck zu unserem positiven Ergebnis - wir haben den Modul $\tilde{U}^{*}$ mit der annullierenden Invariante in der zweiten Potenz von $\langle X^{2},Y^{2},XY \rangle$ wiedergefunden! Damit k"onnen wir das vielleicht in\-teressan\-teste Resultat dieser Arbeit formulieren:

\begin{samepage}
\begin{Bsp}[phsop und Kozyklus im Grad 2]
F"ur Charakteristik $p=3$ ist mit
\[
V^{*}:=\langle X,Y \rangle \oplus \langle X^{3},Y^{3} \rangle \oplus_{i=1}^{3} \langle X^{2},Y^{2},XY \rangle \quad (\cong V)
\]
der Invariantenring \INVSL2 nicht Cohen-Macaulay, wobei $\dim V=13$. Da alle Summanden selbstdual sind, ist auch $V=V^{*}$ selbstdual. Der nicht\-triviale Kozyklus liegt im Grad $2$ im zum Tensorprodukt isomorphen Unter\-modul von $S^{2}\left( \langle X,Y \rangle \oplus \langle X^{3},Y^{3} \rangle\right)$. Das annullierende phsop liegt eben\-falls im Grad $2$ und ist durch die drei Kopien der $X^{2}Y^{2}-(XY)^{2}$ gegeben. \qed
\end{Bsp}
\end{samepage}

\noindent Wenn wir dieses Beispiel mit \emph{IsNotCohenMacaulay} untersuchen wollen (Datei {\tt Beispiel-6.10.txt}), m"ussen wir bereits {\tt dmax:=8=2+2+2+2} setzen, denn nun liegen sowohl phsop als auch Kozyklus in Grad $2$. Ohne den be\-schleu\-nig\-ten Algorithmus zur Berechnung der Invarianten ist es v"ollig undenkbar, hier in akzeptabler Zeit ein Ergebnis zu erhalten - denn selbst der beschleunigte Algorithmus (der bei den bisherigen Beispielen bis zu $500$ mal schneller war) ben"otigt hier gut eineinhalb Stunden! Als weitere Be\-schleunigungs\-ma\ss nah\-me wurde "ubrigens der Parameter {\tt maxdp:=2} gesetzt, d.h. ein phsop wird blos bis zum Grad $2$ gesucht (die Berechnung der Dimensionen geschieht n"amlich "uber Gr"obner-Basen und ist bei diesem gro\ss en Beispiel ebenfalls recht langsam). Im Grad $8$ werden dann $366$ Invarianten gesch"atzt, w"ahrend es nur $362$ sind - also gibt es wieder mehrere Relationen.

\subsection{Beispiele f"ur \GL2}
In diesem Abschnitt wollen wir die Beispiele f"ur \SL2 so modifizieren, dass daraus Beispiele f"ur \GL2 werden. Zun"achst wollen wir nochmal daran erinnern, dass die Elemente der \emph{algebraischen} Gruppe \GL2 (vorerst $\textrm{char } K=p$ beliebig) gegeben sind durch $\sigma=\left( \begin{array}{cc} a &b\\ c &d \end{array} \right)_{e}$ mit $(ad-bc)e=1$. Das Inverse ist dann gegeben durch $\left( \begin{array}{cc} a &b\\ c &d \end{array} \right)_{e}^{-1}=\left( \begin{array}{cc} ed &-eb\\ -ec & ea \end{array} \right)_{ad-bc}$. Wir definieren den eindimensionalen Modul $\langle E \rangle$ durch die Operation, die durch Multiplikation mit der Inversen der Determinanten gegeben ist, also
\[
\left( \begin{array}{cc} a &b\\ c &d \end{array} \right)_{e} \cdot E := e \cdot E =\frac{1}{ad-bc}E.
\]
Beim Versuch der "Ubertragung von Beispielen f"ur \SL2 auf \GL2 kommt es oft vor, dass bei der Berechnung der Darstellungsmatrix ein Eintrag der Form $ad-bc$ (evtl. in h"ohere Potenz) auftaucht. F"ur \SL2 kann man dies durch $1$ ersetzen, und gegebenenfalls erh"alt man so einen Kozyklus oder eine Invariante. Als Beispiel betrachte man f"ur $p=2$ die Darstellung von $\langle X^{2},Y^{2},XY \rangle$ f"ur \SL2 und \GL2 in Tabelle \ref{Darstellungen}. Der Trick, hier eine $1$ zu erzeugen, besteht aus Tensorierung mit $\langle E \rangle$. F"ur $p=2$ hat dann zum Beispiel
\[
\langle X^{2},Y^{2},XY \rangle \otimes \langle E \rangle = \langle X^{2} \otimes E,Y^{2}\otimes E,XY\otimes E \rangle
\]
f"ur \GL2 die Darstellung
\begin{equation} \label{XYGL}
\sigma=\left( \begin{array}{cc} a &b\\ c &d \end{array} \right)_{e} \mapsto
\left( \begin{array}{ccc}
ea^{2} & eb^{2} & eab\\
ec^{2}&ed^{2} &ecd\\
0&0&1
\end{array} \right).
\end{equation}
Damit erh"alt man dann einen Kozyklus in $\langle X^{2},Y^{2}\rangle \otimes \langle E \rangle$, der nichttrivial ist, da ja bereits die Beschr"ankung auf \SL2 (also $e=1$) nichttrivial ist.

Kozyklen im Grad $1$ und (zugeh"orige annullierende Invarianten) werden f"ur die uns interessierenden Beispiele jedoch sowieso von Satz \ref{PKoz} geliefert, denn dort war ja auch \GL2 als Gruppe zugelassen. Wir m"ussen uns hier also nur noch damit besch"aftigen, wie die im vorigen Abschnitt gefundenen Beispiele mit phsop bzw. Kozyklen in h"oheren Potenzen zu Tensorieren sind, so dass auch hier wieder die h"oheren Potenzen Untermoduln der in Satz \ref{PKoz} gefun\-denen Moduln mit Kozyklus bzw. annullierender Invarianten sind. Wir wer\-den dies im Folgenden lediglich f"ur die beiden $10$ bzw. $13$-dimensionalen Haupt\-beispiele in Charakteristik $2$ bzw. $3$ durchf"uhren.

\subsubsection{Beispiele f"ur \GL2 in Charakteristik 2}
Wir exerzieren Satz \ref{PKoz} f"ur $p=2$ und \GL2 durch, um den dort versprochenen nichttrivialen Kozyklus zu erhalten: Wir setzen
\[
V:=\langle X^{2},Y^{2},XY \rangle
\]
mit Darstellung
\[
\sigma=\left( \begin{array}{cc} a &b\\ c &d \end{array} \right)_{e} \mapsto
A_{\sigma}=\left( \begin{array}{ccc}
a^{2} & b^{2} & ab\\
c^{2} & d^{2} & cd\\
0 & 0 &ad+bc
\end{array} \right)
\]
und Untermodul $W:=\langle X^{2},Y^{2} \rangle$, und $U:=\textrm{Hom}_{K}(V,W)_{0}$. Wir berechnen die Darstellung von $\tilde{U}$; Dabei beachten wir, dass wir $A_{\sigma}^{-1}$ durch auswerten an $\left( \begin{array}{cc} ed &eb\\ ec & ea \end{array} \right)_{ad+bc}$ erhalten. Damit erhalten wir aus

\begin{eqnarray*}
\lefteqn{ \left( \begin{array}{cc}
a^{2} & b^{2}\\
c^{2}&d^{2}
\end{array} \right)
\left( \begin{array}{ccc}
0 &0&1\\
0&0&0
\end{array} \right)A_{\sigma^{-1}}}\\
&=&
\left( \begin{array}{ccc}
0 &0&a^{2}\\
0&0&c^{2}
\end{array} \right)
\left( \begin{array}{ccc}
e^{2}d^{2} & e^{2}b^{2} & e^{2}bd\\
e^{2}c^{2} & e^{2}a^{2} & e^{2}ac\\
0 & 0 &e^{2}ad+e^{2}bc=e
\end{array} \right)\\
&=&
\left( \begin{array}{ccc}
0 &0&ea^{2}\\
0&0&ec^{2}
\end{array} \right) \cong (ea^{2},ec^{2},0)^{T}
\end{eqnarray*}
sowie

\begin{eqnarray*}
\lefteqn{\left( \begin{array}{cc}
a^{2} & b^{2}\\
c^{2}&d^{2}
\end{array} \right)
\left( \begin{array}{ccc}
0 &0&0\\
0&0&1
\end{array} \right)
\left( \begin{array}{ccc}
e^{2}d^{2} & e^{2}b^{2} & e^{2}bd\\
e^{2}c^{2} & e^{2}a^{2} & e^{2}ac\\
0 & 0 &e
\end{array} \right)}\\
&=&
\left( \begin{array}{ccc}
0 &0&b^{2}\\
0&0&d^{2}
\end{array} \right)
\left( \begin{array}{ccc}
e^{2}d^{2} & e^{2}b^{2} & e^{2}bd\\
e^{2}c^{2} & e^{2}a^{2} & e^{2}ac\\
0 & 0 &e
\end{array} \right)=
\left( \begin{array}{ccc}
0 &0&eb^{2}\\
0&0&ed^{2}
\end{array} \right) \cong (eb^{2},ed^{2},0)^{T}
\end{eqnarray*}

und
\begin{eqnarray*} 
\lefteqn{\left( \begin{array}{cc}
a^{2} & b^{2}\\
c^{2}&d^{2}
\end{array} \right)
\left( \begin{array}{ccc}
1 &0&0\\
0&1&0
\end{array} \right)
\left( \begin{array}{ccc}
e^{2}d^{2} & e^{2}b^{2} & e^{2}bd\\
e^{2}c^{2} & e^{2}a^{2} & e^{2}ac\\
0 & 0 &e
\end{array} \right)}\\
&=&\left( \begin{array}{ccc}
a^{2} & b^{2} &0\\
c^{2}&d^{2}&0
\end{array} \right)
\left( \begin{array}{ccc}
e^{2}d^{2} & e^{2}b^{2} & e^{2}bd\\
e^{2}c^{2} & e^{2}a^{2} & e^{2}ac\\
0 & 0 &e
\end{array} \right)\\
&=&
\left( \begin{array}{ccc}
1 &0&e^{2}a^{2}bd+e^{2}ab^{2}c=e^{2}ab(ad+bc)=eab\\
0&1&e^{2}bc^{2}d+e^{2}acd^{2}=e^{2}cd(ad+bc)=ecd
\end{array} \right) \cong (eab,ecd,1)^{T}
\end{eqnarray*}
f"ur $\tilde{U}$ die Darstellung (\ref{XYGL}) von $\langle X^{2} \otimes E,Y^{2}\otimes E,XY\otimes E \rangle$ - also genau dass, was wir erwartet haben. Insbesondere ist 
\[
U =\langle X^{2} \otimes E,Y^{2}\otimes E \rangle = \langle X^{2} ,Y^{2} \rangle \otimes \langle E \rangle.
\]
Wir berechnen die Darstellung des Duals $\tilde{U}^{*}$, indem wir (\ref{XYGL}) transponieren und an der Stelle $\left( \begin{array}{cc} ed &eb\\ ec & ea \end{array} \right)_{ad+bc}$  auswerten zu
\[
\sigma \mapsto
\left( \begin{array}{ccc}
e^{2}d^{2}(ad+bc) & e^{2}c^{2}(ad+bc) & 0\\
e^{2}b^{2}(ad+bc) & e^{2}a^{2}(ad+bc) &0\\
e^{2}bd(ad+bc) & e^{2}ac(ad+bc) & 1
\end{array} \right)
=
\left( \begin{array}{ccc}
ed^{2} & ec^{2} & 0\\
eb^{2} & ea^{2}&0\\
ebd & eac & 1=e(ad+bc)
\end{array} \right).
\]
F"ur \SL2 hatten wir $\tilde{U}^{*}$ in der zweiten Potenz von $\langle X,Y\rangle \oplus \langle X,Y\rangle$ wiedergefunden. Da in obiger Darstellung $e$ jeweils nur in erster Potenz vorkommt, liegt es nahe, einen der beiden Summanden mit $\langle E \rangle$ zu tensorieren. In Gleichung (\ref{XXYYbasis}) auf S. \pageref{XXYYbasis} haben wir die Darstellung von $\langle X,Y\rangle \otimes \langle X,Y\rangle$ bzgl. der Basis
\[
\mathcal{B}=\left( X \otimes X, Y \otimes Y, X\otimes Y+Y\otimes X,Y\otimes X \right)
\]
berechnet. Wir tensorieren diese Darstellung nochmals mit $\langle E \rangle$ und erhalten nach Lemma \ref{TensDarst} die Darstellung
\[
\sigma \mapsto
\left( \begin{array}{cccc}
ea^{2} & eb^{2} & 0 & eab\\
ec^{2} & ed^{2} & 0 & ecd\\
eac & ebd & e(ad+bc)=1 & ebc\\
0 & 0 & 0 & e(ad+bc)=1\\
\end{array} \right).
\]
Nach Lemma \ref{UVW} ist dies auch die Darstellung von $\langle X,Y\rangle \otimes (\langle X,Y\rangle \otimes \langle E \rangle)$ bzgl. der Basis
\[
\left( X \otimes (X\otimes E), Y \otimes (Y\otimes E), X\otimes (Y\otimes E)+Y\otimes (X \otimes E),Y\otimes (X\otimes E) \right).
\]
Wir finden also $\tilde{U}^{*}$ (nach Vertauschen der ersten beiden Zeilen/Spalten)  in der zweiten Potenz von $\langle X,Y\rangle \oplus (\langle X,Y\rangle \otimes \langle E \rangle)$ wieder, wobei die annullierende Invariante gegeben ist durch $X(Y\otimes E)+Y(X \otimes E)$. Wir bemerken noch, dass $\langle X\otimes E,Y \otimes E \rangle$ isomorph ist zum Dual $\langle X^{*},Y^{*} \rangle$ des nat"urlichen Moduls, mit Darstellung $\sigma \mapsto \left( \begin{array}{cc} ed &ec\\ eb & ea \end{array} \right)$ - man hat also die Entsprechungen $X\otimes E \cong Y^{*}, Y\otimes E \cong X^{*}$. Damit haben wir also das Hauptresultat f"ur $p=2$ und \GL2:\\

\begin{samepage}
\begin{Bsp}
F"ur $p=2$ und die Gruppe \GL2 ist mit
\begin{eqnarray*}
V^{*} &=& \left( \langle X^{2},Y^{2} \rangle \otimes \langle E\rangle \right) \bigoplus_{i=1}^{2} \left( \langle X_{i},Y_{i} \rangle \oplus \langle X_{i}^{*},Y_{i}^{*} \rangle \right)\\
&\cong& \left( \langle X^{2},Y^{2} \rangle \otimes \langle E\rangle \right)  \bigoplus_{i=1}^{2} \left( \langle X_{i},Y_{i} \rangle \oplus \left( \langle X_{i},Y_{i} \rangle \otimes E \right) \right)
\end{eqnarray*}
der Invariantenring \INVGL2 nicht Cohen-Macaulay, wobei $\dim V=10$. Der Kozyklus liegt im Grad $1$ im ersten Summanden, ein annullierendes phsop im Grad $2$ ist etwa gegeben durch $X_{1}X_{1}^{*}+Y_{1}Y_{1}^{*},X_{2}X_{2}^{*}+Y_{2}Y_{2}^{*},X_{1}X_{2}^{*}+Y_{1}Y_{2}^{*}$. Die Summanden $\langle X_{i},Y_{i} \rangle \oplus \langle X_{i}^{*},Y_{i}^{*} \rangle$ sind selbstdual. \qed
\end{Bsp}
\end{samepage}

\noindent Anstatt mit $\langle E \rangle$ zu tensorieren, kann man auch einfach $\langle E \rangle$ hinreichend oft addieren - dann liegen der Kozyklus bzw. das phsop noch einen Grad h"oher. Allerdings bringt die Addition nat"urlich eine unerw"unschte Di\-men\-sions\-erh"ohung mit sich. Insbesondere muss man $\langle E \rangle$ so oft addieren, dass sich noch ein phsop ergibt. Bei dem Kozyklus ist dies dagegen nicht kritisch, es kann hier durchaus ein Summand $\langle E \rangle$ verwendet werden, der auch schon f"ur eine Invariante verwendet wurde. Man kann also z.B.
\[
V^{*}:=\langle X^2,Y^2 \rangle \bigoplus_{i=1}^{4} \langle X_{i},Y_{i}\rangle \bigoplus_{i=1}^{3}\langle E_{i} \rangle
\] 
setzen. Jedes Element im phsop bekommt hier ein $E_{i}$ (deshalb braucht man drei Summanden), und ein $E_{i}$ l"asst sich f"ur den Kozyklus im Grad 2 recyclen. Also\\

\begin{Bsp}
Mit 
\[
V^{*}:=\langle X^2,Y^2 \rangle \bigoplus_{i=1}^{4} \langle X_{i},Y_{i}\rangle \bigoplus_{i=1}^{3}\langle E_{i} \rangle
\]
ist \INVGL2 nicht Cohen-Macaulay, $\dim V=13$. Ein Kozyklus im Grad $2$ liegt etwa in 
$S^{2}\left(\langle X^2,Y^2 \rangle \oplus \langle E_{1} \rangle \right)$, und ein phsop an annullierenden In\-va\-ri\-an\-ten im Grad $3$ (welches keine regul"are Sequenz ist), ist gegeben durch $X_{i}Y_{j}E_{k}+X_{j}Y_{i}E_{k}$ mit $(i,j,k) \in \left\{(1,2,1),(2,3,2),(3,4,3) \right\}$. \qed
\end{Bsp}

\noindent Man kann hier auch noch weiter kombinieren, etwa
\begin{Bsp}
Mit
\[
V^{*}:=\langle X^2,Y^2 \rangle \bigoplus_{i=1}^{3} \langle X_{i},Y_{i}\rangle \bigoplus \left( \langle X_{4},Y_{4}\rangle \otimes \langle E \rangle \right) \bigoplus \langle E \rangle
\]
ist \INVGL2 nicht Cohen-Macaulay, $\dim V=11$. Hier liegt der Kozyklus im Grad $2$ (in $S^{2}\left( \langle X^2,Y^2 \rangle \oplus \langle E \rangle \right)$, und ein phsop in den Graden $2,2,3$ ist gegeben durch $X_{1}(Y_{4} \otimes E)+Y_{1}(X_{4} \otimes E),X_{2}(Y_{4} \otimes E)+Y_{2}(X_{4} \otimes E) ,X_{2}Y_{3}E+X_{3}Y_{2}E$. \qed
\end{Bsp}

\noindent Wir kommen zum Test dieser Beispiele mit \emph{IsNotCohenMacaulay} (zusammen\-gefasst in Datei {\tt Beispiel-6.11-13.txt}). Im ersten Beispiel gibt es keine gro\ss e "Uberraschung - man muss {\tt dmax:=7} setzen, und es wird sich dann um eine Invariante ($17$ statt $16$ in Grad $7$) versch"atzt. Im zweiten Beispiel m"usste man {\tt dmax:=2+3+3+3=11} setzen - zuviel f"ur meine Geduld. Bei dem dritten Beispiel hat man {\tt dmax:=2+2+2+3=9} zu setzen - auch hier wird sich im Grad $9$ um eine Invariante versch"atzt ($70$ statt $69$).

\subsubsection{Beispiele f"ur \GL2 in Charakteristik 3}
Auch f"ur $p=3$ wollen wir die Beispiele f"ur \SL2 "ubertragen. Die Ver\-all\-gemeinerung des Moduls mit dem Kozyklus $U$ liefert wieder Satz \ref{PKoz}, der auch f"ur \GL2 gilt. Der Vollst"andigkeit halber sind hier nochmals alle Rechnungen so weit n"otig ausgef"uhrt, auch wenn sie mit denen f"ur \SL2 fast identisch sind. In den meisten F"allen gen"ugt es, die Dar\-stellungs\-ma\-tri\-zen f"ur die Gruppe \SL2 mit geeigneten Potenzen von $e$ zu multiplizieren.
Wir gehen jedoch systematisch vor und exerzieren zun"achst wieder Satz \ref{PKoz} durch:
\[
\begin{array}{rcl}
V&:=&S^{3}(\langle X,Y \rangle)=\langle X^{3},Y^{3},X^{2}Y,XY^{2} \rangle\\
\textrm{und }W&:=&\langle X^{3},Y^{3} \rangle.
\end{array}
\]
Als erstes ben"otigen wir die Darstellung auf $V$. In (\ref{SL2p3}) auf S. \pageref{SL2p3} haben wir bei der Berechnung der Darstellung f"ur die Gruppe \SL2 jeweils erst in der dritten Zeile jeder Gleichung von der Eigenschaft $ad-bc=1$ Gebrauch gemacht. Die Betrachtung der zweiten Zeilen liefert f"ur die Darstellung auf $V$
\[
\sigma \mapsto A_{\sigma}=
\left( \begin{array}{cccc}
a^{3}&b^{3}&a^{2}b&ab^{2}\\
c^{3}&d^{3}&c^{2}d&cd^{2}\\
0&0&a(ad-bc)&-b(ad-bc)\\
0&0&-c(ad-bc)&d(ad-bc)
\end{array} \right).
\]
Die linke obere Teilmatrix beschreibt dabei die Darstellung auf $W$, die rechte untere Teilmatrix die auf $V/W$. Durch Auswerten an 
\[
\left( \begin{array}{cc} a&b\\ c&d \end{array} \right)_{e}^{-1}=\left( \begin{array}{cc} ed&-eb\\ -ec&ea \end{array} \right)_{ad-bc}
\]
 und transponieren erhalten wir f"ur $(V/W)^{*}$ die Darstellung 
\begin{equation} \label{VWstar}
\sigma \mapsto 
\left( \begin{array}{cc} 
ed(e^{2}ad-e^{2}bc)&ec(e^{2}ad-e^{2}bc)\\ 
eb(e^{2}ad-e^{2}bc)&ea(e^{2}ad-e^{2}bc) \end{array} \right)=
\left( \begin{array}{cc} 
e^{2}d&e^{2}c\\ 
e^{2}b&e^{2}a \end{array} \right),
\end{equation}
was die Darstellung von $\langle Y,X \rangle \otimes \langle E^{2} \rangle$ ist (wobei $\langle E^{2} \rangle:=S^{2}\left( \langle E \rangle\right)$). Damit haben wir also  $(V/W)^{*}\cong \langle X,Y \rangle  \otimes \langle E^{2} \rangle$, und f"ur den Modul mit dem Kozyklus $U$ nach Satz \ref{PKoz} gilt dann mit Satz \ref{Hom0}
\[
U=\textrm{Hom}_{K}(V,W)_{0} \cong W \otimes (V/W)^{*}=\langle X^{3},Y^{3} \rangle \otimes\left( \langle X,Y \rangle  \otimes \langle E^{2} \rangle \right).
\]
Dieses Tensorprodukt werden wir wieder durch die direkte Summe ersetzen, so dass wir wieder einen Kozyklus im Grad $2$ erhalten. Den Faktor $\langle E^{2} \rangle$ k"onnen wir dabei auch auf die andere Seite bringen, symmetrisch aufspalten (also beide Seiten mit $\langle E \rangle$ tensorieren) oder auch als zus"atzlichen direkten Summand dazu nehmen. Dies liefert dann einen Kozyklus im Grad $3$ oder $4$ (falls man nur $\langle E \rangle$ addiert statt $\langle E^{2} \rangle$), aber nat"urlich auch eine h"ohere Dimension.

Wir berechnen nun die Darstellungen von $U$ und $\tilde{U}$. Dazu identifizieren wir
\[
U:=\textrm{Hom}_{K}(V,W)_{0}=\{ f \in \textrm{Hom}_{K}(V,W): f|_{W}=0\}
\]
wieder mit Darstellungsmatrizen der Form
\[
\left( \begin{array}{cccc}
0 & 0 & x_{1} & x_{2}\\
0 & 0 & x_{3} & x_{4}
\end{array} \right).
\]
Die entsprechende Operation hierauf ist dann gegeben durch
\[
\sigma \cdot
\left( \begin{array}{cccc}
0 & 0 & x_{1} & x_{2}\\
0 & 0 & x_{3} & x_{4}
\end{array} \right)=
\left( \begin{array}{cc} a^{3}&b^{3}\\ c^{3}&d^{3} \end{array} \right)
\left( \begin{array}{cccc}
0 & 0 & x_{1} & x_{2}\\
0 & 0 & x_{3} & x_{4}
\end{array} \right) A_{\sigma^{-1}}.
\]
Von $A_{\sigma^{-1}}$ ist hierf"ur nur der rechte untere Block interessant (und das transponierte dieses Blocks haben wir in (\ref{VWstar}) berechnet), so dass wir letztendlich die Operation
\[
\begin{array}{rcl}
\sigma \cdot
\left( \begin{array}{cc}
x_{1} & x_{2}\\
x_{3} & x_{4}
\end{array} \right)&=&
\left( \begin{array}{cc} a^{3}&b^{3}\\ c^{3}&d^{3} \end{array} \right)
\left( \begin{array}{cc}
x_{1} & x_{2}\\
x_{3} & x_{4}
\end{array} \right) 
\left( \begin{array}{cc} 
e^{2}d&e^{2}c\\ 
e^{2}b&e^{2}a \end{array} \right)^{T}
\end{array}
\]
erhalten. Wir haben also nach Lemma \ref{Tensor} eine Darstellung von $\langle X^{3},Y^{3} \rangle \otimes \left( \langle Y,X \rangle \otimes  \langle E^{2} \rangle \right)$. Nach Lemma \ref{TensDarst} erhalten wir also die folgende Darstellung bez"uglich der Koordinaten $x_{1},x_{2},x_{3},x_{4}$, die wir einfach durch Multiplikation der Darstellung (\ref{S4XY}) (S. \pageref{S4XY})  der Gruppe \SL2 mit $e^{2}$ erhalten:
\[
\sigma \mapsto
\left( \begin{array}{cc} a^{3}&b^{3}\\ c^{3}&d^{3} \end{array} \right)
\otimes
\left( \begin{array}{cc} e^{2}d&e^{2}c\\ e^{2}b&e^{2}a \end{array} \right)
=
\left( \begin{array}{cccc}
e^{2}a^{3}d&e^{2}a^{3}c&e^{2}b^{3}d&e^{2}b^{3}c\\
e^{2}a^{3}b&e^{2}a^{4}&e^{2}b^{4}&e^{2}ab^{3}\\
e^{2}c^{3}d&e^{2}c^{4}&e^{2}d^{4}&e^{2}cd^{3}\\
e^{2}bc^{3}&e^{2}ac^{3}&e^{2}bd^{3}&e^{2}ad^{3}
\end{array} \right).
\]
Nun berechnen wir den Kozyklus:
\begin{eqnarray*}
\lefteqn{\sigma \cdot  
\left( \begin{array}{cccc}
1 & 0 & 0 & 0\\
0 & 1 & 0 & 0
\end{array} \right)=
\left( \begin{array}{cc} a^{3}&b^{3}\\ c^{3}&d^{3} \end{array} \right)
\left( \begin{array}{cccc}
1 & 0 & 0 & 0\\
0 & 1 & 0 & 0
\end{array} \right) A_{\sigma^{-1}}}\\
&=&
\left( \begin{array}{cccc}
a^{3} & b^{3} & 0 & 0\\
c^{3} & d^{3} & 0 & 0
\end{array} \right)
\left( \begin{array}{cccc}
e^{3}d^{3}&-e^{3}b^{3}&-e^{3}bd^{2}&e^{3}b^{2}d\\
-e^{3}c^{3}&e^{3}a^{3}&e^{3}ac^{2}&-e^{3}a^{2}c\\
0&0&e^{2}d&e^{2}b\\
0&0&e^{2}c&e^{2}a
\end{array} \right)\\
&=&
\left( \begin{array}{cccc}
1 & 0 & e^{3}\left(ab^{3}c^{2}-a^{3}bd^{2} \right) & e^{3}\left(a^{3}b^{2}d-a^{2}b^{3}c\right)\\
0 & 1 & e^{3}\left(ac^{2}d^{3}-bc^{3}d^{2}\right) & e^{3}\left(b^{2}c^{3}d-a^{2}cd^{3}\right)
\end{array} \right)\\
&=&
\left( \begin{array}{cccc}
1 & 0 & -e^{2}ab(ad+bc) & e^{2}a^{2}b^{2}\\
0 & 1 & e^{2}c^{2}d^{2} & -e^{2}cd(ad+bc)
\end{array} \right)
\end{eqnarray*}
Dies liefert schlie\ss lich die folgende Darstellung auf $\tilde{U}$:
\[
\sigma \mapsto
\left( \begin{array}{ccccc}
e^{2}a^{3}d&e^{2}a^{3}c&e^{2}b^{3}d&e^{2}b^{3}c&-e^{2}ab(ad+bc)\\
e^{2}a^{3}b&e^{2}a^{4}&e^{2}b^{4}&e^{2}ab^{3}&e^{2}a^{2}b^{2}\\
e^{2}c^{3}d&e^{2}c^{4}&e^{2}d^{4}&e^{2}cd^{3}&e^{2}c^{2}d^{2}\\
e^{2}bc^{3}&e^{2}ac^{3}&e^{2}bd^{3}&e^{2}ad^{3}&-e^{2}cd(ad+bc)\\
0&0&0&0&1
\end{array} \right).
\]
Damit hat $\tilde{U}^{*}$ die Darstellung (transponieren und auswerten an der Stelle $\sigma^{-1}=\left( \begin{array}{cc} ed&-eb\\ -ec&ea \end{array} \right)_{ad-bc}$) 
\[
\sigma \mapsto
\left( \begin{array}{ccccc}
e^{2}ad^{3} & -e^{2}bd^{3} & -e^{2}ac^{3} & e^{2}bc^{3} & 0\\
-e^{2}cd^{3} & e^{2}d^{4} & e^{2}c^{4} & -e^{2}c^{3}d & 0\\
-e^{2}ab^{3} & e^{2}b^{4} & e^{2}a^{4} & -e^{2}a^{3}b & 0\\
e^{2}b^{3}c & -e^{2}b^{3}d & -e^{2}a^{3}c & e^{2}a^{3}d &0\\
e^{2}bd(ad+bc) & e^{2}b^{2}d^{2} & e^{2}a^{2}c^{2} & e^{2}ac(ad+bc) & 1
\end{array} \right).
\]
F"ur \SL2 haben wir $\tilde{U}^{*}$ in zweiter Potenz von $\langle X^{2},Y^{2},XY \rangle$ wiederge\-fun\-den. Da hier nun der Faktor $e^{2}$ auftaucht, versuchen wir die zweite Potenz von $M:=\langle X^{2},Y^{2},XY \rangle \otimes \langle E \rangle$, und zwar berechnen wie die Darstellung von $S^{2}(M)$  bez"uglich derselben Basis wie f"ur \SL2 (nur tensoriert mit $E$).
Wir inspizieren die Berechnung der Darstellung von $S^{2}(M)$ auf Seite \pageref{S2M}. Bei der Berechnung der Operation auf den ersten vier Basisvektoren haben wir die Eigenschaft $ad-bc=1$ gar nicht benutzt - die ersten vier Spalten der Darstellungsmatrix (\ref{DarstS2M}) k"onnen wir also einfach mit $e^{2}$ multipliziert "uber\-nehmen. Bei der Berechnung von $\sigma \cdot (X^{2}Y^{2}-(XY)^{2})$ tritt der Faktor $(ad-bc)^{2}$ auf. Nach Multiplikation mit $e^{2}$ wird dieser wieder zu $1$. Die Darstellung auf $S^{2}(M)$ ist also gegeben durch
\[
\sigma \mapsto
\left( \begin{array}{cccccc}
e^{2}ad^{3} & -e^{2}bd^{3} & -e^{2}ac^{3} & e^{2}bc^{3} & 0&*\\
-e^{2}cd^{3} & e^{2}d^{4} & e^{2}c^{4} & -e^{2}c^{3}d & 0&*\\
-e^{2}ab^{3} & e^{2}b^{4} & e^{2}a^{4} & -e^{2}a^{3}b & 0&*\\
e^{2}b^{3}c & -e^{2}b^{3}d & -e^{2}a^{3}c & e^{2}a^{3}d &0&*\\
e^{2}bd(ad+bc) & e^{2}b^{2}d^{2} & e^{2}a^{2}c^{2} & e^{2}ac(ad+bc) & 1 &*\\
0&0&0&0&0&*
\end{array} \right),
\]
und wir erkennen $\tilde{U}^{*}$ als Untermodul von $S^{2}(M)$. Damit erhalten wir das Hauptresultat f"ur \GL2 in Charakteristik $3$:

\begin{Bsp}
F"ur $p=3$ ist mit
\[
V^{*}:=\left( \langle X,Y \rangle \otimes\langle E^{2}\rangle \right) \oplus  \langle X^{3},Y^{3}\rangle \bigoplus_{i=1}^{3} \left( \langle X^{2},Y^{2},XY\rangle\otimes \langle E\rangle \right)
\]
der Invariantenring \INVGL2 nicht Cohen-Macaulay, $\dim V=13$. Der Kozyklus liegt im Grad $2$ im zum Tensorprodukt isomorphen Untermodul von $S^{2}\left(\left( \langle X,Y \rangle \otimes\langle E^{2}\rangle \right) \oplus  \langle X^{3},Y^{3}\rangle\right)$. Ein annullierendes phsop im Grad $2$ ist gegeben durch die drei Kopien der $(X^{2}\otimes E)(Y^{2}\otimes E)-(XY\otimes E)^{2}$.  \qed
\end{Bsp}

\noindent Von \emph{IsNotCohenMacaulay} (Datei {\tt Beispiel-6.14.txt}, {\tt dmax:=2+2+2+2=8, maxdp:=2}) wird dann im Grad $8$ eine Basis aus $220$ Invarianten gesch"atzt, w"ahrend es nur $210$ sind. An Rechenzeit sollte man etwa eine halbe Stunde einplanen. "Uberraschenderweise braucht das Beispiel f"ur \SL2 deutlich l"anger (und auch mehr Speicher).
\newpage

\subsection{Beispiele f"ur \SO in Charakteristik 2 und endliche Gruppen}
Wir besprechen zum Ende noch einige Beispiele f"ur nicht so interessante, weil nicht zusammenh"angende oder endliche Gruppen. Zun"achst zur Gruppe \SO. Man kann zeigen, dass diese f"ur $p=3$ linear reduktiv ist, dass es dort also nach Hochster und Roberts keine Beispiele geben kann. Wir betrachten daher nur den Fall $p=2$. Zun"achst ist zu beachten, dass die \SO f"ur diese Charakteristik anders definiert ist, n"amlich als die Unter\-grup\-pe der Elemente aus \SL2, f"ur die das Element $XY$ von $\langle X^{2},Y^{2},XY \rangle$ invariant ist, also $\sigma \cdot XY=XY$. Aus Tabelle \ref{Darstellungen} (S. \pageref{Darstellungen}) folgen damit die Bedingungen $ab=cd=0$. Damit ist \SO also eine lineare algebraische Gruppe mit den drei Relationen $ad-bc-1=0, ab=0, cd=0$. Falls $a=0$, folgt wegen $ad-bc=1$ jedenfalls $bc=1$, und damit $d=0$, also  $\sigma =\left( \begin{array}{cc} 0 & b\\b^{-1} & 0 \end{array}\right)$. Analog folgt im Falle $b=0$ sofort $ad=1$ und $c=0$, also $\sigma =\left( \begin{array}{cc} a & 0\\0 & a^{-1} \end{array}\right)$. Folglich ist die \SO semidirektes Produkt des Torus mit der Gruppe $Z_{2}$. Genauer ist durch
\[
K^{*} \rtimes Z_{2} \rightarrow \textrm{SO}_{2}(K): (a,0) \mapsto \left( \begin{array}{cc} a & 0\\0 & a^{-1} \end{array}\right), \quad (a,1) \mapsto \left( \begin{array}{cc} 0 & a\\a^{-1} & 0 \end{array}\right)
\]
ein Gruppenisomorphismus gegeben, wobei
\[
\begin{array}{rl}
(a,0)(b,0)=(ab,0) & (a,1)(b,0)=(ab^{-1},1)\\
(a,0)(b,1)=(ab,1) &  (a,1)(b,1)=(ab^{-1},0).
\end{array}
\]
Ausserdem ist durch
\[
\textrm{SO}_{2}(K) \rightarrow Z_{2}: \left( \begin{array}{cc} a & b\\c & d \end{array}\right)  \mapsto bc \quad (\in \{0,1\})
\]
bzw.
\[
K^{*} \rtimes Z_{2} \rightarrow Z_{2}: (a,x) \mapsto x
\]
ein Gruppenhomomorphismus gegeben.

Nun ist aber $K^{*}$ linear reduktiv, d.h. wenn wir ein nicht \CM Beispiel f"ur die Gruppe \SO angeben, kommt das im wesentlichen von der (nicht interessanten, da endlichen) Gruppe $Z_{2}$, und f"ur diese sind bereits seit geraumer Zeit Beispiele bekannt. Dennoch wird sich hier ein interessanter Effekt im Zu\-sam\-men\-hang mit dem Hauptsatz zeigen.

\subsubsection{Beispiele f"ur \SO und $Z_{2}$}
Wir betrachten die Darstellung des Tensorprodukt $\langle X\otimes X,\ldots,Y\otimes Y \rangle$ in Gleichung (\ref{TensBasis}) auf Seite \pageref{TensBasis} und schr"anken auf \SO ein. Da f"ur $\sigma=\left( \begin{array}{cc} a & b\\c & d \end{array}\right)$ aus \SO in jedem Fall $ab=cd=ac=bd=0$ gilt, erhalten wir f"ur die Einschr"ankung die Darstellung
\begin{equation} \label{TensSO}
\sigma \mapsto \left( \begin{array}{cccc}
a^{2} & b^{2} & 0 & 0\\
c^{2} & d^{2} & 0 & 0\\
0 & 0 & 1 & bc\\
0 & 0 & 0 & 1\\
\end{array} \right),
\end{equation}
und erkennen in der rechten unteren Teilmatrix die Darstellung eines Unter\-moduls, von nun an bezeichnet mit
\[
 M_{1bc}:=\langle \pi, \tau \rangle \textrm{ und Darstellung }\sigma \mapsto \left( \begin{array}{cc} 1 & bc\\0 & 1 \end{array}\right).
\]
Der triviale Modul $I:=\langle \pi \rangle$ hat also einen nichttrivialen Kozyklus!\\

{\small \noindent
Da die linke obere $3 \times 3$ Teilmatrix von (\ref{TensSO}) die Einschr"ankung der Darstellung von $\langle \mu,\nu,\pi \rangle$ auf \SO beschreibt, erkennen wir, dass auch dieser Modul einen nichttrivialen Kozyklus enth"alt. Damit k"onnen wir schon mal alle Beispiele aus Abschnitt \ref{munupi} (nach Einschr"ankung auf \SO) "ubernehmen, was aber nicht sonderlich interessant ist.

Ausserdem sehen wir mittels des angegebenen Gruppen\-homomorphismus auf die $Z_{2}$, dass es sich hier im Wesentlichen um eine Darstellung dieser Gruppe handelt.\\}

\noindent Da offenbar $\tilde{I}=M_{1bc}$  selbstdual ist (Darstellung des Duals ist $\sigma \mapsto \left( \begin{array}{cc} 1 & 0\\bc & 1 \end{array}\right)$), erhalten wir nach Korollar \ref{NCM} folgendes Beispiel:
\begin{Bsp}
F"ur \SO und $p=2$ ist mit
\[
V^{*}:=\langle \pi \rangle \oplus_{i=1}^{3} \langle \pi_{i}, \tau_{i} \rangle \quad(\cong V)
\]
der Invariantenring \INVSO2 nicht Cohen-Macaulay, $\dim V=7$. Der Kozyklus liegt im trivialen Modul im Grad $1$, und ein annullierendes phsop ist gegeben durch $\pi_{1},\pi_{2},\pi_{3}$. \qed
\end{Bsp}

Ferner ist $M_{1bc}$ als Untermodul des Tensorprodukts $\langle X\otimes X,\ldots,Y\otimes Y \rangle$ auch Untermodul von $S^{2}(\langle X,Y \rangle)$, und wir erhalten
\begin{Bsp}
F"ur \SO und $p=2$ ist mit
\[
V^{*}:=\langle \pi \rangle \oplus_{i=1}^{4} \langle X_{i}, Y_{i} \rangle \quad(\cong V)
\]
der Invariantenring \INVSO2 nicht Cohen-Macaulay, $\dim V=9$. Der Kozyklus liegt im trivialen Modul im Grad $1$, und ein annullierendes phsop im Grad $2$ ist gegeben durch $X_{i}Y_{j}+X_{j}Y_{i}$ mit $(i,j)=(1,2),(2,3),(3,4)$. \qed
\end{Bsp}

\noindent Beim Testen mit \emph{IsNotCohenMacaulay} (Datei {\tt Beispiel-6.15-18.txt}) er\-leben wir zum ersten Mal eine "Uberraschung: Gem"a\ss{} unserer Theorie sollte man jeweils {\tt dmax:=4} bzw. {\tt dmax:=7} setzen, aber {\tt IsNotCohenMacaulay} ver\-sch"atzt sich jeweils schon in den Graden $3$ bzw. $6$! (Derselbe Effekt tritt bei den Einschr"ankungen der Beispiele aus
Abschnitt \ref{munupi} auf \SO auf). Da man den trivialen Modul $\langle \pi \rangle$ gef"uhlsm"a\ss ig sowieso als "uberfl"ussig empfindet, lassen wir ihn weg und "uberpr"ufen das mit {\tt IsNotCohenMacaulay} - wir erhalten tats"achlich ein positives Ergebnis, und es wird sich auch hier in den Graden $3$ und $6$ versch"atzt. Damit erhalten wir die (einzigen!) auf rein experimentellem Weg gefundenen Beispiele - wir geben nachher noch die theoretische Be\-gr"un\-dung an.

\begin{Bsp}
F"ur \SO und $p=2$ ist mit
\[
V^{*}:= \oplus_{i=1}^{3} \langle \pi_{i}, \tau_{i} \rangle \quad(\cong V)
\]
der Invariantenring \INVSO2 nicht Cohen-Macaulay, $\dim V=6$. Es gibt keinen (konstruierten) Kozyklus mehr. Ein phsop, das keine regul"are Sequenz ist, ist gegeben durch $\pi_{1},\pi_{2},\pi_{3}$.
\end{Bsp}

\begin{Bsp}
F"ur \SO und $p=2$ ist mit
\[
V^{*}:=\oplus_{i=1}^{4} \langle X_{i}, Y_{i} \rangle \quad(\cong V)
\]
der Invariantenring \INVSO2 nicht Cohen-Macaulay, $\dim V=8$. Ein phsop im Grad $2$, das keine regul"are Sequenz ist, ist gegeben durch $X_{i}Y_{j}+X_{j}Y_{i}$ mit $(i,j)=(1,2),(2,3),(3,4)$.
\end{Bsp}

Wir m"ussen noch {\it beweisen}, warum die angegebenen phsops keine regu\-l"aren Sequenzen sind (wir machen dies O.E. f"ur das phsop im Grad 1. Sofern man \Magma traut, ist nat"urlich auch das Ergebnis von \emph{IsNotCohenMacaulay} ein Beweis). Unser Kozyklus in $I:=\langle v_{1}:=\pi_{0} \rangle$ ist von der Form $g_{\sigma} =(\sigma-1)\tau_{0}$, wobei $\tilde{I}=M_{1bc}=\langle v_{1}:=\pi_{0}, v_{2}:=\tau_{0} \rangle$. Dann ist $\tilde{I}^{*}=\langle v_{1}^{*}=\tau_{i}, v_{2}^{*}=\pi_{i} \rangle$ ($i=1,2,3$) , und $v_{2}^{*}=\pi_{i}$ ist die annullierende Invariante. 
Im Beweis von Proposition \ref{anni} (S. \pageref{anni}) haben wir genau berechnet, wie annulliert wird: Nach den Gleichungen (\ref{xanni}) und (\ref{piganni}) gilt n"amlich mit $x:=-v_{1}^{*} \cdot v_{1}=\tau_{i} \cdot \pi_{0}$, dass 
\begin{equation} \label{annispez}
\pi_{i} \cdot g_{\sigma}=(\sigma-1)x=(\sigma-1)(\tau_{i}\pi_{0}).
\end{equation}
Wir betrachten nun den Beweis des Hauptsatzes (S. \pageref{HSatz}), und setzen
\[
V_{1}^{*}:=\langle \pi_{0} \rangle \oplus_{i=1}^{3} \langle \tau_{i}, \pi_{i} \rangle
\]
und
\[
V_{2}^{*}:=\oplus_{i=1}^{3} \langle \tau_{i}, \pi_{i} \rangle \subseteq V_{1}^{*}.
\] 
Die $a_{i}$ sind die $\pi_{i}$, und in Gleichung (\ref{bi}) (S. \pageref{bi}) werden die $b_{i}$ nach Gleichung (\ref{annispez}) definiert als $b_{i}=\tau_{i}\pi_{0}$. Als n"achstes wird in (\ref{uij}) definiert
\[
u_{ij}:=\pi_{i}b_{j}-\pi_{j}b_{i}=\pi_{0}(\pi_{i}\tau_{j}-\pi_{j}\tau_{i})=:\pi_{0}\cdot u_{ij}'.
\]
Dabei sind die $u_{ij}$ in $K[V_{1}]^{G}$ und die $u_{ij}'$ in $K[V_{2}]^{G} \subseteq K[V_{1}]^{G}$
Als n"achstes wird dann gezeigt, dass
\[
u_{23}\pi_{1}-u_{13}\pi_{2}+u_{12}\pi_{3}=0
\]
gilt - eine Relation in $K[V_{1}]^{G}$. Nach vorheriger Gleichung gilt dann nach Division durch $\pi_{0}$ jedenfalls auch
\[
u_{23}'\pi_{1}-u_{13}'\pi_{2}+u_{12}'\pi_{3}=0,
\]
was eine Relation in $K[V_{2}]^{G}$ ist. W"are nun $(\pi_{1},\pi_{2},\pi_{3})$ regul"ar in $K[V_{2}]^{G}$, so g"abe es also $f_{1}',f_{2}' \in  K[V_{2}]^{G} \subseteq K[V_{1}]^{G}$ mit
\[
u_{12}'=f_{1}'\pi_{1}+f_{2}'\pi_{2}.
\]
Nach Multiplikation dieser Gleichung mit $\pi_{0}$ und setzen von $f_{i}:=\pi_{0}f_{i}' \in K[V_{1}]^{G}$ erhalten wir daraus
\[
u_{12}=f_{1}\pi_{1}+f_{2}\pi_{2}.
\]
Diese Annahme (siehe Gleichung (\ref{u12Ann})) wird jedoch dann im weiteren Beweis des Hauptsatzes zu einem Widerspruch gef"uhrt, indem gezeigt wird, dass dann der Kozyklus $g$ trivial w"are. Also war unsere Annahme, dass das phsop $\pi_{1},\pi_{2},\pi_{3}$ in $K[V_{2}]^{G}$ eine regul"are Sequenz ist falsch - also ist $K[V_{2}]^{G}$ nicht Cohen-Macaulay.
Da wir die nichttriviale Relation im Hauptsatz also durch $\pi_{0}$ k"urzen k"onnen - also eine Relation einen Grad tiefer bekommen, kann man f"ur \emph{IsNotCohenMacaulay} auch jeweils {\tt dmax} um einen Grad niedriger ansetzen. \qed

\subsubsection{Beispiele f"ur $G=({\mathbb F}_{q},+)$}
Wir verallgemeinern nun noch die Beispiele des letzten Abschnitts f"ur die Gruppe $Z_{2} \cong ({\mathbb F}_{2},+)$ auf die endliche, also reduktive Gruppe  $G=({\mathbb F}_{q},+)$. Dabei ist $q=p^{n}$ und $p=\textrm{char }K$. Man beachte, dass $({\mathbb F}_{q},+)$ als Nullstellen\-menge von $X^{q}-X$ "uber $K$ (algebraisch abgeschlossen) tats"achlich algebraisch ist. Wir definieren den Modul $M:=\langle \pi, \tau \rangle$ durch die Darstellung
\[
({\mathbb F}_{q},+) \rightarrow K^{2 \times 2}: \sigma \mapsto \left( \begin{array}{cc} 1 & \sigma\\ 0 & 1 \end{array} \right),
\]
und erkennen, dass der triviale Modul $I:=\langle \pi \rangle$ einen nichttrivialen Kozyklus besitzt, und $\tilde{I}=M$. Die Darstellung von $M^{*}$ ist gegeben durch $\sigma \mapsto \left( \begin{array}{cc} 1 & 0\\ -\sigma & 1 \end{array} \right)$, und durch Multiplizieren der ersten Zeile/Spalte mit $-1$ und danach ver\-tau\-schen der beiden Zeilen/Spalten erkennen wir $M$ als selbstdual. Damit ist $\pi \in M=\tilde{I}^{*}$ die annullierende Invariante des Kozyklus in $I$. Mit (fast) w"ortlich dem selben Beweis wie zuvor (man muss in ungerader Charakteristik einige Minuszeichen einf"ugen, die jedoch insgesamt nur einen unwesentlichen Pro\-por\-tion\-alit"atsfaktor ergeben), k"onnen wir den Modul mit dem Kozyklus in Korollar \ref{NCM} weglassen, und wir erhalten

\begin{Bsp}
F"ur $G=({\mathbb F}_{q},+)$ mit $q=p^{n}, p=\textrm{char }K$ ist mit
\[
V^{*}:= \oplus_{i=1}^{3} \langle \pi_{i}, \tau_{i} \rangle \quad(\cong V)
\]
der Invariantenring $K[V]^{G}$ nicht Cohen-Macaulay, $\dim V=6$. Es gibt keinen (konstruierten) Kozyklus mehr. Ein phsop, das keine regul"are Sequenz ist, ist gegeben durch $\pi_{1},\pi_{2},\pi_{3}$. \qed
\end{Bsp}

\noindent Beim Test mit \emph{IsNotCohenMacaulay} (Datei {\tt Beispiel-6.19.txt}) gen"ugt es wieder, {\tt dmax:=3} zu setzen. F"ur {\tt q} kann man eine beliebige Primzahlpotenz w"ahlen.\\

\noindent Man beachte, dass man dieses Beispiel \emph{nicht} so ohne weiteres auf $G_{a}=(K,+)$ verallgemeinern kann, denn diese Gruppe ist nicht reduktiv!

\newpage
\appendix
\section{\mbox{Zusammenfassung der wichtigsten Beispiele}}
Wie stellen in diesem Anhang die wichtigsten Beispiele tabellarisch zusammen und skizieren teilweise auch den Weg, wie wir sie erhalten haben. Zur schneller\-en "Ubersicht geben wir auch nochmal kurz die verwendeten Notationen und Bezeichnungen an.

\subsection{Notation}
Wir bezeichnen mit $K$ stets einen algebraisch abgeschlossenen K"orper der Charakteristik $p>0$. Ein Element der Gruppe \GL2 bezeichnen wir stets mit 
$\sigma = \left( \begin{array}{cc}
a & b\\
c & d
\end{array} \right)$. Wollen wir betonen, dass es sich hierbei um ein Element aus einer linearen algebraischen Gruppe handelt, die Nullstellenmenge von $(ad-bc)e-1=0$ ist (also eine Teilmenge von $K^{5}$), so  schreiben wir auch ausf"uhrlicher
$\sigma = \left( \begin{array}{cc}
a & b\\
c & d
\end{array} \right)_{e}$ - die Inverse der Determinante $e$ geh"ort zum Gruppenelement dazu.
Im Folgenden bezeichnen wir mit $G$ stets eine der klassischen Untergruppen von \GL2.\\ 
Mit $\langle X,Y\rangle $ bezeichnen wir dann den Modul mit der nat"urlichen  Darstellung der Gruppe $G$, also
\[
\begin{array}{l}
\sigma \cdot X = aX+cY\\
\sigma \cdot Y = bX+dY.
\end{array}
\] 
Seine zweite symmetrische Potenz bezeichnen wir mit 
\[\langle X^{2},Y^{2},XY\rangle :=S^{2}(\langle X,Y\rangle )\]
 und analog f"ur h"ohere Potenzen. Den Dual von $\langle X^{2},Y^{2},XY\rangle $ f"ur $G=\textrm{SL}_{2}(K)$ in Charakteristik 2 bezeichnen wir mit 
\[
\langle \mu, \nu, \pi\rangle \quad :=\quad \langle X^{2},Y^{2},XY\rangle ^{*},
\]
 wobei $\mu:=(X^2)^{*},\nu:=(Y^2)^{*},\pi:=(XY)^{*}$. Er hat bzgl. dieser Basis die (Links-)Darstellung
\[
\left(
\begin{array}{ccc}
d^2 & c^2 & 0\\
b^2 & a^2 & 0\\
bd & ac & 1
\end{array} \right),
\]
d.h. $\pi$ ist eine Invariante.\\
Ist $V$ ein $G$-Modul, so bezeichnen wir den zugeh"origen Polynomring mit $K[V]:=S(V^{*})$, also die direkte Summe aller symmetrischen Potenzen des {\it Duals}. Daher geben wir f"ur die folgenden Beispiele stets den Dual $V^{*}$ an, so dass der Invariantenring $K[V]^{G}$ nicht Cohen-Macaulay ist. Gegebenenfalls Klammern wir die Basisvektoren aus $V^{*}$ um anzuzeigen, dass es sich dabei um Monome ersten Grades in $K[V]$ handelt. Beispielsweise sind f"ur\\
 $V^{*}:=\langle X^{2},Y^{2},XY\rangle $ durch $(X^{2})(Y^{2})$ und $(XY)^{2}$ zwei verschiedene Elemente in $S^2(V^{*})$ gegeben.

\subsection{Beispiele f"ur \SL2 in Charakteristik 2}
\begin{Bsp}[Das Beispiel aus dem Projekt bzw. aus Satz \ref{PKoz}]
\end{Bsp}
\begin{tabular}{|l|c|}
\hline
Modul & $V^{*}:=\langle X^2,Y^2\rangle  \bigoplus_{i=1}^{3} \langle \mu_{i}, \nu_{i}, \pi_{i}\rangle $\\
\hline
Dimension & $11$\\
\hline
phsop   & $\pi_{1},\pi_{2},\pi_{3}$\\
in Graden:  & $1,1,1$\\
\hline
Kozyklus im Summanden & $\langle X^2,Y^2\rangle $\\
im Grad & 1\\ 
\hline
\end{tabular}\\
\\
Der Kozyklus kann interpretiert werden als $g_{\sigma}:=(\sigma-1)(XY)$, wobei\\ $\langle X^{2},Y^{2}\rangle $ als Untermodul von $\langle X^{2},Y^{2},XY\rangle $ zu sehen ist. \\
Man fin\-det $\langle \mu,\nu,\pi\rangle $ als Untermodul von $\langle X\otimes X,...,Y\otimes Y\rangle $. \\
Dies f"uhrt zu\\

\begin{Bsp}[Der Weg zum Hauptbeispiel]
\end{Bsp} 
\begin{tabular}{|l|c|}
\hline
Modul & $V^{*}:=\langle X^2,Y^2\rangle  \bigoplus_{i=1}^{3} \langle X\otimes X,...,Y\otimes Y\rangle $\\
\hline
Dimension & $14$\\
\hline
phsop   & dreimal  $X\otimes Y + Y \otimes X$ \\
in Graden:  & $1,1,1$\\
\hline
Kozyklus im Summanden & $\langle X^2,Y^2\rangle $\\
im Grad & 1\\ 
\hline
\end{tabular}\\
\\
$V$ ist selbstdual.\\
Es ist  $\langle X\otimes X,...,Y\otimes Y\rangle $ ein Untermodul von $S^{2}(\langle X,Y\rangle  \oplus \langle X,Y\rangle )$. Damit gelangen wir zu 

\begin{samepage}
\begin{Bsp}[Das Hauptbeispiel in Dimension 10]
\end{Bsp} 
\begin{tabular}{|l|c|}
\hline
Modul & $V^{*}:=\langle X^2,Y^2\rangle  \bigoplus_{i=1}^{4} \langle X_{i},Y_{i}\rangle $\\
\hline
Dimension & $10$\\
\hline
phsop   & $X_{i}Y_{j}+X_{j}Y_{i}$ mit $(i,j) \in \left\{(1,2),(2,3),(3,4) \right\}$ \\
in Graden:  & $2,2,2$\\
\hline
Kozyklus im Summanden & $\langle X^2,Y^2\rangle $\\
im Grad & 1\\ 
\hline
\end{tabular}\\
\\
Man ben"otigt den Summanden $\langle X,Y\rangle $ viermal, weil man sonst kein phsop in $K[V]$ erhalten w"urde.\\ $V$ ist selbstdual (da alle Summanden es sind).
\end{samepage}

\begin{Bsp}[Das erste wesentlich neue Beispiel nach dem Projekt] \label{neu1}
\end{Bsp} 
\begin{tabular}{|l|c|}
\hline
Modul & $V^{*}:=\langle \mu,\nu,\pi\rangle  \bigoplus_{i=1}^{3} \langle X\otimes X,...,Y\otimes Y\rangle $\\
\hline
Dimension & $15$\\
\hline
phsop   & dreimal $X\otimes Y + Y \otimes X$ \\
in Graden:  & $1,1,1$\\
\hline
Kozyklus im Summanden & $\langle \mu,\nu,\pi\rangle $\\
im Grad & 1\\ 
\hline
\end{tabular}\\

\begin{Bsp}[... und die Ab"anderung] \label{neu2}
\end{Bsp} 
\begin{tabular}{|l|c|}
\hline
Modul & $V^{*}:=\langle \mu, \nu, \pi\rangle  \bigoplus_{i=1}^{4} \langle X_{i},Y_{i}\rangle $\\
\hline
Dimension & $11$\\
\hline
phsop   & $X_{i}Y_{j}+X_{j}Y_{i}$ mit $(i,j) \in \left\{(1,2),(2,3),(3,4) \right\}$ \\
in Graden:  & $2,2,2$\\
\hline
Kozyklus im Summanden & $\langle \mu, \nu, \pi\rangle $\\
im Grad & 1\\ 
\hline
\end{tabular}\\
\\
\\
Die Summanden mit den Invarianten kann man nat"urlich auch mischen, was aber mehr eine Spielerei ist. Damit bekommt man ein phosp in verschiedenen Graden. Wir geben nur ein
\begin{samepage}
\begin{Bsp}[langweilig]
\end{Bsp} 
\begin{tabular}{|l|c|}
\hline
Modul & $V^{*}:=\langle X^2,Y^2\rangle  \bigoplus_{i=1}^{3} \langle X_{i},Y_{i}\rangle  \bigoplus \langle X \otimes X,...,Y \otimes Y\rangle $\\
\hline
Dimension & $12$\\
\hline
phsop   & $X_{1}Y_{2}+X_{2}Y_{1},X_{2}Y_{3}+X_{3}Y_{2} , X\otimes Y+ Y\otimes X$ \\
in Graden:  & $2,2,1$\\
\hline
Kozyklus im Summanden & $\langle X^2,Y^2\rangle $\\
im Grad & 1\\ 
\hline
\end{tabular}\\
\end{samepage}
\\
Die Invariante $\pi^{2}$ von $S^{2}(\langle \mu,\nu,\pi\rangle )=:\tilde{U}^{*}$ annulliert den Kozyklus von $U$, und es stellt sich heraus, dass $U \cong \langle \mu\pi,\nu\pi,\mu^{2},\nu^{2},\pi^{2}\rangle $. Ferner ist 
\[
\langle \mu,\nu,\pi\rangle  \cong \langle \mu\pi,\nu\pi,\pi^{2}\rangle  \quad \le \quad \langle \mu\pi,\nu\pi,\mu^{2},\nu^{2},\pi^{2}\rangle .
\]
Man kann sich also einen Summanden sparen, weil $\langle \mu\pi,\nu\pi,\mu^{2},\nu^{2},\pi^{2}\rangle $ sowohl den Kozyklus als auch gleich in zweiter Potenz eine annullierende Invariante enth"alt. Dies liefert

\newpage
\begin{Bsp}[Kozyklus und Invariante in einem]
\end{Bsp} 
\begin{tabular}{|l|c|}
\hline
Modul & $V^{*}:=\langle \mu_{1}\pi_{1},\nu_{1}\pi_{1},\mu_{1}^{2},\nu_{1}^{2},\pi_{1}^{2}\rangle    \bigoplus_{i=2}^{3} \langle \mu_{i},\nu_{i},\pi_{i}\rangle $\\
\hline
Dimension & $11$\\
\hline
phsop   & $(\pi_{1}^{2})^{2},\pi_{2}^{2},\pi_{3}^{2} $ \\
in Graden:  & $2,2,2$\\
\hline
Kozyklus im Summanden & $\langle \mu_{1}\pi_{1},\nu_{1}\pi_{1},\mu_{1}^{2},\nu_{1}^{2},\pi_{1}^{2}\rangle $\\
im Grad & 1\\ 
\hline
\end{tabular}\\
\\
\\
Falls der Summand mit dem Kozyklus eine nicht zur Annullation ben"otigte Invariante besitzt, und es eine annullierende Invariante in erster Potenz gibt, kann man diese Verheften: Sind
\[
\left(
\begin{array}{cc}
A_{\sigma} & 0\\
g_{\sigma} & 1
\end{array} \right), \quad \quad
\left(
\begin{array}{cc}
B_{\sigma} & 0\\
h_{\sigma} & 1
\end{array} \right)
\]
Darstellungen von Moduln mit einer Invarianten, so ist
\[
\left(
\begin{array}{ccc}
A_{\sigma} & 0 &0\\
0 & B_{\sigma} &0\\
g_{\sigma} & h_{\sigma} & 1
\end{array} \right)
\]
Darstellung eines Moduls, der die beiden anderen als Untermodul enth"alt, aber eine Dimension kleiner ist als die direkte Summe.\\
Ausgehend von dem Modul
\[
V^{*}:=\langle \mu,\nu,\pi\rangle  \oplus \langle X\otimes X,...,Y\otimes Y\rangle  \bigoplus_{i=1}^{3} \langle X_{i},Y_{i}\rangle 
\]
mit nicht \CM Invariantenring bezeichnen wir mit 
\[
U \textrm{ die Verheftung von } \langle \mu,\nu,\pi\rangle  \textrm{ und } \langle X\otimes X,...,Y\otimes Y\rangle  
\]
an den Invarianten $\pi$ bzw. $X \otimes Y + Y \otimes X$, wobei wir die Invariante von $U$ wieder mit $\pi$ bezeichnen. Damit haben wir

\begin{Bsp}[mit verhefteten Invarianten]
\end{Bsp} 
\begin{tabular}{|l|c|}
\hline
Modul & $V^{*}:=U \bigoplus_{i=1}^{3} \langle X_{i},Y_{i}\rangle $\\
\hline
Dimension & $12$\\
\hline
phsop   & $\pi,X_{1}Y_{2}+X_{2}Y_{1},X_{2}Y_{3}+X_{3}Y_{2}$ \\
in Graden:  & $1,2,2$\\
\hline
Kozyklus im Summanden & $U$\\
im Grad & 1\\ 
\hline
\end{tabular}\\

\subsection{Beispiele f"ur \SL2 in Charakteristik 3}
Wir erinnern zun"achst nochmal an die Konstruktion aus Satz \ref{PKoz}: Sei 
\[
M:=\langle X^{3},Y^{3},X^{2}Y,XY^{2}\rangle  \textrm{ und } W:=\langle X^{3},Y^{3}\rangle  \quad \le \quad  M
\]
sowie
\[
\begin{array}{l}
U:=\left\{ f \in \textrm{Hom}_{K}(M,W) : f|_{W}=0 \right\}\\
\iota \in \textrm{Hom}_{K}(M,W) \textrm{ mit } \iota|_{W}=\textrm{id}_{W} \textrm{ und } \iota(X^{2}Y)=\iota(XY^{2})=0\\
\tilde{U}:=U \oplus K \cdot \iota.
\end{array}
\]
Dies ergab dann folgendes

\begin{Bsp}[Das Beispiel aus dem Projekt bzw. nach Satz \ref{PKoz}]
\end{Bsp}
\begin{tabular}{|l|c|}
\hline
Modul & $V^{*}:=U \bigoplus_{i=1}^{3} \tilde{U}^{*}$\\
\hline
Dimension & $19$\\
\hline
phsop   & $\pi_{1},\pi_{2},\pi_{3}$\\
in Graden:  & $1,1,1$\\
\hline
Kozyklus im Summanden & $U$\\
im Grad & 1\\ 
\hline
\end{tabular}\\
\\
\\
Nun findet man aber, wenn man obige Konstruktion von $U$ f"ur einen beliebigen Modul $M$ mit Untermodul $W$ durchf"uhrt, dass
\[
U \cong W \otimes (M/W)^{*}
\]
gilt. Dann ist jedoch
\[
U \quad \le \quad W \otimes (M/W)^{*} \bigoplus S^{2}(W) \bigoplus S^{2}((M/W)^{*}) \quad \cong \quad S^{2}(W \oplus (M/W)^{*}).
\]
Da die Summe direkt ist, bleibt der Kozyklus nichttrivial, also kann man bei der Konstruktion f"ur $V^{*}$ schon mal $U$ durch die direkte Summe $W \oplus (M/W)^{*}$ ersetzen, welche dann einen Kozyklus in zweiter Potenz enth"alt. Diese Konstruktion l"asst sich bei allen nach Satz \ref{PKoz} konstruierten Beispielen f"ur die Gruppe \SL2 durchf"uhren (f"ur beliebiges $\textrm{char }K=p$). Dies ergibt eine Di\-men\-sionsre\-duk\-tion von $8p-5$ auf $7p-2$. Im Falle $p=2$ ist dann aber der zweite Summand $(M/W)^{*}$ trivial (gibt also bei Tensorierung nichts neues), und im Falle $p=3$, dem wir uns nun wieder zuwenden, ist die Dimension immer noch $7 \cdot 3-2=19$. Allerdings findet man dort
\[
\begin{array}{lrcl}
&\tilde{U}^{*} & \le & S^{2}(\langle X^{2},Y^{2},XY\rangle )\\
\textrm{sowie} &  (M/W)^{*} &\cong & \langle X,Y\rangle .
\end{array}
\]
Dies alles ergibt

\begin{Bsp}[Hauptbeispiel mit Kozyklus und Invarianten im Grad $2$]
\end{Bsp}
\begin{tabular}{|l|c|}
\hline
Modul & $V^{*}:=\langle X,Y\rangle  \bigoplus \langle X^{3},Y^{3}\rangle  \bigoplus_{i=1}^{3} \langle X_{i}^{2},Y_{i}^{2},X_{i}Y_{i}\rangle $\\
\hline
Dimension & $13$\\
\hline
phsop   & $X_{i}^{2}Y_{i}^{2}-(X_{i}Y_{i})^{2}$\\
in Graden:  & $2,2,2$\\
\hline
Kozyklus in Potenz & $S^{2}(\langle X,Y\rangle  \oplus \langle X^{3},Y^{3}\rangle )$\\
im Grad & 2\\ 
\hline
\end{tabular}\\
$V$ ist selbstdual.\\
\\
\subsection{Beispiele f"ur \GL2 mit $p=2,3$}
Satz \ref{PKoz} liefert auch f"ur \GL2 Moduln $U$ mit nichttrivialem Kozyklus. Allerdings ist nicht ganz offensichtlich, wie man die Beispiele f"ur \SL2 ver"andern muss, um wieder $\tilde{U}^{*}$ in einer zweiten Potenz wiederzufinden. Im wesentlichen geht das durch Tensorieren mit Potenzen der Inversen der Determinante, aber die Frage ist, wo und mit welcher Potenz. Dazu etwas Notation: Wir bezeichnen mit $\langle E\rangle $ den Modul, auf dem \GL2 mit dem Inversen der Determinante operiert, also
\[
\sigma \cdot E =
\left( \begin{array}{cc}
a & b\\
c & d
\end{array} \right)_{e} \cdot E=
\frac{1}{ad-bc}\cdot E=e\cdot E,
\]
und wir sehen, dass die Operation tats"achlich algebraisch ist.\\
Mit der leicht zu pr"ufenden Isomorphie
\[
\langle X,Y\rangle  \otimes \langle E\rangle  \quad \cong \quad \langle X,Y\rangle^{*} \quad =: \quad \langle X^{*},Y^{*}\rangle 
\]
(Dabei ist $X^{*} \rightarrow Y\otimes E, Y^{*} \rightarrow -X\otimes E$)
kann man dann formulieren

\begin{Bsp}[Hauptbeispiel f"ur \GL2, $\textrm{char }K=2$]
\end{Bsp}
\begin{tabular}{|l|c|}
\hline
Modul & $V^{*}:=\left( \langle X^2,Y^2\rangle \otimes\langle E\rangle \right) \bigoplus_{i=1}^{2}( \langle X_{i},Y_{i}\rangle  \oplus \langle X_{i}^{*},Y_{i}^{*}\rangle )$\\
\hline
Dimension & $10$\\
\hline
phsop   & $X_{i}X_{j}^{*}+Y_{i}Y_{j}^{*}$ mit $(i,j) \in \left\{(1,1),(1,2),(2,2) \right\}$ \\
in Graden:  & $2,2,2$\\
\hline
Kozyklus im Summanden & $\langle X^2,Y^2\rangle \otimes \langle E\rangle $\\
im Grad & 1\\ 
\hline
\end{tabular}\\

Anstatt mit $\langle E\rangle $ zu tensorieren, kann man auch einfach $\langle E\rangle $ hinreichend oft addieren - dann liegen der Kozyklus bzw. das phsop noch einen Grad h"oher. Allerdings bringt die Addition nat"urlich eine unerw"unschte Di\-men\-sions\-erh"ohung mit sich. Insbesondere muss man $\langle E\rangle $ so oft addieren, dass sich noch ein phsop ergibt. Bei dem Kozyklus ist dies dagegen nicht kritisch, es kann hier durchaus ein Summand $\langle E\rangle $ verwendet werden, der auch schon f"ur eine Invariante verwendet wurde. Man kann also z.B.
\[
V^{*}:=\langle X^2,Y^2\rangle  \bigoplus_{i=1}^{4} \langle X_{i},Y_{i}\rangle  \bigoplus_{i=1}^{3}\langle E_{i}\rangle 
\] setzen. Jedes Element im phsop bekommt hier ein $E_{i}$ (deshalb braucht man drei Summanden), und ein $E_{i}$ l"asst sich f"ur den Kozyklus im Grad 2 recyclen. Also\\
\\

\begin{Bsp}[Aufgebl"ahtes Beispiel f"ur \GL2, $\textrm{char }K=2$]
\end{Bsp}
\begin{tabular}{|l|c|}
\hline
Modul & $V^{*}:=\langle X^2,Y^2\rangle  \bigoplus_{i=1}^{4} \langle X_{i},Y_{i}\rangle  \bigoplus_{i=1}^{3}\langle E_{i}\rangle 
$\\
\hline
Dimension & $13$\\
\hline
phsop   & $X_{i}Y_{j}E_{k}+X_{j}Y_{i}E_{k}$ mit $(i,j,k) \in \left\{(1,2,1),(2,3,2),(3,4,3) \right\}$ \\
in Graden:  & $3,3,3$\\
\hline
Kozyklus in & $S^{2}(\langle X^2,Y^2\rangle \oplus  \langle E_{1}\rangle )$\\
im Grad & 2\\ 
\hline
\end{tabular}\\

"Ahnlich gehen wir im Fall $p=3$ vor:
\begin{Bsp}[Hauptbeispiel f"ur \GL2, $\textrm{char }K=3$]
\end{Bsp}
\begin{tabular}{|l|c|}
\hline
Modul & $V^{*}:=\left( \langle X,Y\rangle \otimes\langle E^{2}\rangle  \right) \bigoplus \langle X^{3},Y^{3}\rangle  \bigoplus_{i=1}^{3} \langle X_{i}^{2},Y_{i}^{2},X_{i}Y_{i}\rangle \otimes\langle E_{i}\rangle$\\
\hline
Dimension & $13$\\
\hline
phsop   & $ (X_{i}^{2}\otimes E_{i})(Y_{i}^{2}\otimes E_{i})-(X_{i}Y_{i}\otimes E_{i})^{2}  $ \\
in Graden:  & $2,2,2$\\
\hline
Kozyklus in  & $\ S^{2}(\langle X,Y\rangle \otimes \langle E^{2}\rangle  \bigoplus \langle X^{3},Y^{3}\rangle )$\\
im Grad & 2\\ 
\hline
\end{tabular}\\

\subsection{Beispiele f"ur \SO, $\textrm{char }K=2$}
Die Elemente von \SO werden durch die Relation $\sigma \cdot XY=XY$ aus 
\SL2 ausgesondert. Es folgt 
$\sigma = \left( \begin{array}{cc}
a & 0\\
0 & a^{-1}
\end{array} \right)$
oder 
$\sigma = \left( \begin{array}{cc}
0 & b\\
b^{-1} & 0
\end{array} \right)$. Damit hat diese Gruppe jedoch als Zu\-sammen\-hangs\-kom\-ponen\-te den Torus und ist daher nicht so interessant. Wir geben dennoch einige Beispiele an. Wir k"onnen Beispiel \ref{neu1} und \ref{neu2} "ubernehmen. Man findet nun im Tensorprodukt \\
$\langle X\otimes X,...,Y\otimes Y\rangle $ einen Untermodul, bezeichnet mit 
\[
M_{1bc}=:\langle \pi,\tau\rangle  \textrm{ und der Darstellung }
\left( \begin{array}{cc}
1 & bc\\
0 & 1
\end{array} \right).\]
Damit hat der triviale Modul $\langle \pi\rangle $ schon einen nicht trivialen Kozyklus, der von der Invariante von $M_{1bc}$ (selbstdual) annulliert wird, und wir haben
\begin{Bsp}[Ein Beispiel in Dimension 9]
\end{Bsp} 
\begin{tabular}{|l|c|}
\hline
Modul & $V^{*}:=\langle \pi\rangle  \bigoplus_{i=1}^{4} \langle X_{i},Y_{i}\rangle $\\
\hline
Dimension & $9$\\
\hline
phsop   & $X_{i}Y_{j}+X_{j}Y_{i}$ mit $(i,j) \in \left\{(1,2),(2,3),(3,4) \right\}$ \\
in Graden:  & $2,2,2$\\
\hline
Kozyklus im Summanden & $\langle \pi \rangle $\\
im Grad & 1\\ 
\hline
\end{tabular}\\
\begin{Bsp}[Ein Beispiel in Dimension 7]
\end{Bsp} 
\begin{tabular}{|l|c|}
\hline
Modul & $V^{*}:=\langle \pi\rangle  \bigoplus_{i=1}^{3} \langle \pi_{i},\tau_{i}\rangle $\\
\hline
Dimension & $7$\\
\hline
phsop   & $\pi_{1},\pi_{2},\pi_{3}$ \\
in Graden:  & $1,1,1$\\
\hline
Kozyklus im Summanden & $\langle \pi \rangle $\\
im Grad & 1\\ 
\hline
\end{tabular}\\

Nach dem Beweis zum Hauptsatz kann man aber in der auftretenden Relation f"ur das annullierende phsop die Koeffizienten durch $\pi$ k"urzen, so dass der triviale Modul f"ur das Erzeugen dieser Relation gar nicht gebraucht wird. Er kann also weggelassen werden und wir erhalten:

\begin{Bsp}[Ein Beispiel in Dimension 8]
\end{Bsp} 
\begin{tabular}{|l|c|}
\hline
Modul & $V^{*}:=\bigoplus_{i=1}^{4} \langle X_{i},Y_{i}\rangle $\\
\hline
Dimension & $8$\\
\hline
phsop   & $X_{i}Y_{j}+X_{j}Y_{i}$ mit $(i,j) \in \left\{(1,2),(2,3),(3,4) \right\}$ \\
in Graden:  & $2,2,2$\\
\hline
Kozyklus & -\\
im Grad & -\\ 
\hline
\end{tabular}\\
\newpage
\begin{Bsp}[Ein Beispiel in Dimension 6]
\end{Bsp} 
\begin{tabular}{|l|c|}
\hline
Modul & $V^{*}:=\bigoplus_{i=1}^{3} \langle \pi_{i},\tau_{i}\rangle $\\
\hline
Dimension & $6$\\
\hline
phsop   & $\pi_{1},\pi_{2},\pi_{3}$ \\
in Graden:  & $1,1,1$\\
\hline
Kozyklus & -\\
im Grad & -\\ 
\hline
\end{tabular}\\

Im Prinzip kommt hier aber die nicht \CM Eigenschaft von der Gruppe $Z_{2}$, denn
$\textrm{SO}_{2}(K) \cong K^{*} \rtimes Z_{2}$, und $K^{*}$ ist linear reduktiv.

\newpage
\section{Funktionsweise des (modifizierten) Bayer-Algorithmus}
Der Algorithmus von Bayer dient zur Berechnung einer $K$-Vektorraumbasis von $K[V]^{G}_{d}$ zu gegebenem Grad $d$. Dabei ist $V$ durch eine Darstellung (des Duals) gegeben und die Gruppe $G$ durch ihr zugeh"origes Ideal.
In diesem Anhang wollen wir nur kurz das Prinzip beschreiben, nach dem der Bayer-Algorithmus und seine Modifikation funktioniert. Die Idee f"ur die Modifikation geht auf Kemper zur"uck, der diese zun"achst f"ur seinen Algorithmus \cite{Kem1} vor\-schlug (der denselben Zweck wie der Bayer-Algorithmus erf"ullt). Da sich jedoch zumindest f"ur die in dieser Arbeit untersuchten Moduln der Bayer-Algorithmus stets als um einiges schneller erwies, geben wir die Modifikation nur f"ur diesen an. F"ur Einzelheiten und Be\-weise sei auf die Arbeit \cite{Bayer} verwiesen (die dortigen Beweise gehen auch f"ur die Modifikation durch). Ausserdem sei nochmals betont, dass ohne die Modifikation die Mehrheit der Beispiele aufgrund einer dann \emph{viel} zu langen Rechenzeit \emph{nicht} mit \emph{IsNotCohenMacaulay} h"atte getestet werden k"onnen.
\subsection{Berechnung von Torus-Invarianten}
Der modifizierte Bayer-Algorithmus berechnet zuerst die Invarianten einer sehr einfachen Untergruppe der \SL2, n"amlich des Torus
\[
T:=\left\{ 
\left( \begin{array}{cc}
a & 0\\
0 & a^{-1}
\end{array} \right): a\in K\setminus\{0\} \right\}.
\]
Wie bereits in Abschnitt \ref{ImplMagma} erw"ahnt, reduziert sich die Darstellungsmatrix eines $G$-Moduls (mit $T \subseteq G$) $V^{*}=\langle X_{1},\ldots,X_{n} \rangle$ f"ur die Elemente $\sigma \in T$  sehr oft auf die Form
\[
A_{\sigma}=\left( \begin{array}{cccc}
a^{w_{1}} & 0 & \cdots &0 \\
0 & a^{w_{2}} & &\vdots\\
\vdots & &\ddots &0\\
0 & \cdots &0 & a^{w_{n}}
\end{array} \right) \quad   \textrm{ d.h. } \sigma \cdot X_{i}=a^{w_{i}}X_{i},
\]
mit $w_{1},\ldots,w_{n} \in {\mathbb Z}$ fest, und kann ggf. mit Hilfe eines Basiswechsels stets auf diese Form gebracht werden. Der Gewichtsvektor $w=\left(w_{1},\ldots,w_{n}\right)$ beschreibt die Einschr"ankung der Darstellung auf $T$ vollst"andig.
Sei nun
\[
f=\sum f_{i_{1},\ldots,i_{n}}X^{i_{1}}\dots X^{i_{n}} \in K[V]^{G}_{d} \textrm{ oder } K[V]^{T}_{d},\quad \textrm{ mit } f_{i_{1},\ldots,i_{n}} \in K.
\]
F"ur $\sigma \in T \subseteq G$ gilt dann also $f=\sigma \cdot f$, oder
\[
\sum f_{i_{1},\ldots,i_{n}}X^{i_{1}}\dots X^{i_{n}}=\sum f_{i_{1},\ldots,i_{n}}a^{w_{1}i_{1}+\ldots+w_{n}i_{n}}       X^{i_{1}}\dots X^{i_{n}}.
\]
Es folgt $w_{1}i_{1}+\ldots+w_{n}i_{n}=0$ f"ur $f_{i_{1},\ldots,i_{n}} \ne 0$ (denn $K$ ist unendlich), oder anders ausgedr"uckt: Die Invarianten von $K[V]^{G}_{d}$ oder $K[V]^{T}_{d}$ sind eine Linearkombination \emph{monomialer} Invarianten von $K[V]^{T}_{d}$, n"amlich der $X^{i_{1}}\dots X^{i_{n}}$ mit $w_{1}i_{1}+\ldots+w_{n}i_{n}=0$. F"ur $K[V]^{T}_{d}$ ist eine Basis also gegeben durch
\[
\left\{
X^{i_{1}}\dots X^{i_{n}}: \quad i_{1}+\ldots+i_{n}=d,\quad w_{1}i_{1}+\ldots+w_{n}i_{n}=0
\right\}.
\]
Diese l"asst sich mit \Magma am effizientesten durch Auswahl aus der Menge aller Monome vom Grad $d$ bestimmen, etwa durch folgende Zeilen ({\tt P} ist dabei ein Modell f"ur $K[V]$, also ein Polynomring in $n$ Variablen):\\
{\tt mons:=MonomialsOfDegree(P,d);\\
 mons:=[f: f in mons| \&+[Exponents(f)[i]*w[i]: i in [1..n]] eq 0];}\\
(Siehe die Funktion {\tt TorusInvariantsSL2} in der Datei {\tt CMTest.txt}.)

\subsection{Das Prinzip des Bayer-Algorithmus}
Die Gruppe $G$ sei gegeben als Nullstellenmenge eines (Radikal-)Ideals $I_{G} \le K[S_{1},\ldots,S_{r}]$. Die Operation auf $V^{*}=\langle X_{1},\ldots,X_{n} \rangle$ ist f"ur $\sigma \in G \subseteq K^{r}$ gegeben durch eine Darstellungsmatrix $A_{\sigma}:=\left(a_{ij}(\sigma)\right)_{i,j=1..n}$ mit Polynomen $a_{ij} \in K[S_{1},...,S_{r}]$ mit  $i,j=1..n$. Die Polynome $a_{ij}$ seien dabei homogen vom gleichen Grad (dies ist keine Einschr"ankung, denn mit einer zus"atzlichen Variable $S_{r+1}$ und der Hinzunahme der Polynoms $S_{r+1}-1$ zur Definition von $I_{G}$ kann man die Polynome $a_{ij}$ \emph{homogenisieren}). Dann sind die Polynome
\[
\Psi_{i}:=\sum_{j=1}^{n}a_{ij}X_{j} \in K[S_{1},...,S_{r},X_{1},\ldots,X_{n}]
\]
ebenfalls alle homogen mit einem gemeinsamen Grad $\delta$. Ferner sei $I_{G}^{h}$ die Homogenisierung des Ideals $I_{G}$. Dies ist dasjenige Ideal in $K[S_{1},...,S_{r},h]$, das von den Homogenisierungen (zum kleinst m"oglichen Grad) jeden Elements aus $I_{G}$ mittels der neuen Variablen $h$ erzeugt wird. Sei nun
\[
M:=\left\{(i_{1},\ldots,i_{n}) \in {\mathbb N}_{0}^{n}: i_{1}+\ldots+i_{n}=d \right\}
\]
die Menge aller Exponenten von Monomen vom Grad $d$. Bayer zeigt nun in Proposition 1 in \cite{Bayer} folgendes Resultat, nach welchem klar ist, wie der Algorithmus zur Berechnung der Invarianten auszusehen hat:\\

\begin{samepage}
{\noindent \it Ist $\{f_{1},\ldots,f_{k}\}$ eine Gr"obner-Basis des Ideals $I \cap K[X_{1},\ldots,X_{n},h]$ mit
\[
I:=\left\langle \left\{\Psi_{1}^{i_{1}}\cdot \dots \cdot \Psi_{n}^{i_{n}}: (i_{1},\ldots,i_{n}) \in M 
 \right\} \cup I_{G}^{h} \right\rangle \le K[S,X,h],
\]
so ist
\[
\left\{f_{i}(X_{1},\ldots,X_{n},1): 1\le i \le k,\quad \deg f_{i}=d \cdot \delta 
 \right\}
\]
eine $K$-Vektorraumbasis von $K[X_{1},\ldots,X_{n}]^{G}_{d}$.}\\
\end{samepage}

Der dort gegebene Beweis funktioniert jedoch w"ortlich auch dann, wenn man $M$ durch die Menge
\[
M_{T}:=\left\{(i_{1},\ldots,i_{n}) \in {\mathbb N}_{0}^{n}: i_{1}+\ldots+i_{n}=d, \quad w_{1}i_{1}+\ldots+w_{n}i_{n}=0\right\}
\]
der Exponenten der Torus-Invarianten ersetzt. Bei der Beweis\-rich\-tung $\subseteq$ "andert sich gar nichts, und f"ur $\supseteq$ ist lediglich zu beachten, dass sich jede Invariante als Linearkombination der Torus-Invarianten berechnen l"asst. Bei der Implemen\-ta\-tion in Magma ersetzt man f"ur die Modifikation also lediglich die Zeile\\
{\tt mons:=MonomialsOfDegree(P,d);}\\
bei der die zu den Exponenten $M$ geh"origen Monome berechnet werden, durch die Zeile\\
{\tt mons:=TorusInvariantsSL2(w,d,P);}\\
 wo dann die zu $M_{T}$ geh"origen Monome berechnet werden. Wie bereits erw"ahnt, f"uhrt dies zu einer enormen Geschwindigkeitssteigerung.

\end{document}